\documentclass{amsart}

\usepackage{psfrag,  epsfig, color, graphicx}

\usepackage{amsmath, amssymb,  amsthm,amsfonts}

\newtheorem{theorem}{Theorem}[section]

\theoremstyle{definition}

\theoremstyle{remark}
\newtheorem{remark}[theorem]{Remark}

\numberwithin{equation}{section}

\swapnumbers \theoremstyle{plain}

\newtheorem{thm}{Theorem}[section]

\newtheorem{lem}[thm]{Lemma}

\theoremstyle{remark}

\theoremstyle{definition}
\newtheorem{defn}[thm]{Definition}

\newcommand{\T}{\mathcal T}

\newcommand{\bdy}{\partial}
\newcommand{\bbb}{\mathbb}
\newcommand{\rppp}{\mathbb{R}P^3}
\newcommand{\rpp}{\mathbb{R}P^2}

\newcommand{\td}{\tilde}
\newcommand{\open}[1]{\stackrel{\circ}{#1}}
\newcommand{\ve}{\varepsilon}
\newcommand{\lra}{\leftrightarrow}

\begin{document}

\title{Inflations of Ideal Triangulations}

\author{William Jaco}

\thanks{The first author was partially supported by NSF/DMS Grants,
The Grayce B. Kerr Foundation, The American Institute of Mathematics
(AIM), and and The Visiting Research Scholar Program at University
of Melbourne (Australia)}

\author{J. Hyam Rubinstein}
\thanks{The second author was partially supported by The Australian Research
Council and The Grayce B. Kerr Foundation.}

\begin{abstract}  

Starting with an ideal triangulation of $\open{M}$, the interior of a compact 3--manifold M with boundary, no component of which is a 2--sphere, we provide a construction, called an inflation of the ideal triangulation,  to obtain a strongly related triangulations of M itself.  Besides a step-by-step algorithm for such a construction, we provide examples of an inflation of the two-tetrahedra ideal triangulation of the complement of the {\it figure-eight} knot in $S^3$, giving a minimal triangulation, having ten tetrahedra, of the {\it figure-eight} knot exterior. As another example,  we provide an inflation of the one-tetrahedron Gieseking manifold giving a minimal triangulation, having seven tetrahedra, of a nonorientable compact 3-manifold with Klein bottle boundary.  Several applications of inflations are discussed.
\end{abstract}


\maketitle


\date{\today}

\section{Introduction} Triangulations play a central role in the the study and understanding of 3-manifolds. They are used directly or indirectly for the major work on a census of 3-manifolds \cite{bur-enumer, mat-comp-6, mat-tabulate, mar-pet-9} and are fundamental to most of our advances on decision problems, algorithms, and issues of computational complexity.  Triangulations naturally give rise to classes of surfaces called normal and almost normal surfaces, these surfaces in turn have been used in constructions of decompositions and recognition algorithms for 3--manifolds. A triangulation of a 3--manifold can be thought of as a combinatorial analog of a metric on the manifold and just as we try to deform metrics to gain geometric and topological information about a 3--manifold, we can similarly hope to gain geometric and topological information about a 3-manifold by deforming a given triangulation to a 'good' triangulation of the 3--manifold. This work contributes to constructions that can be used to modify one triangulation to another that exhibits desirable properties. 

It is well know that triangulations contain many normal surfaces that are not very interesting topologically but are artifacts of the triangulation; on the other hand, with certain modification, we can often arrive at a triangulation where there are useful connections between the geometry and topology of the manifold and the normal surfaces in the triangulation. In our work on 0--efficient  triangulations \cite{jac-rub-0eff}, the aim was to control normal surfaces with positive Euler characteristic; this leads to a very nice algorithm for the connected sum decomposition of a 3-manifold \cite{jac-rub-0eff} and triangulations that lend themselves nicely to the 3--sphere recognition algorithm \cite{rubin1, tho, jac-rub-0eff}.  For many algorithms and structure problems it is very desirable to control normal surfaces with zero Euler characteristic; our work on such triangulations includes the work presented here. One of its applications is control of normal annuli in 3--manifolds with boundary \cite{jac-rub-annular}.  Angle structures in ideal triangulations give interesting examples of connections between normal surfaces and the geometry and topology of 3--manifolds. In fact the space of normal surfaces of an ideal triangulation forms a natural dual object of the space of angle structures \cite{kan-rub-II, luo-till-angle}.

In our study of positive Euler characteristic normal surfaces in a 3--manifold \cite{jac-rub-0eff}, we developed a technique for crushing a triangulation of a 3--manifold along a normal surface; in this work and subsequent work \cite{jac-rub-annular}, we extend these techniques to manifolds with boundary, crushing a triangulation along the boundary (crushing the boundary to a point) and arriving at a related ideal triangulation of the interior of the 3--manifold. In our considerations of surface with zero Euler characteristic, we discovered an operation on ideal triangulations that is dual to the operation of crushing a triangulation of a 3--manifold with boundary along its boundary. We call this operation on an ideal triangulation an inflation of the ideal triangulation. Starting from an ideal triangulation of the interior of a compact 3--manifold with boundary, an inflation gives a strongly related triangulation of the compact 3-manifold itself, which, in turn, admits a crushing along its boundary returning to the original ideal triangulation.

In Section 3 of this paper, we review the construction of crushing a triangulation of a 3--manifold along a normal surface and apply these techniques to this work, which we distinguish by saying we crush a triangulation along a normal boundary.  Theorem 3.1 can be considered the Fundamental Theorem for Crushing Triangulations along a normal surface.  There can be obstructions to crushing a triangulation along a normal surface; in fact, there are two such obstructions which can be manifested in the natural cell-decomposition coming from splitting a triangulation along a normal surface. We provide examples of the obstructions that can occur in Figure 4: (A) demonstrates what we refer to as ``too many product blocks" and (B) demonstrates ``a cycle of prisms."  We follow these examples with an example in Figure 5 of crushing a triangulation along a normal surface for which there are no obstructions.  Full details and a proof of Theorem 3.1 can be found in \cite{jac-rub-0eff}.  We end Section 3 by introducing the new notion of a {\it combinatorial crushing} of a triangulation along a normal surface. By definition a combinatorial crushing has no obstructions and has a very discrete aspect that may not be the situation in more general crushing without obstructions. 

Section 4 introduces inflations of ideal triangulations.  If $M$ is a compact 3--manifold with boundary, $\T$ is a triangulation of $M$ with all of its vertices in $\partial M$, then by ``crushing the triangulation $\T$ along $\partial M$ means the crushing of the triangulation $\T$ along a \underline{normal} surface that is the frontier of a small regular neighborhood of $\partial M$.  For it to be possible to crush a triangulation along $\bdy M$ it is necessary that a small regular neighborhood of $\bdy M$ be normally isotopic to a normal surface and that the hypothesis of Theorem 3.1 be satisfied.  If these conditions are satisfied, then the process takes some proper subcollection of the tetrahedra of $\T$ and uses the face identifications of $\T$ to give face identifications to the specific subcollection of tetrahedra resulting in an ideal triangulation of $\open{M}$, the interior of $M$, necessarily with fewer tetrahedra than those in $\T$. 

Suppose  $M$ is a compact 3--manifold with boundary, no component of which is a 2--sphere,  and $\T^*$ is an ideal triangulation of $\open{M}$. A triangulation $\T$ of $M$, having all of its vertices in $\bdy M$, is called an {\it inflation} of the ideal triangulation $\T^*$, if there is a combinatorial crushing of $\T$ along $\bdy M$ giving the ideal triangulation $\T^*$. While a combinatorial crushing of a triangulation $\T$ along $\bdy M$ gives a unique ideal triangulation $\T^*$ of $\open{M}$, the operation of inflation, which is dual to crushing, can consist of many choices which lead to possibly inequivalent triangulations of $M$. The construction of an inflation $\T$ of an ideal triangulation $\T^*$ uses all of the tetrahedra of $\T^*$ along with some number of \underline{new} tetrahedra; the precise number of the new tetrahedra necessary to the construction can be determined at the beginning.  Furthermore, crushing the inflated triangulation $\T$ along its boundary eliminates precisely the new tetrahedra that were added to the tetrahedra of $\T^*$ in the inflation construction and gives back the triangulation $\T^*$.  Finally, we remark that any ideal triangulation $\T^*$ of the interior of a compact 3--manifold with boundary, no component of which is a 2--sphere, admits an inflation.  All of the details in the construction of an inflation of any ideal triangulation are given in Section 4.  The precise statement is given in Theorem 4.3. 

In Section 5. we provide two examples of the inflation construction. The first is an inflation of the two-tetrahedron ideal triangulation of the {\it figure-eight}--knot complement. The example has a minimal complexity for the inflation and produces a minimal triangulation of the {\it figure--eight} knot exterior.  recall, that it is necessary that a minimal triangulation of a knot exterior in $S^3$ have precisely one vertex and it must be in its boundary.  The second example is an inflation of the one-tetrahedron ideal triangulation of the Giesking manifold. This is a non-orientable 3--manifold that is double covered by the ideal triangulation of the {\it figure-eight} knot complement.  The inflation in this example gives a compact, non-orientable 3--manifold with a Klein Bottle boundary; it is a seven-tetrahedron triangulation and is, again, a minimal triangulation. However, inflations, even of minimal ideal triangulations, do not need to be minimal. It is not know if a minimal triangulation can always be constructed as an inflation.  

In the Appendix we give the standard ideal triangulation of the Whitehead link complement; we use this to exhibit in Section 4 how certain steps in the construction take care of inflations having multiple ideal vertices.

Applications of the inflation construction are given in \cite{jac-inflate-curves, jac-rub-surface, 
jac-rub-annular}.  In \cite{jac-inflate-curves} we provide a relationship between inflations  and adding two-
handles to the boundary of a 3--manifold. In particular, this
construction, called {\it inflation along a curve}, when used in the inflation of an ideal 
triangulation of the interior of a compact 3--manifold with a torus boundary results 
in a Dehn filling of the compact manifold along the slope of the curve used in the 
inflation.  In \cite{jac-rub-surface}, we provide a relationship between inflations and the (closed)
normal surfaces in an ideal triangulation and the closed normal surfaces in 
any inflation. In particular, we prove that if $\T^*$ is an ideal triangulation of the 
interior of the compact 3--manifold $M$ with boundary, no component of which is a 
2--sphere, and $\T$ is an inflation of $\T^*$, then there is a bijective 
correspondence between the closed normal surfaces in $\T^*$ and those in $\T$. In 
particular, all inflations of an ideal triangulation have isomorphic collections of 
closed normal surfaces. In \cite{jac-rub-annular}, we use the inflation construction as a main 
tool to show that any triangulation of a compact, orientable, irreducible, $\bdy$-irreducible, and anannular $3$--manifold can be modified to an annular-efficient triangulation; i.e., a 0--efficient triangulation so that the only normal annuli with essential boundary are edge-linking. A result of this work, also in \cite{jac-rub-annular}, is that in any annular efficient triangulation of the compact 3--manifold $M$,  there are only finitely many boundary slopes for connected normal surfaces in $\T$ of bounded Euler characteristic.

\section{Triangulations} We follow the notation and basic results of
\cite{jac-rub-0eff} on (pseudo-) triangulations, ideal
triangulations, and normal surface theory.

Suppose $\boldsymbol{\Delta} =
\{\td{\Delta}_1,\ldots,\td{\Delta}_t\}$ is a pairwise-disjoint
collection of compact, convex, linear 3--cells and $\Phi$ is a set
of face pairings on the faces of the cells in $\boldsymbol{\Delta}$
so that if $\phi\in\Phi$, then $\phi$ is an affine isomorphism from
a face $\sigma_i\in\td{\Delta}_i$ to a face
$\sigma_j\in\td{\Delta}_j$, possibly $i = j$.  A face appears in at
most one face pairing and  the natural quotient map
$p:\boldsymbol{\Delta}\rightarrow \boldsymbol{\Delta}/\Phi$ is
injective on the interior of each simplex of each dimension.

Under these conditions, the quotient space
$\boldsymbol{\Delta}/\Phi$ is a $3$--manifold, except possibly at
the image of a vertex or at the image of the midpoint of an edge.  We collect all
this information into a single symbol $\T$ and call $\T$ a {\it
cell-decomposition} of $\boldsymbol{\Delta}/\Phi$, if
$\boldsymbol{\Delta}/\Phi$ is a manifold, or {\it ideal
cell-decomposition} of $\boldsymbol{\Delta}/\Phi$, if $\boldsymbol{\Delta}/\Phi$ is a manifold except possible at the image of a vertex. If
each cell in $\boldsymbol{\Delta}$ is a tetrahedron, we call $\T$ a
{\it triangulation} or {\it ideal triangulation} of
$\boldsymbol{\Delta}/\Phi$. A {\it cell (tetrahedron), face, edge},
or {\it vertex} in this cell decomposition is, respectively, the
image under $p:\boldsymbol{\Delta}\rightarrow
\boldsymbol{\Delta}/\Phi$ of a cell (tetrahedron), face, edge, or
vertex from the collection $\boldsymbol{\Delta} =
\{\td{\Delta}_1,\ldots,\td{\Delta}_t\}$. We will denote the image of
the faces by $\T^{(2)}$, the image of the edges by $\T^{(1)}$ and
the image of the vertices by $\T^{(0)}$. We call $\T^{(i)}$ the {\it
$i$--skeleton of $\T$}; but, generally, we just refer to these as
the faces, edges, or vertices of $\T$. We will denote the image of
$\td{\Delta}_i$ by $\Delta_i$ and call $\td{\Delta}_i$ the {\it lift
of $\Delta_i$}. A cell is the quotient of a unique cell and a face
is the quotient of one or two faces; edges and vertices may be the
quotient of a number of edges or vertices, respectively.  We define
the {\it degree} of an edge $e$ of $\T$ to be the number of edges in
$p^{-1}(e)$.

The collection of normal triangles made up of precisely one normal
triangle of each type forms a normal surface; a component is called
a {\it vertex-linking surface}. $\boldsymbol{\Delta}/\Phi$ is a
3--manifold if and only if each vertex-linking surface is a
2--sphere or a 2--cell (in the latter case $M$ has boundary and the
vertex is in $\bdy M$). Typically, for an ideal triangulation, no
vertex-linking surface is a 2--sphere and all vertex-linking
surfaces are closed; however, such restrictions are not necessary.
The {\it index of an ideal vertex }is the genus of its
vertex-linking surface.

For a triangulation of a $3$--manifold with boundary, if the
frontier of a small regular neighborhood of the boundary is normally
isotopic to a normal surface, then we say the {\it triangulation has
normal boundary}.
In general, for a triangulation of a 3-manifold with boundary,  it is not
necessary that the frontier of a small regular neighborhood of the
boundary be normally isotopic to a normal surface.  For example, if
we layer a tetrahedron along an edge in the boundary of a
triangulation, then the resulting triangulation will not have normal
boundary; in particular, layered triangulations of handlebodies
\cite{jac-rub-layered} do not have normal boundaries. 

We recall some well-known results about triangulations of
3--manifolds.

\begin{thm} A closed 3--manifold admits a triangulation with
precisely one-vertex.\end{thm}

\begin{thm} A compact 3--manifold with boundary, no component of
which is a 2--sphere, admits a triangulation with all vertices in
the boundary and then precisely one vertex in each boundary
component.\end{thm}

In each such situation, we say the manifold has a {\it
minimal-vertex triangulation}.

\begin{thm} The interior of a compact 3--manifold with boundary, no component of
which is a 2--sphere, admits an ideal triangulation.\end{thm}

\section{Crushing triangulations}\label{sec:crushing}

In \cite{jac-rub-0eff} we introduced the notion of ``crushing a
triangulation along a normal surface" and stated and proved the
fundamental theorem for crushing.

Crushing a triangulation of a $3$--manifold along a normal surface
provides a global method for modifying the triangulation. It can be
used to reduce the number of tetrahedra in a given triangulation
\cite{jac-rub-0eff}, construct the prime decomposition of a
$3$--manifold \cite{jac-rub-0eff, jac-rub-sed-find}, construct ideal
triangulations \cite{jac-rub-0eff, jac-rub-annular}, and gain a
better understanding of the normal surfaces in a triangulation
\cite{jac-rub-layered, jac-rub-annular}. In this section we give
definitions and state a special case of the fundamental theorem on
crushing triangulations, which is applicable to our needs in this
work. This version is the inverse of an inflation of an ideal
triangulation. The latter is the main purpose of this paper and is
described in Section \ref{sec:inflations}. In fact, understanding
crushing in this special case provided the motivation and
understanding for developing the inflation construction in the next
section.

\begin{figure}[htbp]
            \psfrag{X}{$X$}
            \psfrag{s}{\small tetrahedron}
            \psfrag{f}{\small face}

             \psfrag{c}{{\tiny crush}}
            \psfrag{e}{\small edge}
            \psfrag{t}{\begin{tabular}{c}
          {\small truncated-tetrahedron}\\
            \end{tabular}}
            \psfrag{p}{\begin{tabular}{c}
            {\small truncated-prism}\\
            \end{tabular}}
            \psfrag{q}{\begin{tabular}{c}
          {\small trianglular}\\
          {\small product block}\\
            \end{tabular}}
            \psfrag{r}{\begin{tabular}{c}
          {\small quadrilateral}\\
          {\small product block}\\
            \end{tabular}}

        \vspace{0 in}
        \begin{center}
        \includegraphics[width=3.5 in]{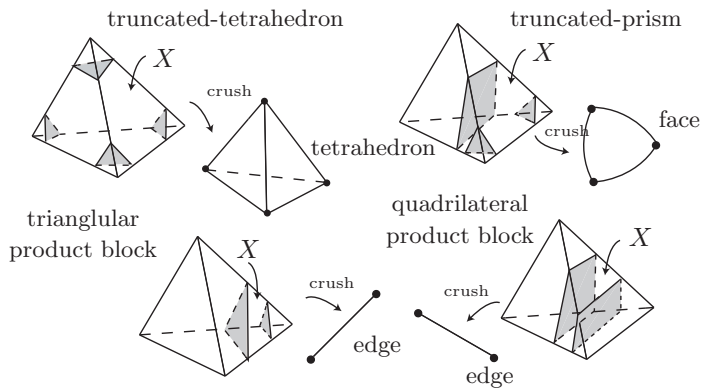}
        \caption{Cells in induced cell-decomposition of $X$ and ideal triangulation of $\open{X}$.}
        \label{f-cell-decomp}
        \end{center}
\end{figure}

Suppose $\T$ is a triangulation of the compact $3$--manifold $M$ or
an ideal triangulation of the interior of $M$. Suppose $S$ is a
closed normal surface embedded in $M$ and $X$ is the closure of a
component of the complement of $S$ in $M$ that does not contain any
vertices of $\T$. For our purposes in this paper, $X$ will be
homeomorphic to $M$. In this situation, we want to use the
tetrahedra of $\T$ to construct a particularly nice ideal
triangulation of $\open{X}$, the interior of $X$. Since none of the
vertices of $\T$ are in $X$, we observe that $X$ has a nice
cell-decomposition, $\mathcal{C}$, consisting of at most four types of cells:
{\it truncated tetrahedra, truncated prisms, triangular product
blocks}, and {\it quadrilateral product blocks}. See Figure
\ref{f-cell-decomp}.

The boundary of each $3$--cell in $\mathcal{C}$ has an induced cell
decomposition in which some of the cells are in $S$ and some are
not. The edges and faces in the decomposition $\mathcal{C}$ are called {\it
horizontal} if their interiors are in $S$ and {\it vertical} if
their interiors are not in $S$. The quadrilateral vertical
$2$--cells are called {\it trapezoids}; there are two in a
truncated-prism, three in a triangular block, and four in a
quadrilateral block. The non-trapezoidal vertical $2$--cells are in
truncated-prisms and truncated-tetrahedra and are hexagons.

We define $\mathbb{P}(\mathcal{C})$ as the union, $\mathbb{P}(\mathcal{C})= \{$vertical edges of $\mathcal{C}\} \cup \{$trapezoids$\}\cup\{$triangular
blocks$\}\cup\{$quadrilateral blocks$\}$. $\mathbb{P}(\mathcal{C})$ is
called the {\it combinatorial product for $\mathcal{C}$}.

Each component of $\bbb{P}(\mathcal{C})$ is an $I$--bundle. Suppose
$\bbb{P}(\mathcal{C})\not= X$ and each component is a product $I$--bundle.
Under these assumptions, a component of $\bbb{P}(\mathcal{C})$ is a product
$\bbb{P}_i = K_i\times[0,1]$, where $K_i$ is isomorphic to a
subcomplex in the induced normal cell structure on $S$, $i =
1,\ldots,k$, and $k$ is the number of components of $\bbb{P}(\mathcal{C})$.
Let $K_{i}^\ve = K_i\times \ve, \ve = 0$ or $1$. Then $K_{i}^0$ and
$K_{i}^1$ are disjoint, isomorphic subcomplexes of the induced
normal cell structure on $S$.

\begin{figure}[htbp]

            \psfrag{o}{\large{or}}
            \psfrag{t}{\begin{tabular}{c}
            truncated-tetrahedron\\
       chain terminates\\
            \end{tabular}}
             \psfrag{p}{\begin{tabular}{c}
            truncated-prism\\
       chain continues\\
            \end{tabular}}
            \psfrag{b}{\begin{tabular}{c}
            face in $\bdy M$\\
       chain terminates\\
            \end{tabular}}
        \vspace{0 in}
        \begin{center}
        \epsfxsize=2.5 in
        \epsfbox{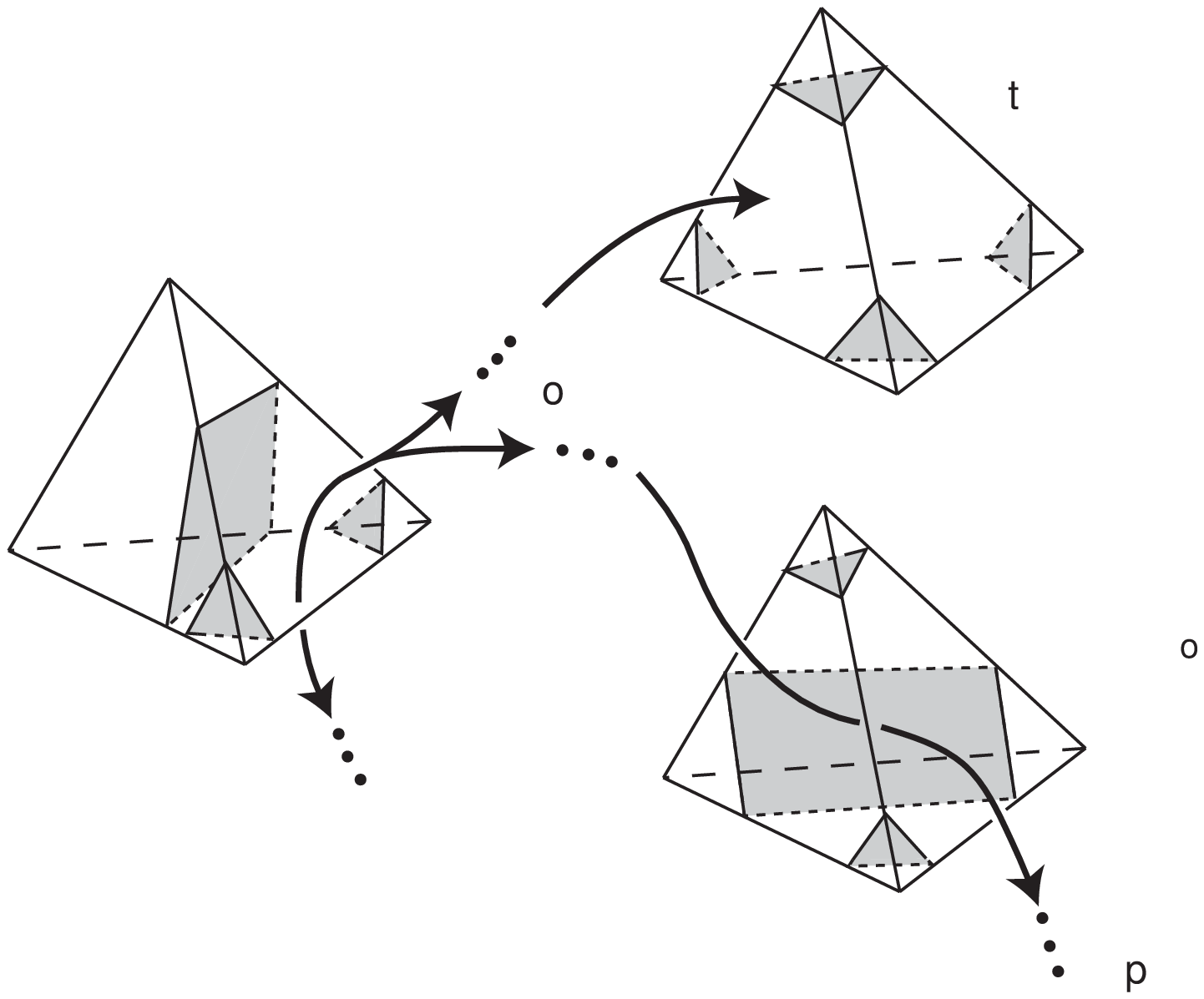}
        \caption{Chain of truncated prisms.}
        \label{f-chain}
        \end{center}
\end{figure}

Now, consider the truncated-prisms in $\mathcal{C}$. Each truncated-prism has
two hexagonal faces. In $\mathcal{C}$, these hexagonal faces are identified
via the face identifications of the given triangulation $\T$ to a
hexagonal face of a truncated-tetrahedron or to a hexagonal face of
truncated-prism. If we follow a sequence of such identifications
through hexagonal faces of truncated-prisms, we trace out a
well-defined arc that terminates at an identification with a
hexagonal face of a truncated-tetrahedron or possibly does not
terminate but forms a complete cycle through hexagonal faces of
truncated-prisms. See Figure \ref{f-chain}. We call a collection of
truncated-prisms identified in this way a {\it chain}. If a chain
ends in a truncated-tetrahedra, we say the chain {\it terminates};
otherwise, we call the chain a {\it cycle of truncated-prisms}.

In general, conditions sufficient for crushing must be
established; for example, in a general situation it may not be true that each $K_{i}\times I$ is a
product $I$-bundle or that each $K_i$ is simply connected, or that
$\bbb{P}(\mathcal{C})\not= X$. However, for the purposes of this work, each
$K_{i}\times I$ is a product $I$-bundle, each $K_i$ is a simply
connected planar complex and hence, cell-like, and $\bbb{P}(\mathcal{C})\not=
X$. Under all these conditions, we say $\bbb{P}(\mathcal{C})$ is a {\it
trivial combinatorial product}. Furthermore, in this work, there are
no cycles of truncated-prisms. In the general situation, we can
allow cycles of truncated-prisms but they must be consumed by more
general product regions than we are considering here. Again, we
refer the reader to \cite{jac-rub-0eff}.

By our assumptions, there are truncated-tetrahedra in the cell
decomposition $\mathcal{C}$ of $X$ (there are not too many product blocks,
$X\ne \bbb{P}(\mathcal{C})$), and there are not too many truncated-prisms (no
cycles of truncated-prisms). To go from the cell decomposition $\mathcal{C}$
of $X$ to an ideal triangulation of $\open{X}$, it is necessary to
crush cells (or collections of cells) of $\mathcal{C}$, arriving at an ideal
triangulation of $\open{X}$. In particular, each component of $S$ is
crushed to a point (distinct points for distinct components), all
products $K_i\times I$ are crushed to arcs (edges) so that if
$K_i\times I$ is crushed to the edge $e_i$, then the crushing
projection coincides with the projection of $K_i\times I$ onto the
$I$ factor.  Each vertical edge, each trapezoid, and each product
block in $\mathcal{C}$ becomes an edge. Each truncated-prism becomes a face
and each truncated-tetrahedron becomes a tetrahedron. See Figure
\ref{f-cell-decomp}.

Let $\{\overline{\Delta}_1,\ldots,\overline{\Delta}_n\}$ be the
collection of truncated-tetrahedra in $\mathcal{C}$. Notice that each
truncated-tetrahedron in $X$ has its triangular faces in $S$. If we
crush each such triangular face of a truncated-tetrahedron to a
point (for the moment, distinct points for each triangular face), we
get a tetrahedron. We use the notation $\td{\Delta}_i^*$ for the
tetrahedron coming from the truncated-tetrahedron
$\overline{\Delta}_i$ after identifying  the triangular faces of
$\overline{\Delta}_i$ to points. If $\overline{\sigma}_i$ is a
hexagonal face in $\overline{\Delta}_i$, then $\overline{\sigma}_i$
is identified to a triangular face, say $\td{\sigma}_i^*$, of
$\td{\Delta}_i^*$.

\begin{figure}[htbp]

          \psfrag{o}{of prisms}
           \psfrag{c}{chain}
             \psfrag{I}{$\overline{\Delta}_i$}
            \psfrag{J}{$\overline{\Delta}_j$}
            \psfrag{i}{$\td{\Delta}_i^*$}
            \psfrag{j}{$\td{\Delta}_j^*$}
            \psfrag{1}{$\overline{\sigma}_i$}
            \psfrag{2}{$\overline{\sigma}_j$}
            \psfrag{3}{$\td{\sigma}_j$}
            \psfrag{4}{$\td{\sigma}_i$}
            \psfrag{a}{\begin{tabular}{c}
            after\\
        crush\\
            \end{tabular}}
            \psfrag{f}{\begin{tabular}{c}
            faces\\      identified\\
            \end{tabular}}
       \vspace{0 in}
        \begin{center}
        \epsfxsize = 3.5 in
        \epsfbox{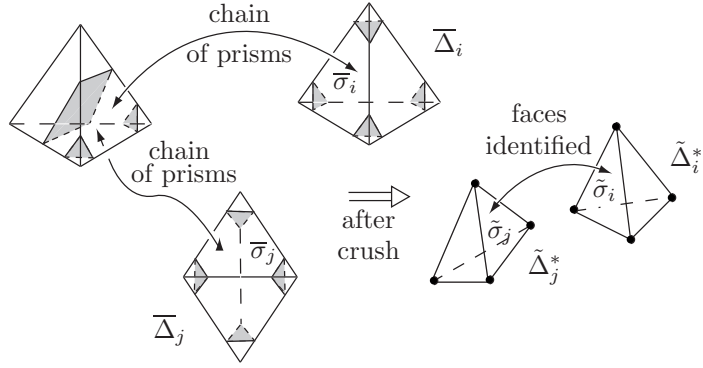}
        \caption{Face identifications induced through a chain of truncated-prisms.}
        \label{f-ident-prisms}
        \end{center}

\end{figure}

Let $\td{\Delta}^* =
\{\td{\Delta}_1^*,\ldots,\td{\Delta}_n^*\}$ be the tetrahedra
obtained from the collection of truncated-tetrahedra
$\{\overline{\Delta}_1,\ldots,\overline{\Delta}_n\}$ following the
crushing of the normal triangles in the surface $S$ to points. It
follows that there is a family $\td{\Phi}^*$ of face-pairings
induced on the collection of tetrahedra
$\td{\Delta}^*$ by the face-pairings of $\mathcal{C}$ (coming
from the face-pairings of $\T$) as follows (see Figure
\ref{f-ident-prisms}):
\begin{itemize}\item[-] if the face $\overline{\sigma}_i$ of
$\overline{\Delta}_i$ is paired with the face $\overline{\sigma}_j$
of $\overline{\Delta}_j$, then this pairing induces the pairing of
the face $\td{\sigma}_i^*$ of $\td{\Delta}_i^*$  with the face
$\td{\sigma}_j^*$ of $\td{\Delta}_j^*$ ;\item [-] if the face
$\overline{\sigma}_i$ of $\overline{\Delta}_i$ is paired with a face
of a truncated-prism in a chain of truncated-prisms and the face
$\overline{\sigma}_j$ of the truncated-tetrahedron
$\overline{\Delta}_j$ is also paired with a face of this chain of
truncated-prisms, then the face $\td{\sigma}_i^*$ of
$\td{\Delta}_i^*$ has an induced pairing with the face
$\td{\sigma}_j^*$ of $\td{\Delta}_j^*$ through the chain of
truncated-prisms.\end{itemize}

Hence, we get a $3$--complex
$\boldsymbol{\td{\Delta}}^*/\boldsymbol{\td{\Phi}}^*$, which is a
$3$--manifold except, possibly, at its vertices. We will denote the
associated ideal triangulation by $\T^*$. We call $\T^*$ the ideal
triangulation obtained by {\it crushing the triangulation $\T$ along
$S$}. We denote the image of a tetrahedron $\td{\Delta}^*_i$ by
$\Delta^*_i$ and, as above, call  $\td{\Delta}^*_i$ the lift of
$\Delta^*_i$.

We have the following theorem.

\begin{thm}\label{thm:crush} Suppose $\T$ is a triangulation of a compact $3$--manifold
or an ideal triangulation of the interior of a compact
$3$--manifold, $M$. Suppose $S$ is a normal surface embedded in $M$,
$X$ is the closure of a component of the complement of $S$, $X$ does
not contain any vertices or ideal vertices of $\T$, and $\mathcal{C}$ is the
induced cell-decomposition on $X$. Let $\bbb{P}(\mathcal{C})$ denoted the
combinatorial product region for $X$. If
\begin{enumerate}
\item[i)] $X\ne \bbb{P}(\mathcal{C})$,  \item[ii)] $\bbb{P}(\mathcal{C})$ is a trivial
product region for $X$, and \item[iii)] there are no cycles of
truncated-prisms in $X$,
\end{enumerate}
then the triangulation $\T$ can be crushed along $S$ giving a unique
ideal triangulation $\T^*$ of $\open{X}$.
\end{thm}
\begin{proof} Using the above notation, we have that $\T$ induces a
cell-decomposition $\mathcal{C}$ on $X$. The truncated-tetrahedra in $\mathcal{C}$ (by
hypothesis, there must be some) determine a collection of tetrahedra
$\td{\Delta}^* =
\{\td{\Delta}_1^*,\ldots,\td{\Delta}_n^*\}$ and, as described above,
the face-pairings of $\T$, along with our hypothesis that there are
no cycles of truncated-prisms, determine a family
$\boldsymbol{\td{\Phi}}^*$ of face-pairings for
$\boldsymbol{\td{\Delta}}^*$. The underlying point set for the
triangulation $\T^*$,
$\boldsymbol{\td{\Delta}}^*/\boldsymbol{\td{\Phi}}^*$, is obtained
from $X$ by identifying each component of $S$ to a point (distinct
points for distinct components), identifying each component, 
$K_i\times [0,1]$ of $\bbb{P}(\mathcal{C})$, of the product region for $X$ to
an edge $e_i$ (distinct edges for distinct components; see
Figure \ref{f-cell-decomp}), and identifying each chain of
truncated-prisms to a face (see Figure \ref{f-ident-prisms}). If we
look at this identification map we have the inverse image of a point
in the interior of a tetrahedron $\Delta^*_i$ is just a point in the
interior of the truncated-tetrahedron $\overline{\Delta}_i$; the
inverse image of a point in the interior of a face is either a point
or an arc, the latter in the case a chain of truncated-prisms is
identified to a face; and the inverse image of a point in the
interior of an edge is a copy $K_j\times x$ for some $j$ and $x\in
[0,1]$. Notice that in the identification of a chain of
truncated-prisms to a face; the associate identification of the
edges is through a band of trapezoids and so there are no new
identifications not already made in $K_j\times [0,1]$ for some $j$.
Thus the identification map on $\open{X}$ is a cell-like map. It
follows by \cite{arm, sie}, that $\T^*$ is an ideal triangulation of
$\open{X}$. Furthermore, there are no choices for the
truncated-tetrahedra; they are completely determined by $\T$ and
$S$. Under our assumptions the truncated-tetrahedra are crushed to
tetrahedra and face identifications of $\T^*$ are completely (and
uniquely) determined by the face-pairings of $\T$. We conclude that
$\T^*$ is uniquely determined by $\T$ and $S$.
\end{proof}

\begin{figure}[htbp]
\psfrag{0}{\footnotesize$0$}\psfrag{1}{\footnotesize$1$}\psfrag{2}{\footnotesize$2$}\psfrag{3}{\footnotesize$3$}
\psfrag{A}{(A) too many product blocks}\psfrag{B}{(B) cycle of
prisms}

            \psfrag{s}{$\sigma'$}
            \psfrag{t}{$\sigma$}\psfrag{x} {\begin{tabular}{c}
            $\sigma(013)\leftrightarrow\sigma'(231)$\\$\sigma(213)\leftrightarrow\sigma'(031)$\\
            $[\sigma(012)\leftrightarrow\sigma'(012)]$\\$[\sigma(023)\leftrightarrow\sigma'(023)]$\\\end{tabular}}
\psfrag{y} {\begin{tabular}{c}
            $\sigma(012)\leftrightarrow\sigma'(023)$\\$\sigma(023)\leftrightarrow\sigma'(031)$\\
            $\sigma(031)\leftrightarrow\sigma'(012)$\\$[\sigma(123)\leftrightarrow\sigma'(123)]$\\\end{tabular}}

\vspace{.25 in}
        \begin{center}
\epsfxsize =5 in 
\epsfbox{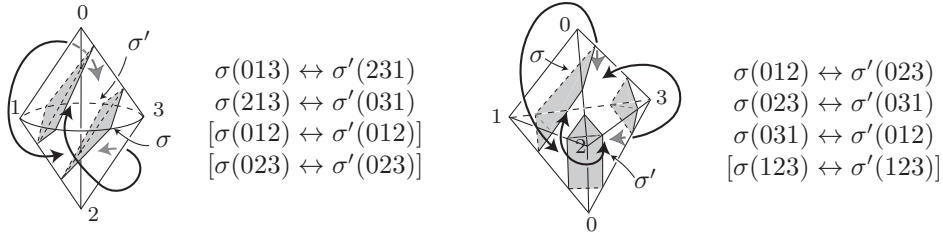} 
\caption{In(A) crushing a two-tetrahedron triangulation of $L(2,1)=\rppp$ along a normal $S^2$ with obstruction $X=\mathcal{P}(X)$. In (B) crushing a two-tetrahedron triangulation of $L(3,1)$ along a normal 2--sphere with obstruction a cycle of prisms.} 
\label{fr-RP3-L31-crush}

\end{center}
\end{figure}

\begin{figure}[htbp]
\psfrag{X}{$(0)$}\psfrag{Y}{$(1)$}\psfrag{Z}{$(2)$}\psfrag{W}{$(3)$}

            \psfrag{x}{$(\overline{0})$}
            \psfrag{y}{$(\overline{1})$}\psfrag{z}{$(\overline{2})$}\psfrag{w}{$(3^*)$}\psfrag{1}{\footnotesize $1$}
            \psfrag{0}{\footnotesize $0$}
\psfrag{2}{\footnotesize $2$}\psfrag{3}{\footnotesize $3$}
\psfrag{L}{$L(4,1)$}\psfrag{a}{\footnotesize
a}\psfrag{b}{\footnotesize b}\psfrag{c}{\footnotesize
c}\psfrag{d}{\footnotesize d}\psfrag{e}{\footnotesize e}

\vspace{.25 in}
        \begin{center}
\epsfxsize =3.5 in 
\epsfbox{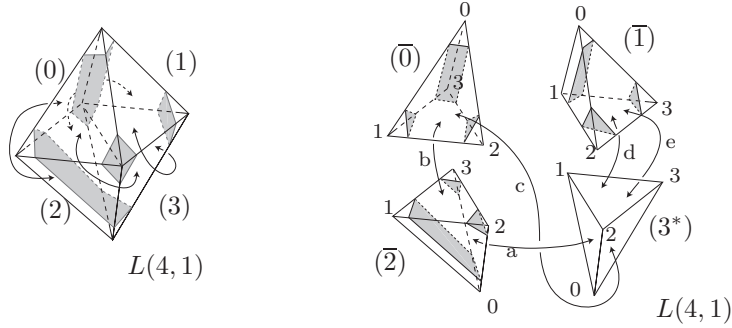}
\caption{Crushing a four-tetrahedron triangulation of $L(4,1)$ along a normal $S^2$ giving the one-tetrahedron triangulation of $L(4,1)$.} 
\label{f-crush-L41-inflate}
\end{center}
\end{figure}

Following are three elementary examples exhibiting the construction
of crushing a triangulation along a normal $2$--sphere. In  Figure
\ref{fr-RP3-L31-crush}(A) the construction terminates when the
induced cell decomposition of $X$ ``has too many product blocks". In
this case, we have that $X =\mathcal{P}(X)$ is a twisted I-bundle
over $\rpp$ and $M$ is $\rppp$. In Figure \ref{fr-RP3-L31-crush}(B)
the construction terminates when the induced cell decomposition of
$X$ ``has a cycle of prisms", giving that $M$ is the 3--manifold
$L(3,1)$. In Figure \ref{f-crush-L41-inflate} the construction
crushes a four-tetrahedron, two-vertex triangulation of $L(4,1)$ to
the one-tetrahedron, one-vertex, minimal triangulation of $L(4,1)$.
Note when there are no obstructions to crushing along a normal
$2$--sphere, the ideal triangulation in the conclusion of Theorem
\ref{thm:crush} gives a triangulation.

\vspace{.125 in}\noindent{\bf Example.} Obstructions when crushing a
triangulation along a normal 2--sphere.  See Figure
\ref{fr-RP3-L31-crush}.

\vspace{.125 in}\noindent{\bf Example.} Crushing a triangulation
along a normal 2--sphere.  See Figure \ref{f-crush-L41-inflate}.

\vspace{.25in}A four-tetrahedron triangulation of $L(4,1)$.
\begin{tabbing}\hspace{.5 in}\=\bf {tet}\hspace{.55
in}\=\hspace{.175in}$(012)$\hspace{.55 in}\=\hspace{.175in}$(013)$
\hspace{.55 in}\=\hspace{.175in}$(023)$\hspace{.55 in}\=\hspace{.175in}$(123)$\\
\>$(0)$\>$(3)(023)$\>$(2)(021)$\>$(1)(021)$
\>$(2)(123)$\\
\>$(1)$\>$(0)(032)$\>$(2)(013)$\>$(3)(031)$\>$(3)(123)$\\
\>$(2)$\>$(0)(031)$\>$(1)(013)$\>$(3)(021)$
\>$(0)(123)$\\
\>$(3)$\>$(2)(032)$\>$(1)(032)$\>$(0)(012)$
\>$(1)(123)$\\
\end{tabbing}

The cell-decomposition of $X$ in this example consists of three
truncated-prisms and one truncated-tetrahedron, $(\overline{3})$.
The three truncated-tetrahedra form two chains, one having two
truncated-prisms denoted $(\overline{2})$ and $(\overline{0})$ and
the other having just one truncated-prism denoted $(\overline{1})$.

The new face identifications after crushing are given as:

 \begin{tabbing}\=$(3^*)(012)\stackrel{a}{\leftrightarrow}(\overline{2})(032)\stackrel{{\footnotesize crush}}{\leftrightarrow}$\=
$(\overline{2})(132)\stackrel{b}{\leftrightarrow}(\overline{0})(132)\stackrel{{\footnotesize
crush}}{\leftrightarrow}$\\
\>\>$\stackrel{{\footnotesize
crush}}{\leftrightarrow}(\overline{0})(102)\stackrel{c}{\leftrightarrow}(3^*)(203)$\\
\>  and\>\\
\>$(3^*)(123)\stackrel{d}{\leftrightarrow}(\overline{1})(123)\stackrel{{\footnotesize
crush}}{\leftrightarrow}(\overline{1})(023)\stackrel{e}{\leftrightarrow}(3^*)(031)$\\\end{tabbing}

The new triangulation of $L(4,1)$ after crushing is the
one-tetrahedron, one-vertex triangulation $\T^*$.

\begin{tabbing}\hspace{.5 in}\=\bf {tet}\hspace{.55
in}\=\hspace{.25in}$(012)$\hspace{.55 in}\=\hspace{.25in}$(013)$
\hspace{.55 in}\=\hspace{.25in}$(023)$\hspace{.55 in}\=\hspace{.25in}$(123)$\\
\>$(3^*)$\>$(3^*)(203)$\>$(3^*)(132)$\>$(3^*)(102)$
\>$(3^*)(031)$\\
\end{tabbing}

We end this section with a definition and an observation. If $M$ is
a 3--manifold, $\T$ is a triangulation or ideal triangulation of
$M$, $S$ is a normal surface, and $X$ is the closure of a component
of the complement of $S$ meeting no vertices of $\T$, then under the
special conditions $X\ne \bbb{P}(X)$, the combinatorial product
$\bbb{P}(X)$ is trivial, and there are no cycles of
truncated-prisms, we have from Theorem \ref{thm:crush} that the
triangulation $\T$ admits a crushing along $S$. In this special
situation, we say the triangulation $\T$ admits a
\emph{combinatorial crushing along $S$}. More general conditions for
crushing are given in \cite{jac-rub-0eff}.  In the case of a combinatorial crushing along $S$, the tetrahedra of the
ideal triangulation $\T^*$ are in one-one correspondence with the
truncated-tetrahedra of the cell-decomposition $\mathcal{C}$ of $X$. Hence,
if $t$ is the number of tetrahedra of $\T$ and $t^*$ is the number
of tetrahedra of $\T^*$, then $t^*\le t$ with equality if and only
if $S$ is a vertex-linking surface, in which case, $\T = \T^*$.

\section{Inflations of ideal triangulations}\label{sec:inflations}
Suppose $X$ is a compact 3--manifold with boundary and $\T$ is a
triangulation of $X$ with normal boundary.  If the triangulation
$\T$ can be crushed along the normal surface that is the frontier of
a small regular neighborhood of the boundary, we say $\T$ {\it
admits a crushing along} $\bdy X$.

\begin{defn}If $\T^*$ is an ideal triangulation of $\open{X}$, the
interior of the compact $3$--manifold $X$, an {\it inflation of
$\T^*$} is a minimal-vertex triangulation $\T$ of $X$ with a normal
boundary that admits a combinatorial crushing along $\bdy X$ giving
the ideal triangulation $\T^*$.\end{defn}

 In this section we
provide an algorithm for constructing an inflation of any given
ideal triangulation of the interior of a compact $3$--manifold.

\subsection{Frames} 
\begin{figure}[htbp]

            \psfrag{L}{\Large{$\xi$}}
        \vspace{0 in}
        \begin{center}
\epsfxsize =2.75 in 
\epsfbox{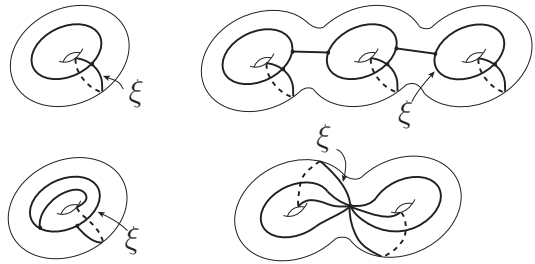} 
\caption{On the left are the only two possible frames for the torus. On the right are two examples of frames: one for the genus two surface and the other for a genus three surface.} 
\label{f-frames}
\end{center}
\end{figure}

\begin{figure}[htbp]

            \psfrag{L}{$\xi = \langle 1\rangle\cup\langle 9,3,\overline{6},\overline{4}\rangle$}
            \psfrag{K}{$\xi = \langle 5\rangle\cup \langle 2\rangle\cup \langle 3\rangle$}
            \psfrag{1}{\footnotesize $1$}\psfrag{2}{\footnotesize$2$}
            \psfrag{3}{\footnotesize$3$}
            \psfrag{4}{\footnotesize$4$}
            \psfrag{5}{\footnotesize$5$}\psfrag{6}{\footnotesize$6$}
            \psfrag{7}{\footnotesize$7$}
            \psfrag{9}{\footnotesize$9$}
        \vspace{0 in}
        \begin{center}
\epsfxsize =3 in 
\epsfbox{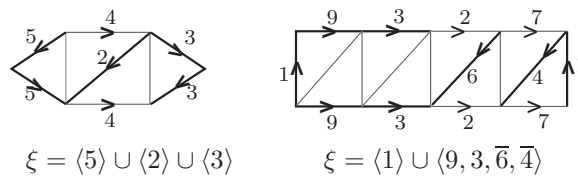} 
\caption{A double index $3$ frame with three branches for the vertex-linking Klein bottle in the ideal triangulation of the Gieseking manifold and an index $4$ frame with two branches for the vertex-linking torus in the two-tetrahedron ideal triangulation of $S^3\setminus$({\it figure-eight}).} 
\label{f-example-frames}
\end{center}
\end{figure}

Suppose $S$ is a triangulated surface. A graph in the
$1$--skeleton of the triangulation of $S$ is called a {\it spine} if
each component of its complement in $S$ is an open disk. We say a spine $\xi$ for $S$
is a {\it frame} if it is minimal with respect to set inclusion;
i.e., if $\xi'$ is a spine for $S$ and $\xi'\subset\xi$, then $\xi'
= \xi$.  Note that a frame has only one component of its complement. A vertex on a frame is called a {\it branch point} if its
index is greater than $2$, in which case its index is called its
{\it branching index}. A component of a frame minus its branch
points is an open arc; we call its closure a {\it branch} of the
frame. For a surface of genus $g$ there are only finitely many
configurations, up to graph isomorphism, for branches and branch
points making up a frame. In Figure \ref{f-frames} we show the only
two possible configurations for the torus and give examples for
frames for genus $2$ and genus $3$ surfaces. In the case of the
torus, we refer to the two possible frames as an {\it index $4$
frame} or a {\it double index $3$ frame}. In Figure
\ref{f-example-frames}, we give explicit examples of frames; one is
a double index $3$ frame in the vertex-linking Klein bottle of the
one-tetrahedron ideal triangulation of the Gieseking manifold and
the other is an index $4$ frame in the vertex-linking torus of the
two-tetrahedron ideal triangulation of the {\it figure-eight} knot
complement in $S^3$. In the latter, the frame is the standard
meridian/longitude frame. The bars on $6$ $(\overline{6})$ and on
$4$ $(\overline{4})$ in Figure \ref{f-example-frames} indicate
traversing the edges $6$ and $4$ in the direction opposite that used
 in the face identifications of the triangulation (see
Figure \ref{f-ideal-8}).

\subsection{Existence of inflations} The construction of an inflation of
an ideal triangulation $\T^*$ begins with the choice of frames in
the induced triangulations of the various vertex-linking surfaces of
$\T^*$. We have organized the construction with notation that,
hopefully, may aid in coding the algorithm so that it can be used
more effectively to generate and study examples.

\begin{figure}[htbp]

            \psfrag{a}{\small {$v_0$}}
            \psfrag{b}{\small{$v_1$}}
            \psfrag{1}{\scriptsize{$x_1$}}
            \psfrag{2}{\scriptsize{$x_2$}}
            \psfrag{3}{\scriptsize{$x_{J-1}$}}
            \psfrag{4}{\scriptsize{$x_J$}}\psfrag{5}{\scriptsize{$z_1$}}
\psfrag{6}{\scriptsize{$z_N$}}\psfrag{7}{\scriptsize{$y_K$}}
\psfrag{8}{\scriptsize{$y_{K-1}$}}
\psfrag{9}{\scriptsize{$y_1$}}\psfrag{L}{\Large{$\xi$}}

        \begin{center}
\epsfxsize =2 in
\epsfbox{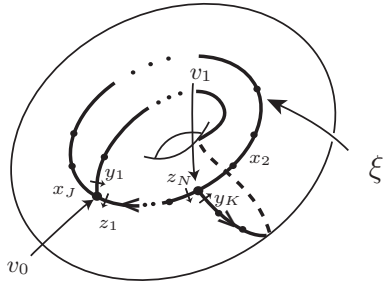} 
\caption{Example of labeling the edges of a frame; the edges in each branch are labeled in succession, beginning (and ending) at a branch point. Transverse directions are shown on two branches.} 
\label{f-label-edges}
\end{center}
\end{figure}

Suppose $\xi$ is a frame in the vertex-linking surface $S$, we label
the edges of $\xi$. We choose to label the edges by selecting some
direction on a branch; the choice of direction is arbitrary. For
such a directed branch, we label the edges successively,
$e_1,e_2,\ldots,e_J$, beginning at the initial branch point
(determined by the chosen direction on the branch) and ending at the
terminal branch point. Of course, the initial and terminal branch
point might be the same. We include additional information (again,
using the direction of the branch) by labeling the initial and
terminal vertices of the edge $e_j$ by $e_j^0$ and $e_j^1$,
respectively; it is possible that these are the same point, for
example, if the branch has only one edge. Hence, if $v_0$ and $v_1$
are the initial and terminal branch points for the branch in
question, we have: $v_0=e_1^0,
e_1^1=e_2^0,\ldots,e_{J-1}^1=e_J^0,e_J^1=v_1$. See Figures
\ref{f-example-frames} and \ref{f-label-edges}; the former gives
examples of actual frames with labeled branches. We label the frames
in all of the vertex-linking surfaces.

In addition to choosing a direction for each branch and labeling its
edges, we choose a transverse direction for each branch. We
consistently choose the transverse direction for a branch by using
the right-hand-rule at its initial vertex; i.e., if at the initial
vertex of the branch, the thumb of the right-hand is pointing in the
direction of the ideal vertex, then the index finger of the right
hand is pointing in the transverse direction to that branch. We then
transport the transverse direction, determined at the initial
vertex, along the branch inducing a transverse direction on each of
its directed edges. For orientable surfaces, the right-hand rule
(described above) can be used at any vertex; however, this is not
the case for non-orientable surfaces. We indicate the transverse
direction by small transverse arrows on two of the branches in
Figure \ref{f-label-edges}; also, see Examples of inflations given
below in Section \ref{blow-up-examples}.

Before we present the inflation construction, we give an overview so
the reader will understand our motivation at various steps in the
construction. An inflation of an ideal triangulation, the very
definition of which involves crushing, is motivated by an attempt to
achieve a model crushing. In crushing a triangulation along a normal
surface, quadrilaterals in the surface lead to prisms or quad
product regions in the tetrahedra containing these quads. The model
situation is to have no cycles of truncated prisms and for quad
products, if there is more than one quad in a tetrahedron, then
there are only two giving a single quad product region between them.
The crushed prisms become faces in the ideal triangulation and the
crushed product regions become edges in the ideal triangulation. If
we consider a copy of the normal surface along which we are crushing
(a parallel, normally isotopic copy), then its image after crushing
becomes a vertex-linking surface; furthermore, each quad in this
surface becomes an edge in the induced triangulation on the
vertex-linking surface. The inflation construction reverses this
model and uses a frame in each of the vertex-linking surfaces as its
guide. In particular, in the inflation construction, we inflate the
vertex-linking surfaces in the ideal triangulation $\T^*$ getting
normal surfaces in the inflation triangulation $\T$. We show that
these surfaces admit a combinatorial crushing that returns to the
starting ideal triangulation $\T^*$; the quads in the induced
triangulation of the surface crush to the edges in the branches of
the frame.

So, how do the frames guide the construction? For each edge in a
frame we add a tetrahedron into the ideal triangulation by what we
call ``an inflation at a face of $\T^*$". The edge in the frame
inflates to a quadrilateral in the inflation of the vertex-linking
surface and is in the added tetrahedron, giving a truncated-prism to
be crushed back to the face, as described in Section
\ref{sec:crushing}. For edges along a branch of the frame, there is,
in general, a unique way to make face identifications for two of the
four faces of each added tetrahedron. The identifications of the
remaining faces are determined at the vertices of the frame, where
an edge in the ideal triangulation meets the vertex-linking surface.
We refer to this part of the construction as ``an inflation at an
edge of $\T^*$". All of the constructions needed for inflating at an
edge of $\T^*$ are combinations of three basic constructions. One is
called generic and is associated with an edge of $\T^*$ that only
meets the frames in a single point of index $2$; another is called a
crossing and is associated with an edge of $\T^*$ that meets the
frames in two distinct points, each of index $2$; the third is
called a branch and is associated with an edge of $\T^*$ meeting the
frames in one point, which is a branch point. In the generic case,
there is only one choice for identification and we do not need to
add any tetrahedra. For a crossing, we need to add a tetrahedron to
make the necessary face identifications of tetrahedra previously
added. For a branch, it is necessary to add a cone over a planar
polygon to make the necessary face identifications of tetrahedra
previously added; then we make some arbitrary choice of subdividing
the polygon (without adding vertices) and cone over the subdivided
polygon to achieve the desired triangulation for the inflation.

\begin{figure}[htbp]

            \psfrag{a}{$\sigma$}
            \psfrag{p}{$(p)$}
            \psfrag{q}{$(p')$}
 \psfrag{S}{\begin{tabular}{c}
           $1$\hspace{.25 in} or \hspace{.25 in}$2$\hspace{.25 in}or\hspace{.25 in} $3$ \\
            edges of frames in the face $\sigma$\\
            \end{tabular}}

        \begin{center}
\includegraphics[width=3.5in]{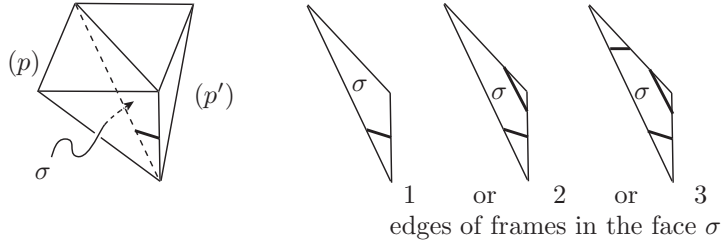}
\caption{The intersection of the frames  with a face $\sigma$ of $\T^*$. The frames can meet $\sigma$ in $1$, $2$, or $3$ edges.} 
\label{f-local-pic-face}
\end{center}
\end{figure}

\begin{figure}[htbp]

            \psfrag{v}{\small $v^*$}
            \psfrag{w}{\small $w^*$}
            \psfrag{E}{\small $E$}\psfrag{1}{\small $D_E^+$}\psfrag{0}{\small $D_E^-$}
\psfrag{2}{\scriptsize $x_j$}
 \psfrag{3}{\scriptsize $x_{j-1}$}
 \psfrag{4}{\scriptsize $\td{x}_j$}
 \psfrag{5}{\scriptsize $\td{x}_{j-1}$}
\psfrag{6}{\scriptsize $y_k$} \psfrag{7}{\scriptsize $y_{k+1}$}
\psfrag{8}{\scriptsize $\td{y}_k$} \psfrag{9}{\scriptsize
$\td{y}_{k+1}$}

        \begin{center}
\includegraphics[width=4in]{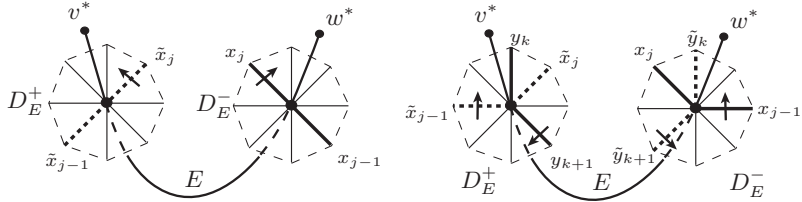}
\caption{The local view of edges and virtual edges of frames at the vertices $E^+$ and $E^-$, where the vertex-linking surfaces meet the edge $E$ of $\T^*$. Here $\td{z}$ denotes the virtual edge associated with $z$.}
\label{f-local-pic-edge}
\end{center}
\end{figure}

In Figure \ref{f-local-pic-face}, we show how the frames can be
viewed in a face of $\T^*$; and in Figure \ref{f-local-pic-edge}, we
show how they can be viewed at an edge of $\T^*$.

An edge $E$ of $\T^*$ meets the vertex-linking surfaces in two
points; we denote these two points by $E^+$ and $E^-$. In
determining these labels, we have implicitly given a direction to
the edge (say, it is directed from $E^-$ to $E^+$); there is no
preferred direction and the choice is arbitrary for each edge of
$\T^*$. However, this direction is important to our construction as
we use it below to assure we get a $3$--manifold when we inflate at
an edge of $\T^*$. Suppose $D_E^+$ and $D_E^-$ are small regular
neighborhoods of $E^+, E^-$, respectively, in the vertex-linking
surfaces. Both $D_E^+$ and $D_E^-$ receive induced subdivisions into
triangles from the triangulation on the vertex-linking surface.

If $\sigma$ is a face of $T^*$ having $E$ as an edge there are
unique edges, one in $D_E^+$ and one in $D_E^-$, lying in $\sigma$
and meeting $E$; we say one of these edges is {\it above} the other
(relative to $E$). If the edge $x$ is above the edge $y$, then $y$
is also above $x$. The combinatorial structures induced on $D_E^+$
and $D_E^-$ by the triangulation of the vertex-linking surfaces are
isomorphic via the correspondence that takes the edge $x$ in $D_E^+$
$(D_E^-)$ to the edge above it in $D_E^-$ $(D_E^+)$. If $x$ is an
edge in a frame $\xi$ of one of the vertex-linking surfaces and $x$
has a vertex at $E^-$ $(E^+)$, then we call the edge in the
vertex-linking surface above $x$ a {\it virtual edge} of $\xi$. It
is possible that a virtual edge  of $\xi$ over $x$ is also an edge
$z$ of a frame (possibly $\xi$); if this is the case, then $x$ is a
virtual edge over $z$. Each edge in a frame has two virtual edges
associated with it. We remark that the local structure of edges and
virtual edges of the frames about an edge can take numerous forms.
Below, we catalog all of the possibilities for an inflation of an
ideal triangulation having only one ideal vertex being of index one.
In Figure \ref{f-local-pic-edge}, we show a local picture of the
vertex-linking surfaces at each end of the edge $E$ along with edges
and virtual edges of the frames meeting $E$. We have presented the
figure with the ideal vertices $v^*$ and $w^*$, as well as a
transverse direction on the edge $x_j$. Note that the transverse
direction for a virtual edge is taken from that induced on the edge
and as such follows the rule for the edge (which in Figure
\ref{f-local-pic-edge} looks like a left-hand-rule on the virtual
edge at $D_E^+$).

\vspace{.125 in}\noindent{\bf Inflation at a face of
$\mathbf{\T^*}$.} Given the triangulation $\T^*$; that is, we have
the collection of tetrahedra  and the associated family of face
identifications of $\T^*$. We shall discard some of the face
identifications, add tetrahedra and make new face identifications.
We  consider the vertices of the tetrahedra in $\T^*$, the ideal
vertices, as being included.

\begin{figure}[htbp]

            \psfrag{s}{$\sigma$}
            \psfrag{p}{ $(p)$}
            \psfrag{q}{$(p')$}\psfrag{d}{\small $d$}\psfrag{a}{\small $a$}
\psfrag{b}{\small $b$}\psfrag{c}{\small $c$}\psfrag{2}{\small
$d'$}\psfrag{0}{\small $a'$} \psfrag{1}{\small
$b'$}\psfrag{3}{\small $c'$}\psfrag{i}{identify}
\psfrag{e}{$(p)(abc)\leftrightarrow(p')(a'b'c')$}

        \begin{center}
\epsfxsize =3 in 
\includegraphics[width=3.5in]{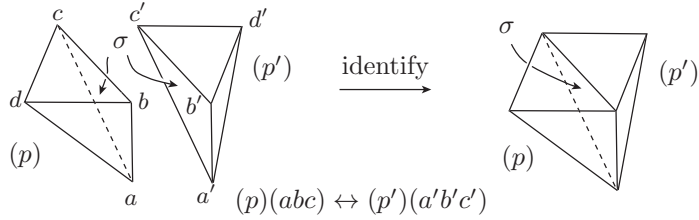} 
\caption{Two tetrahedra of $\T^*$ meeting along the face $\sigma$. We use the notation of {\footnotesize\textsc REGINA} to label simplicies and give face identifications.} 
\label{f-face-ident}
\end{center}
\end{figure}

Suppose $\sigma$ is a face of $\T^*$ and the frames meet $\sigma$.
Let $(p)$ and $(p')$ denote the tetrahedra in $\T^*$ having $\sigma$
as a face and suppose $(p)(abc)=\sigma=(p')(a'b'c')$ is the face
identification. See Figure \ref{f-face-ident}.

As mentioned above, we add a tetrahedron to the triangulation $\T^*$
for each edge of the frame(s). If $x_j$ is an edge of a frame we
denote the tetrahedron to be added by $(x_j)$ and think of it as the
join of two edges, one with vertices $0, 1$ and the other with
vertices $2,3$; hence, the vertices of $(x_j)$ are labeled
$0,1,2,3$. We follow the notational conventions of
{\footnotesize\textsc REGINA} \cite{burton-regina}; i.e., the faces
are $(x_j)(012), (x_j)(013), (x_j)(023),\dots$; edges are
$(x_j)(01), (x_j)(02),\ldots$; vertices are $(x_j)(0),\ldots$, etc.
The choice of the edge with vertices $2$ and $3$ in $(x_j)$ is
arbitrary; however, later, this choice will be significant in our
choice of face identification.

\begin{figure}[htbp]

           \psfrag{1}{\small $1$}\psfrag{0}{\small $0$}
 \psfrag{X}{$(x_j)$}
\psfrag{j}{\small $2$}\psfrag{J}{\small $3$}

        \begin{center}
\includegraphics[width=2.5 in]{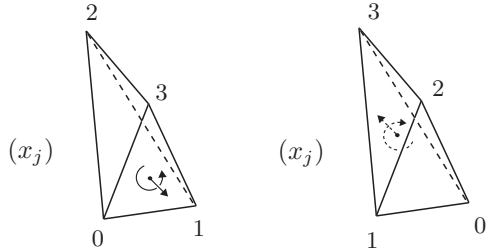} 
\caption{The added tetrahedron $(x_j)$ is a join of edge $(23)$ with $(01)$ where the vertices $0$ and $1$ are chosen after $2$ and $3$, using the ``right-hand rule".} 
\label{f-direct-edges}
\end{center}
\end{figure}

The choice for the vertices $0$ and $1$ is made after that for $2$
and $3$ and uses the convention of labeling the edge $(x_j)(01)$ so
that a right-hand twist while going along the edge $(x_j)(23)$ from
$(x_j)(2)$ to $(x_j)(3)$ moves the vertex $0$ to the vertex $1$ in
the faces $(x_j)(301)$ and $(x_j)(201)$. We indicate this in the
tetrahedra in Figure \ref{f-direct-edges}.

We now give the construction for an inflation at a face $\sigma$ of
$\T^*$, which meets the frames.

\begin{figure}[htbp]

            \psfrag{s}{$\sigma$}
            \psfrag{p}{$(p)$}
            \psfrag{q}{$(p')$}\psfrag{1}{\small $1$}\psfrag{0}{\small $0$}
\psfrag{2}{\small $2$}\psfrag{3}{\small
$3$}\psfrag{a}{$a$}\psfrag{b}{ $b$}
\psfrag{c}{$c$}\psfrag{d}{$d$}\psfrag{-}{\scriptsize{$x_j^0$}}\psfrag{5}{\small
$a'$}\psfrag{6}{\small $b'$} \psfrag{7}{\small
$c'$}\psfrag{8}{\small $d'$}
\psfrag{+}{\scriptsize{$x_j^1$}}\psfrag{x}{\footnotesize{$x_j$}}\psfrag{X}{$(x_j)$}
\psfrag{B}{inflate}\psfrag{f}{face}\psfrag{j}{$2$}\psfrag{J}{$3$}

\psfrag{e}{$(p)(abc)\leftrightarrow(x_j)(032)$ and
            $(x_j)(132)\leftrightarrow(p')(a'b'c')$}

        \begin{center}
\epsfxsize =5 in 
\includegraphics[width=4.5 in]{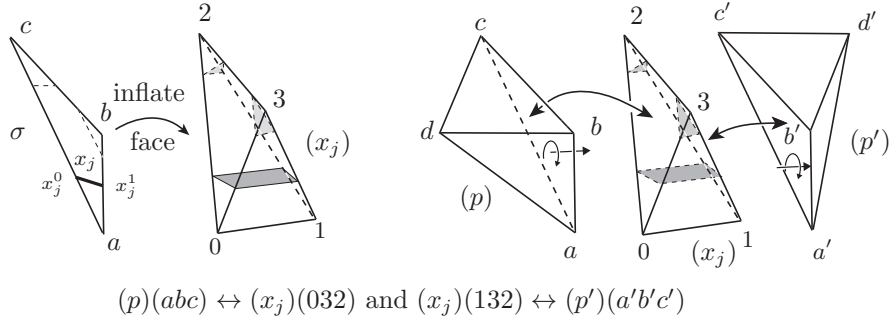} 
\caption{The face $\sigma$ meets the frames in one edge $x_j$; an inflation at the face $\sigma$ adds one tetrahedron.} 
\label{f-blow-up-face1}
\end{center}
\end{figure}

\vspace{.125 in}\noindent{\it One edge in $\sigma$.} Suppose
$\sigma$ contains just one edge of the frame(s). As an edge of a
frame it has been labeled and been given a direction and a
transverse direction. Suppose its label is $x_j$. In this case, we
add one tetrahedron, denoted $(x_j)$ with vertices $0,1,2,3$.  We
discard the face identification $(p)(abc)\leftrightarrow
(p')(a'b'c')$ and add two new face identifications.  The rule for
the new face identifications is that in the face $(p)(abc)$ the edge
of the face opposite the vertex with the edge $x_j$ of the frame is
identified to the edge $(x_j)(23)$ in the same direction as the
directed edge $x_j$; and similarly, in the face $(p')(a'b'c')$ the
edge of the face opposite the vertex with the edge $x_j$ of the
frame is identified to the edge $(x_j)(23)$ in the same direction as
the directed edge $x_j$. In this case, we have the edge $(p)(cb)$
identified with $(x_j)(23)$ and the edge $(p')(c'b')$ also
identified with $(x_j)(23)$. It then needs to be determined which
vertices of $(x_j)$ the vertices $(p)(a)$ and $(p')(a')$ are to be
identified with; one to be identified with $(x_j)(0)$ and the other
with $(x_j)(1)$. The rule for these last identifications is
determined by the transverse direction to the edge $x_j$. If the
transverse direction along $x_j$ points \underline{out of} the
tetrahedron $(p)$, then the vertex $p(a)$ is identified with the
vertex $(x_j)(0)$, leaving the vertex $p'(a')$ to be identified with
the vertex $(x_j)(1)$ and the face identifications are
$(p)(abc)\lra(x_j)(032)$ and $(x_j)(132)\lra (p')(a'b'c')$. In the
exhibited case (Figure \ref{f-blow-up-face1}), the transverse
direction to the edge $x_j$ points \underline{out of} the
tetrahedron $(p)$; hence, the vertex $p(a)$ is identified with the
vertex $(x_j)(0)$, leaving the vertex $p'(a')$ to be identified with
the vertex $(x_j)(1)$ (the transverse direction to the edge $x_j$
points {\it into} the tetrahedron $(p')$). The face identifications
are $(p)(abc)\lra(x_j)(032)$ and $(x_j)(132)\lra (p')(a'b'c')$.  A
quadrilateral is added to the vertex-linking surface for $x_j$ and
two triangles are added for the two virtual edges corresponding to
$x_j$. See Figure \ref{f-blow-up-face1}.

\vspace{.125 in}\noindent{\it Two edges in $\sigma$.} Suppose
$\sigma$ contains two edges of the frames. As edges of frames, they
have labels and directions; suppose their labels are $x_j$ and
$y_k$. In this case, we add two tetrahedra; one denoted $(x_j)$ with
vertices $0,1,2,3$ and the other denoted $(y_k)$ with vertices
$0,1,2,3$. We discard the face identification
$(p)(abc)\leftrightarrow (p')(a'b'c')$ and add three new face
identifications $(p)(abc)\leftrightarrow(x_j)(032)$;
$(x_j)(132)\leftrightarrow(y_k)(312)$; and
$(y_k)(302)\leftrightarrow(p')(a'b'c')$.

\begin{figure}[htbp]

            \psfrag{s}{$\sigma$}
            \psfrag{p}{$(p)$}
            \psfrag{q}{$(p')$}\psfrag{1}{\small $1$}\psfrag{0}{\small $0$}
\psfrag{2}{\small $2$}\psfrag{3}{\small
$3$}\psfrag{a}{$a$}\psfrag{b}{ $b$}
\psfrag{c}{$c$}\psfrag{d}{$d$}\psfrag{-}{\scriptsize{$x_j^0$}}\psfrag{5}{\small
$a'$}\psfrag{6}{\small $b'$} \psfrag{7}{\small
$c'$}\psfrag{8}{\small $d'$}\psfrag{-}{\scriptsize{$x_j^0$}}
\psfrag{+}{\scriptsize{$x_j^1$}}\psfrag{m}{\scriptsize{$y_k^0$}}
\psfrag{n}{\scriptsize{$y_k^1$}}
\psfrag{x}{\footnotesize{$x_j$}}\psfrag{y}{\footnotesize{$y_k$}}\psfrag{X}{$(x_j)$}
\psfrag{Y}{$(y_k)$} \psfrag{j}{$2$}\psfrag{J}{$3$}
\psfrag{k}{$2$}\psfrag{K}{$3$}
\psfrag{e}{$(p)(abc)\leftrightarrow(x_j)(032)$;
$(x_j)(132)\leftrightarrow(y_k)(312)$;
            $(y_k)(302)\leftrightarrow(p')(a'b'c')$}\psfrag{P}{\small{pyramid}}
\psfrag{B}{\begin{tabular}{c}
           inflate \\
           face\\
            \end{tabular}}

        \begin{center}
\epsfxsize =4.5 in 
\includegraphics[width=4.5in]{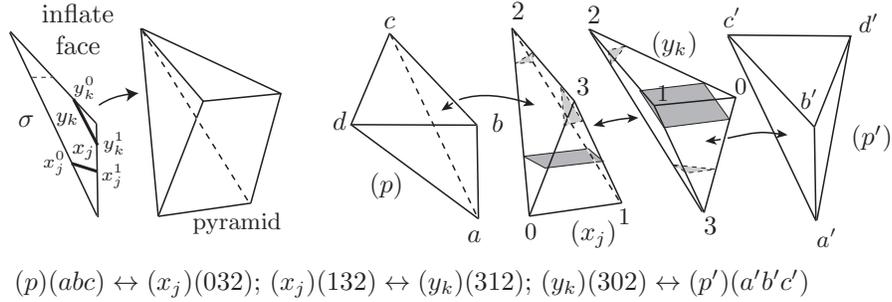} 
\caption{The face $\sigma$ meets the frames in two edges $x_j$ and $y_k$; an inflation at the face $\sigma$ adds two tetrahedra.} 
\label{f-blow-up-face2}
\end{center}
\end{figure}

The rule for the new face identifications is just as that above and
is determined by the direction of the edge in the frame and the
transverse direction. See Figure \ref{f-blow-up-face2}. In the face
$(p)(abc)$ the edge of the face opposite the vertex with the edge
$x_j$ of the frame is $(p)(cb)$; it is identified to the edge
$(x_j)(23)$ in the same direction as the directed edge $x_j$. Since
the transverse direction is pointing out of the tetrahedron $(p)$,
the vertex $(p)(a)$ is identified with $(x_j)(0)$. This gives the
face identification $(p)(abc)\leftrightarrow(x_j)(032)$. The edge
$y_k$ carries its direction and transverse direction to the
tetrahedron $(x_j)$ and we have $(x_j)(21)$ being the edge opposite
$y_k$ and in the same direction. Hence, $(x_j)(21)$ is identified
with $(y_k)(23)$. The transverse direction on $y_k$ is pointing into
$(p)$ (out of $(p')$); this is carried over to $(x_j)$ and we
identify $(x_j)(3)$ with $(y_k)(1)$. This gives the face
identification $(x_j)(132)\leftrightarrow(y_k)(312)$. Finally, we
have $(y_k)(302)\leftrightarrow(p')(a'b'c')$.

Two quadrilaterals are added to the vertex-linking surfaces for
$x_j$ and $y_k$ and four triangles are added for the two virtual
edges corresponding to each. Again, see Figure
\ref{f-blow-up-face2}.

Note that it does not matter in our construction of an inflation
whether we add the tetrahedron $(x_j)$ along $(p)$ or add $(y_k)$
along $(p)$. In fact, one can see that this construction is the same
as adding a single pyramid (also shown in Figure
\ref{f-blow-up-face2}) and then selecting one of the two diagonals
in the quadrilateral face to subdivide it into two tetrahedra.
Making one choice verses the other at this step comes up later when
it may be necessary to add an additional tetrahedron in order to
make necessary face identifications in inflating at the edges of
$\T^*$. Below we consider economy in the number of tetrahedra added
in an inflation and for this the order does matter but we will see a
good choice (economic) is dictated by the frame.

\begin{figure}[htbp]

            \psfrag{s}{$\sigma$}
            \psfrag{p}{$(p)$}
            \psfrag{q}{$(p')$}\psfrag{1}{\small $1$}\psfrag{0}{\small $0$}
\psfrag{2}{\small $2$}\psfrag{3}{\small
$3$}\psfrag{a}{$a$}\psfrag{b}{$b$}
\psfrag{c}{$c$}\psfrag{d}{$d$}\psfrag{5}{\small
$a'$}\psfrag{6}{\small $b'$} \psfrag{7}{\small
$c'$}\psfrag{8}{\small $d'$}\psfrag{-}{\scriptsize{$x_j^0$}}
\psfrag{+}{\scriptsize{$x_j^1$}}\psfrag{m}{\scriptsize{$y_k^0$}}
\psfrag{n}{\scriptsize{$y_k^1$}}\psfrag{r}{\scriptsize{$z_n^0$}}
\psfrag{t}{\scriptsize{$z_n^1$}}
\psfrag{x}{\footnotesize{$x_j$}}\psfrag{y}{\footnotesize{$y_k$}}
\psfrag{z}{\footnotesize{$z_n$}}\psfrag{X}{$(x_j)$}\psfrag{Y}{$(y_k)$}\psfrag{Z}{$(z_n)$}
\psfrag{r}{\scriptsize{$z_n^0$}}
 \psfrag{e}{\begin{tabular}{c}
$(p)(abc)\leftrightarrow(x_j)(032)$;\\
\\
$(x_j)(132)\leftrightarrow(z_n)(231)$;\\
\\
$(z_n)(230)\leftrightarrow(y_k)(203)$;\\
\\
$(y_k)(213)\leftrightarrow(p')(a'b'c')$\\
\end{tabular}}
            \psfrag{P}{\small{prism}}
\psfrag{B}{\begin{tabular}{c}
           \small{inflate} \\
           \small{face}\\
            \end{tabular}}

        \begin{center}
 
\includegraphics[width=4in]{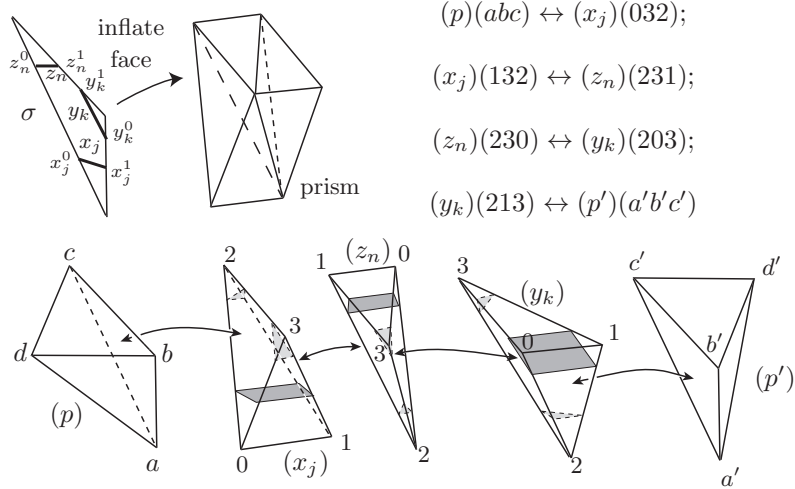} 
\caption{The face sigma meets the frames in three edges $x_j$, $y_k$, and $z_n$; an inflation at the face $\sigma$ adds three tetrahedra.} 
\label{f-blow-up-face3}
\end{center}
\end{figure}

\vspace{.125 in}\noindent{\it Three edges in $\sigma$.} Suppose the
frames meet the face $\sigma$ in three edges. The edges in the
frames have labels and directions; suppose the labels are
$x_j,\thinspace y_k$ and $z_n$. In this case, we add three
tetrahedra, denoted $(x_j), (y_k)$, and $(z_n)$; we label the
vertices of $(x_j)$ as $0,1,2,3$, those of $(y_k)$ as $0,1,2,3$, and
those of $(z_n)$ as $0,1,2,3$. We discard the face identification
$(p)(abc)\leftrightarrow (p')(a'b'c')$ and add four new face
identifications. Again, the choice for the order we add the
tetrahedra is arbitrary. For this demonstration we shall add the new
tetrahedra in the order $(x_j), (z_n)$ and then $(y_k)$. Also, we
need to assume some direction on the edges $x_j, z_n$ and $y_k$, as
well as transverse directions. The new face identifications are:
$(p)(abc)\leftrightarrow(x_j)(032)$;
$(x_j)(132)\leftrightarrow(z_n)(231)$;
$(z_n)(230)\leftrightarrow(y_k)(203)$; and
$(y_k)(213)\leftrightarrow(p')(a'b'c')$. See Figure
\ref{f-blow-up-face3}. In the tetrahedron $(x_j)$ the edge
$(x_j)(13)$ has the direction of $z_n$ and, hence,  is identified to
$(z_n)(23)$; similarly, in $(z_n)$ the edge $(z_n)(20)$ has the
direction of $(y_k$) and, hence, is identified with $(y_k)(23)$.
Three quadrilaterals are added to the vertex-linking surface, one
for each of the edges $x_j$, $y_k$, and $z_n$, and two triangles are
added for the two virtual edges corresponding to each added
tetrahedron, making six triangles added to the vertex-linking
surface(s).

\vspace{.25 in}

\begin{figure}[htbp]

            \psfrag{3}{\begin{tabular}{c}
           \small{$3$ tetrahedra} \\
           \small{to subdivide}\\
            \end{tabular}}
            \psfrag{4}{\begin{tabular}{c}
           \small{$4$ tetrahedra} \\
           \small{to subdivide}\\
            \end{tabular}}

        \begin{center}
\includegraphics[width=2in]{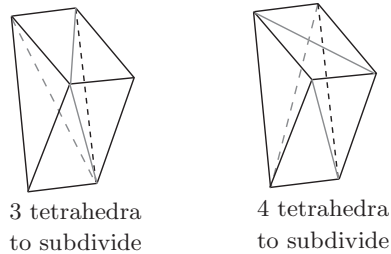}
\caption{The configuration of diagonals in the quadrilateral faces of the prism on the left requires only three tetrahedra to triangulate the prism while the configuration on the right requires four tetrahedra to triangulate.} 
\label{f-prism-diag}
\end{center}
\end{figure}

\vspace{.125 in}  For a face meeting three edges of the frames, one
can see that the construction is the same as adding a single prism
(also shown in Figure \ref{f-blow-up-face3}) and then making a
choice of diagonals in the quadrilateral faces to subdivide the
prism. In this case, however, the choice of diagonals in the quad
faces of a prism impact on further subdivision of the prism; some
choices can be extended to a subdivision requiring only three
tetrahedra and for other choices it is necessary to use four
tetrahedra to subdivide the prism. See Figure \ref{f-prism-diag}. In
adding tetrahedra as we have done,  the diagonals are chosen so that
we only need three tetrahedra; adjustments may then need to be made
when we address necessary face identifications in inflating at the
edges of $\T^*$. As mentioned above in the case of two edges of the
frames in a face, such choices need much more scrutiny when we come
to the issue of economy in adding tetrahedra.

\vspace{.25 in}

\begin{remark} When we inflate $\T^*$ at a face, we
have well-define face identifications for two of the four faces of
each added tetrahedron; at this step, we leave two of the faces in
these tetrahedron unidentified. We have chosen notation so that these
two unidentified faces in the tetrahedron $(x_j)$ are $(x_j)(012)$
and $(x_j)(013)$; the edge common to these faces is $(x_j)(01)$,
which we shall refer to as a {\it free edge}. Below these faces will
be identified to other faces when we inflate at the edges of $\T^*$.
After we have completed the inflation construction, we shall see
that all of the free edges coming from adding the tetrahedra
$(x_1)\ldots,(x_j),\ldots,(x_J)$ to $\T^*$ along a fixed branch of
of a frame having edges labeled $x_1,\ldots,x_j,\ldots,x_J$ are
identified to a single edge $e_x$ in the inflation $\T$ and $e_x$ is
in $\bdy X$.

For later reference, we also note that for each edge in the frame,
we can modify the vertex-linking surface in $\T^*$ by adding a
quadrilateral about the added free edge and two triangles, one at
each end of the edge opposite the free edge in the added
tetrahedron. Hence, there is an ``inflation" of the vertex-linking
surface about each ideal vertex. These inflated surfaces will become
boundary-linking surfaces in the inflated triangulation; i.e., the
frontiers of a small regular neighborhood of the boundary components
in the inflated triangulation. These added quads and triangles are
shown in Figures \ref{f-blow-up-face1}, \ref{f-blow-up-face2}, and
\ref{f-blow-up-face3}.

Finally, we note that at this stage, the vertices of the 3-complex
we have constructed are exactly the ideal vertices we started with
for $\T^*$.
\end{remark}

\vspace{.125 in}\noindent{\bf Inflation at an edge of
$\mathbf{\T^*}$.} Following the inflation at faces of $\T^*$ there
are a number of unidentified faces of the added tetrahedra. These
unidentified faces correspond to the vertices of the frame in the
sense that if $x_j$ is an edge in the frame and $(x_j)$ is a
tetrahedron added when we inflate the face containing $x_j$, then
the faces $(x_j)(012)$ and $(x_j)(013)$ are faces not identified in
that step. The former is associated with the vertex $x_j^0$ of $x_j$
and the latter with the vertex $x_j^1$. The vertices of the
vertex-linking surfaces are precisely where the edges of the
triangulation $\T^*$ meet the vertex-linking surface and the
vertices we are interested in are where the edges of $\T^*$ meet the
frames.

Suppose $E$ is an edge of $\T^*$ and $E$ meets the frames; then
$E^+$ or $E_-$ is a vertex of a frame; it is possible that both are
vertices of frames and that they are vertices of different frames.
Locally about $E$ there are edges and virtual edges of the frames
meeting $E$. The arrangement of edges and virtual edges  in $D_E^+$ meeting
$E^+$  is isomorphic to the arrangement of those in $D_E^-$
meeting $E^-$  via the isomorphism between the
combinatorial structures induced on $D_E^+$ and $D_E^-$ that takes
an edge to the edge above it. In this isomorphism, edges go to
virtual edges and vice-versa.

\begin{figure}[htbp]

            \psfrag{x}{\footnotesize{$x_j$}}\psfrag{y}{\footnotesize{$y_k$}}
            \psfrag{w}{\footnotesize{$\td{w}_l$}}
            \psfrag{v}{\footnotesize{$v_m$}}\psfrag{z}{\footnotesize{$\td{z}_n$}}
            \psfrag{1}{\footnotesize{$\td{x}_j$}}\psfrag{2}{\footnotesize{$\td{y}_k$}}
            \psfrag{3}{\footnotesize{$w_l$}}\psfrag{4}{\footnotesize{$\td{v}_m$}}
            \psfrag{5}{\footnotesize{$z_n$}}
            \psfrag{+}{\footnotesize{$+$}}
            \psfrag{-}{\footnotesize{$-$}}\psfrag{E}{$E$}
            \psfrag{D}{$D_E^+$}\psfrag{d}{$D_E^-$}\psfrag{o}{or}
        \vspace{0 in}
        \begin{center}
\epsfxsize =3.5 in
\includegraphics[width=3.5in]{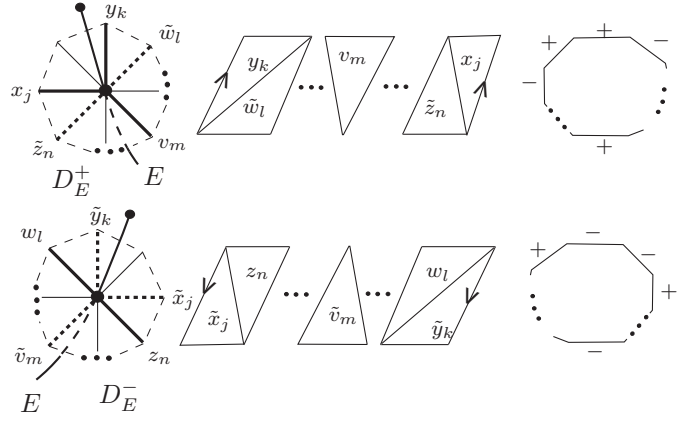} 
\caption{Unidentified faces about the edge $E$ form a band of triangles. The polygons record the pattern and are used in the algorithm to determine both necessary tetrahedra to be added and face identifications for an inflation at the edge $E$ of $\T^*$.} 
\label{f-band-at-E}
\end{center}
\end{figure}

If $x_j$ is an edge of a frame in $D_E^+$, then $E^+$ is a vertex of
$x_j$. If $E^+ = x_j^1$, then the face $(x_j)(013)$ is unidentified
at $E$; and if $E^+ = x_j^0$, then the face $(x_j)(012)$ is
unidentified at $E$. If $\td{z}_n$ is a virtual edge of a frame in
$D_E^+$, then $E^-$ is a vertex of $z_n$. If $E^- = z_n^1$, then the
face $(z_n)(013)$ is unidentified at $E$; and if $E^- = z_n^0$, then
the face $(z_n)(012)$ is unidentified at $E$. Starting at any point
in $D_E^+$ $(D_E^-)$ and making a complete cycle about the boundary
of $D_E^+$ $(D_E^-)$ we can sequentially list the edges and virtual
edges of the frames meeting $E$. Suppose
$x_j,y_k,\td{w}_l,\ldots,v_m,\ldots,\td{z}_n$ is such a listing, say
at $D_E^+$, where  $\td{e}$ denotes the virtual edge of $e$. Then
the unidentified faces of the tetrahedra, added at an inflation of
the faces containing the edges $x_j,y_k,w_l\ldots,v_m,\ldots,z_n$,
form a band of triangles identified as shown in Figure
\ref{f-band-at-E}. It is possible, and most likely, that some of the
edges involved in such a sequence are in the same branch; and if
$E^+$ is an index $2$ vertex in a branch, then consecutive edges of
that branch appear in the sequence but are not necessarily
consecutive in the cyclic order of edges in the frames at $E$.

\begin{figure}[htbp]

            \psfrag{j}{\footnotesize{$x_j$}}
            \psfrag{J}{\footnotesize{$x_{j+1}$}}\psfrag{k}{\footnotesize{$\td{y}_{k-1}$}}
            \psfrag{y}{\footnotesize{$\td{y}_K$}}\psfrag{K}{\footnotesize{$\td{y}_k$}}
            \psfrag{z}{\footnotesize{$\td{z}_N$}}
            \psfrag{w}{\footnotesize{$\td{w}_1$}}\psfrag{1}{\footnotesize{$\td{y}_1$}}
            \psfrag{v}{\footnotesize{$v_m$}}\psfrag{z}{\footnotesize{$\td{z}_N$}}
            \psfrag{+}{\footnotesize{$+$}}
            \psfrag{-}{\footnotesize{$-$}}
            \psfrag{D}{$D_E^+$}\psfrag{d}{$D_E^-$}\psfrag{A}{(A)}\psfrag{B}{(B)}\psfrag{C}{(C)}
        \vspace{0 in}
        \begin{center}
 
\includegraphics[width=3.5in]{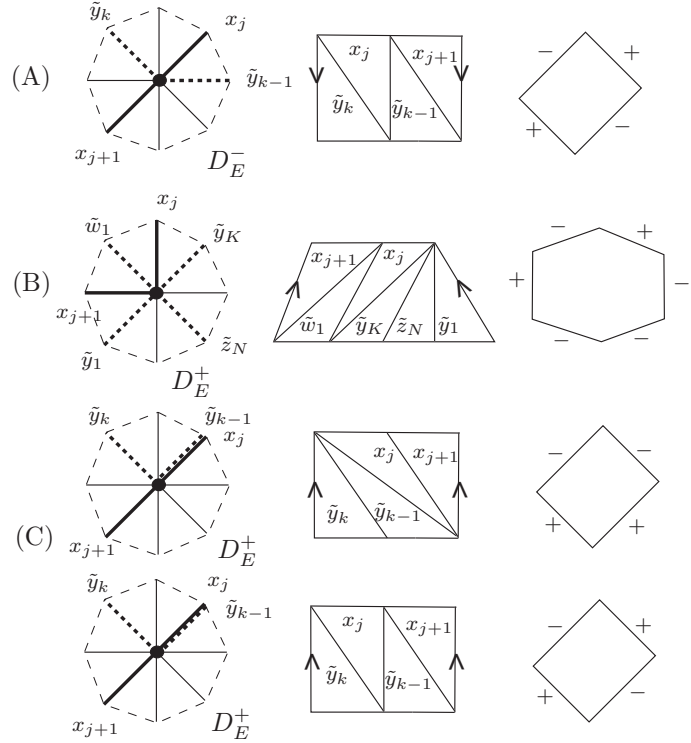} 
\caption{Four examples of the arrangement of edges of frames about the edge $E$ of $\T^*$ and unidentified faces of tetrahedra about $E$. Part (C) demonstrates the arrangement depending on choices made when we inflate at a face of $\T^*$.} 
\label{f-band-examples}
\end{center}
\end{figure}

We have labeled the unidentified faces by using the label of the
associated edge or virtual edge; of course, the role of edge or
virtual edge changes depending on being at $D_E^+$ or $D_E^-$. An
edge in a frame locally about $E$ may also be a virtual edge; this
creates a potential ambiguity but can easily be resolved. If an edge
of a frame about $E$ is also a virtual edge, then there is a face
$\sigma$ of $\T^*$ containing $E$ along with both of these edges in
the frame; so, from our earlier considerations in the inflation of
the face $\sigma$, we made a choice of diagonals in the subdivision
of the pyramid or prism we added. This choice of diagonals
determines the order about the edge $E$ of the faces of the
tetrahedron added for the edge versus that added for the virtual
edge. The possibilities are demonstrated in Part (C) of Figure
\ref{f-band-examples} where the arrangement of faces about the edge
$E$ incorporates the choice of diagonal in our subdivision of a
pyramid in the inflation of the face $\sigma$.

Our algorithm needs information on the arrangement of edges and
virtual edges about an edge $E$. This can be recorded at either end
of $E$ as represented in $D_E^+$ or $D_E^-$; we only need to be
consistent and stay with one choice or the other at the edge $E$.
This enables us to code the arrangement of edges and virtual edges
about the edge $E$ as follows. Make a choice of either $D_E^+$ or
$D_E^-$, then using a planar polygon having the number of sides as
there are edges and virtual edges of the frames about $E$, label its
boundary edges either ``$+$" or ``$-$", where we use ``$+$" to
correspond to an edge of a frame and ``$-$" to correspond to a
virtual edge of a frame. In Figure \ref{f-band-at-E}, we show this
method of coding the arrangement of edges about $E$; later, we shall
need to add information about the transverse directions which can
also be recorded from $D_E^+$ or $D_E^-$. These polygons, along with
recorded information on transverse directions will be called {\it
configuration polygons}.

In Figure \ref{f-band-examples} we provide three examples; in
Example (A) there are $4$ edges of the frames meeting $E$ (at
$D_E^-$ there are two edges and two virtual edges), in Example (B)
there are $6$ edges (at $D_E^+$ there are $2$ edges and $4$ virtual
edges), and (C) is an example where an edge is also a virtual edge
(at $D_E^+$ there are $2$ edges and $2$ virtual edges). In this last
example we provide two possibilities for the order of the edge and
virtual edge ( which determines the choice of diagonal in the
inflation of the face containing the two edges). Below we catalog
all possibilities about edges of $\T^*$ in the case of one ideal
vertex of index one (a cusped manifold with one cusp).

\vspace{.125 in}\noindent{\it Basic inflations at an edge of
$\T^*$}. There are three basic configurations of unidentified faces
about an edge of $\T^*$ and hence, three basic constructions for the
inflation at an edge. All other configurations are decomposed into a
combination of these three; hence, all other inflation constructions
at an edge of $\T^*$ are a combination of these basic constructions.
We give an algorithm in Lemma \ref{divide-poly} to decompose an
arbitrary configuration about an edge to a combination of the basic
configurations.

\begin{figure}[htbp]

            \psfrag{j}{\footnotesize{$x_j$}}\psfrag{J}{\footnotesize{$x_{j+1}$}}
\psfrag{3}{\footnotesize{$3$}}\psfrag{2}{\footnotesize{$2$}}\psfrag{x}{$(x_j)$}
\psfrag{X}{$(x_{j+1})$}
            \psfrag{+}{\footnotesize{$+$}}\psfrag{0}{\footnotesize{$0$}}
            \psfrag{1}{\footnotesize{$1$}}
            \psfrag{-}{\footnotesize{$-$}}
            \psfrag{D}{$D_E^-$}\psfrag{E}{\small{$E$}}\psfrag{v}{\small $v^*$}
            \psfrag{e}{$(x_j)(013)\leftrightarrow
(x_{j+1})(012)$}\psfrag{S}{Generic: no tetrahedra added}
        \vspace{0 in}
        \begin{center}
 
\includegraphics[width= 4 in]{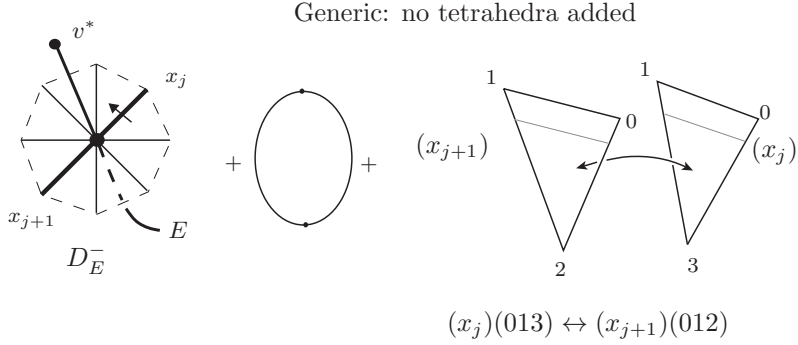} 
\caption{In a generic configuration the edge $E$ meets the frame in only one point, an index $2$ vertex. Free faces can be identified with no tetrahedra added.} 
\label{f-generic-config}
\end{center}
\end{figure}

\vspace{.125 in}\noindent{\bf Generic.} In the {\it generic
configuration} the edge $E$ only meets the frames in a single point
that is a vertex of index $2$ in the frame. Suppose notation is such
that $E$ meets the frame in $E^-$, which is the vertex $x_j^1 =
x_{j+1}^0$ between the edges $x_j$ and $x_{j+1}$ of the frame $\xi$.
Then the two unidentified faces about the edge $E$ are $(x_j)(013)$
and $(x_{j+1})(012)$. In this situation, we make the face
identification $(x_j)(013)\leftrightarrow (x_{j+1})(012)$. Note that
our convention for labeling the vertices of the tetrahedra $(x_j)$
and $(x_{j+1})$ and the identifications when inflating at a face of
$\T^*$ now have the free edges $(x_j)(01)\leftrightarrow
(x_{j+1})(01)$ identified.

\vspace{.125 in}\noindent{\bf Crossing.} In the {\it crossing
configuration} the edge $E$  meets the frames in two points, each is
a vertex of index $2$ in the frames (here the frames may be
different). Suppose notation is such that $E$ meets the frames at
$E^-$ in the vertex $y_{k-1}^1 = y_k^0$ between the edges $y_{k-1}$
and $y_k$ and $E$ meets the frames at $E^+$ in  the vertex $x_j^1 =
x_{j+1}^0$  between the edges $x_j$ and $x_{j+1}$. Then there are
four unidentified faces about the edge $E$: $(x_j)(013)$,
$(y_{k-1})(013)$, $(x_{j+1})(012)$, and $(y_{k})(012)$. In $D_E^-$,
the edge formed by $y_{k-1}$ and $y_k$ crosses the edge formed by
the virtual edges $\td{x}_j$ and $\td{x}_{j+1}$; hence, the name
crossing. It is possible that the edge $E$ meets the frames in two
points where each is a vertex of index $2$ in the frames but the
edges and virtual edges in this case do not cross. If this is the
situation, we do not have a crossing and will see below that we can
decompose this into two distinct generic configurations.

\begin{figure}[htbp]

            \psfrag{j}{\footnotesize{$\td{x}_j$}}\psfrag{J}{\footnotesize{$\td{x}_{j+1}$}}
\psfrag{3}{\footnotesize{$3$}}\psfrag{2}{\footnotesize{$2$}}
\psfrag{y}{\footnotesize{$(y_{k-1})(013)$}}\psfrag{Y}{\footnotesize{$(y_k)(012)$}}
\psfrag{x}{\footnotesize{$(x_{j})(013)$}}\psfrag{X}{\footnotesize{$(x_{j+1})(012)$}}
            \psfrag{+}{\footnotesize{$+$}}\psfrag{0}{\footnotesize{$0$}}
            \psfrag{K}{\footnotesize{$y_k$}}\psfrag{k}{\footnotesize{$y_{k-1}$}}
            \psfrag{1}{\footnotesize{$1$}}\psfrag{2}{\footnotesize{$2$}}
            \psfrag{3}{\footnotesize{$3$}}\psfrag{c}{\footnotesize{$(c)$}}
            \psfrag{-}{\footnotesize{$-$}}\psfrag{+}{\footnotesize{$+$}}
            \psfrag{D}{$D_E^-$}\psfrag{E}{\small $E$}\psfrag{v}{\small $v^*$}
            \psfrag{e}{\begin{tabular}{c}
$(x_j)(013)\leftrightarrow (c)(023);$\\
\\
$(c)(021)\leftrightarrow (x_{j+1})(012);$\\
\\
$(y_{k-1})(013)\leftrightarrow (c)(132);$\\
\\
$(c)(130)\leftrightarrow (y_{k})(012)$\\
\end{tabular}}\psfrag{S}{Crossing: one tetrahedron added}
        \vspace{0 in}
        \begin{center}
\includegraphics[width= 4.5 in]{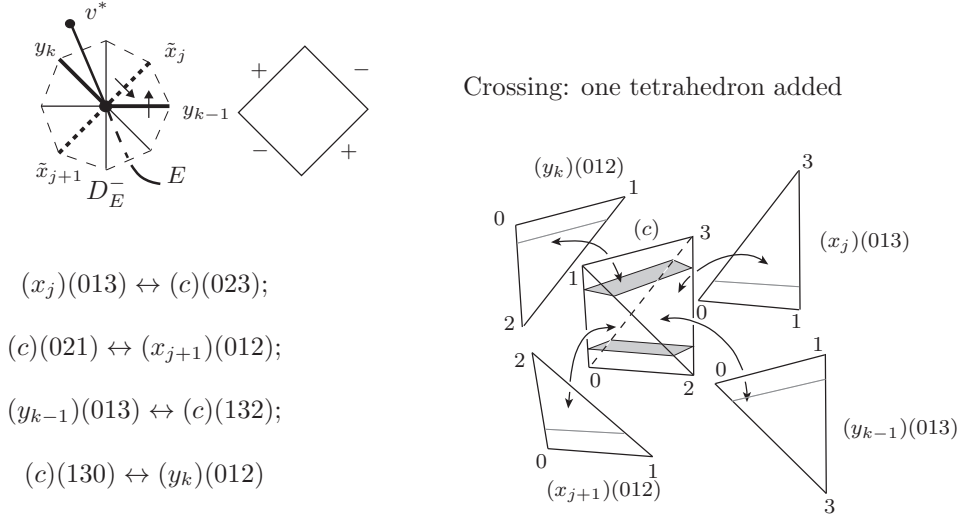} 
\caption{At a crossing the edge $E$ meets the frames in two points, each an index $2$ vertex and the edges and virtual edges form a crossing. A new tetrahedron is added.} 
\label{f-crossing-config}
\end{center}
\end{figure}

In the case of a crossing, it is necessary that we add a new
tetrahedron. We shall denote this tetrahedron by $(c)$ and its
vertices by $0,1,2,3$; we think of it as the join $(c)(02)\ast
(c)(13)$. If there is a crossing at the edge $E$, we shall use the
convention that the edge $(c)(02)$ is always associated with the
vertex at $E^+$ and then the edge $(c)(13)$ is always associated
with the vertex $E^-$; that is, since $x_j$ and $x_{j+1}$ are edges
at $E^+$, then  $(x_j)(01)=(c)(02)=(x_{j+1})(01)$ and since
$y_{k-1}$ and $y_k$ are edges at $E^-$, then
$(y_{k-1})(01)=(c)(13)=(y_k)(01)$. Using these conventions, then the
transverse directions in $D_E^-$ (or $D_E^+$) determine the
remaining vertices for the face identifications of the free faces
about the edge $E$.  Hence, we have the face identification
$(x_j)(013)\leftrightarrow (c)(023)$; $(c)(021)\leftrightarrow
(x_{j+1})(012)$; $(y_{k-1})(013)\leftrightarrow (c)(132)$; and
$(c)(130)\leftrightarrow (y_{k})(012)$. See Figure
\ref{f-crossing-config}. The scheme is that the free edge associated
with the branch containing $x_j$ and $x_{j+1}$ match through the
edge $(c)(02)$ and the free edge associated with the branch
containing $y_{k-1}$ and $y_k$ match through the edge $(c)(13)$; of
course, it is possible these are all the same branch and the same
free edge.

\begin{remark} Again, for later reference, we note that during inflation at a
crossing, we add a tetrahedron and  can modify the vertex-linking
surface in $\T^*$ by adding two quadrilaterals, one about each of
the added free edges. See Figure \ref{f-crossing-config} which shows
the added quads.\end{remark}

\vspace{.125 in}\noindent{\bf Branch.} In the {\it branching
configuration} the edge $E$  meets a frame in one point that is a
branch point of index $b, b\ge 3$. Suppose notation is such that $E$
meets the frame at $E^-$; hence, $E^-$ is a branch point of index
$b$ of a frame. In $D_E^-$ there are $b$ edges in the frame, which
we may assume are ordered cyclically as $x, y, \ldots, z$. Note, it
is possible that two edges may be from the same branch of the frame
and it is possible that $E^-$ is an initial point for some branches
and a terminal point for others. Furthermore, here  we need to use
the transverse directions given to the branches of the frames. There
are $b$ unidentified faces of tetrahedra, which were added to the
triangulation $\T^*$ from the $b$ branches at $E^-$, about the edge
$E$; our conventions give labels to these unidentified faces as
$(x)(01{\varepsilon}_{x})$, $(y)(01\varepsilon_{y})$, \ldots,
 and $(z)(01\varepsilon_{z})$, where
$\varepsilon_e = 2$ or $3$, depending on the direction induced on
the edge $e$ by the directed branch of the frame containing $e$. In
this situation, we let $P_b$ be a planar $b$-gon and form the cone
$(b^*)=0\ast P_b$ over the $b$-gon $P_b$. We label the vertices of
$(b^*)$ as $0,1,2,\ldots,b$, where $1,2,\ldots, b$ are the vertices
of $P_b$ ordered counter-clockwise from the view at the cone point.
The cone $(b^*)$ has $b$ triangular faces:
$(b^*)(120),(b^*)(230),\ldots,(b^*)(b10)$. In making face
identifications we must take into consideration the induced
transverse directions on the edges meeting $E^-$. For demonstration,
see Figure \ref{f-crossing-config},  we assume the edges are
cyclically ordered in $D_E^-$ as $x_1, y_{\tiny{K}}, x_{\tiny{J}},
z_1,\ldots,$. Then we have face identifications
$(x_1)(012)\leftrightarrow (b^*)(b10)$; $(y_K)(013)\leftrightarrow
(b^*)(210)$; $(x_J)(013)\leftrightarrow (b^*)(320)$;
$(z_1)(012)\leftrightarrow (b^*)(340)$, $\ldots$. The cone point of
$(b^*)$ is identified with the ideal vertex at $E^+$.

Note that we may start the cyclical ordering in $D_E^-$ at any edge;
we then start the face labeling of $b^*$ so that we start
$(b^*)(b1)$ at the first edge in our cyclic ordering; it is after
this, we need to follow the transverse directions in making face
identification. If $e=e_1$ is an initial edge of a branch, then the
identification is $(e_1)(012)\leftrightarrow (b^*)(n(n+1)0)$ in the
case of an orientable surface; if $e_N$ is a terminal edge of a
branch, then the identification is $(e_N)(013)\leftrightarrow
(b^*)((n+1)n0)$ in the case of an orientable surface.

\begin{figure}[htbp]

            \psfrag{x}{\footnotesize $(x_1)$}\psfrag{X}{\footnotesize $(x_{\tiny{J}})$}
            \psfrag{Y}{\footnotesize $(y_{\tiny{K}})$}
            \psfrag{z}{\footnotesize $(z_1)$}
\psfrag{J}{\footnotesize $x_{\tiny{J}}$} \psfrag{K}{\footnotesize
$y_{\tiny{K}}$} \psfrag{n}{\footnotesize $z_1$}
\psfrag{j}{\footnotesize $x_1$}
            \psfrag{0}{\footnotesize{$0$}}
            \psfrag{1}{\footnotesize{$1$}}\psfrag{a}{\footnotesize{$2$}}
            \psfrag{c}{\footnotesize{$3$}}\psfrag{b}{\footnotesize{$b$}}
            \psfrag{v}{\footnotesize{$v*$}}
            \psfrag{-}{\footnotesize{$-$}}\psfrag{B}{$(b^*)=0\ast P_b$}
            \psfrag{P}{$P_b$}\psfrag{E}{\small$E$}
            \psfrag{D}{$D_E^-$}
            \psfrag{e}{\begin{tabular}{c}
$(x_1)(012)\leftrightarrow (b^*)(b10);$\\
\\
$(y_K)(013)\leftrightarrow (b^*)(210);$\\
\\
$(x_J)(013)\leftrightarrow (b^*)(320);$\\
\\
$(z_1)(012)\leftrightarrow (b^*)(340)$\\
$\cdots$\\
\end{tabular}}\psfrag{S}{Branch: $(b-2)$ tetrahedra added}

        \vspace{.5 in}

 \begin{center}
\includegraphics[width= 4 in]{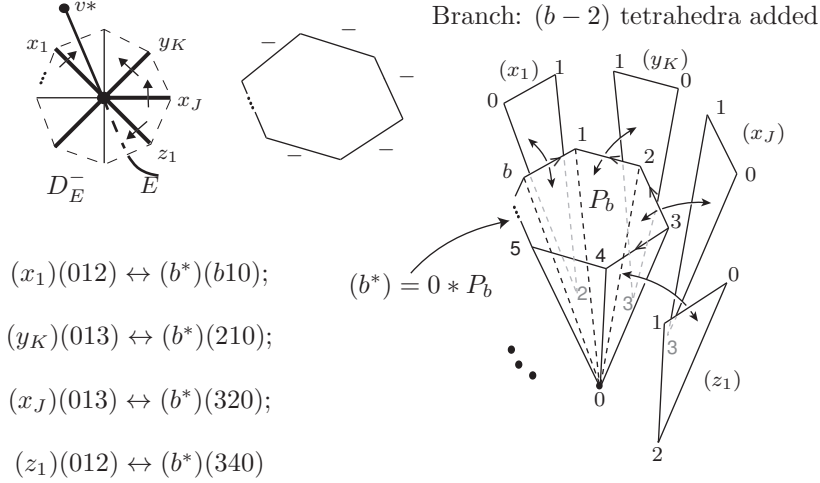} 
\caption{The edge $E$ meets the frames in one point, an index $b$ branch point ($b\ge 3$). The cone over a planar $b$-gon is added and then is subdivided into $(b-2)$ tetrahedra, without adding vertices.}
\label{f-branch-config}
\end{center}
\end{figure}

\begin{remark} In the case of an inflation at a branch
configuration, we have the possibility of subdividing the planar
$b$--gon in a number of ways; hence, we see that even for a fixed
frame, the inflation construction does not lead to a unique
triangulation.

Also, for later reference, we note that during inflation at a branch
point, in each added tetrahedron we modify the vertex-linking
surface in $\T^*$ by adding two parallel triangles, one that is
vertex-linking at the vertex opposite the free face and one that is
of the same normal triangle type but we refer to it as parallel to
the free face. Also, the faces of subdivided polygon make up faces in the boundary of the manifold resulting from inflation. 
\end{remark}

We now have the various pieces in place to put together our
inflation construction. The following lemma provides the algorithm
to reduce any configuration of frames about an edge of the ideal
triangulation $\T^*$  to a composition of the three basic
configurations given above. We assume we are given a configuration
via a planar polygon with its edges marked by either a ``$+$" or a
``$-$". We shall call a maximal sequence of adjacent edges having
the same mark a {\it link}.

\begin{lem}\label{divide-poly} Suppose $P$ is a planar polygon with edges marked either
``$+$" or ``$-$" and there is not just one edge of a given mark.
Then $P$ can be subdivided into sub-polygons where each sub-polygon
has edges marked with either a ``$+$" or ``$-$", extending the
markings on the boundary of $P$, and each sub-polygon has one of the
three forms given above; generic, crossing, or branch.\end{lem}

\begin{figure}[htbp]

            \psfrag{b}{\scriptsize{branch}}\psfrag{A}{(A)}
\psfrag{B}{(B)}\psfrag{C}{(C)}
            \psfrag{g}{\scriptsize{generic}}
            \psfrag{c}{\scriptsize{crossing}}
        \psfrag{2}{\small{add $2$ tetrahedra}}
            \psfrag{3}{\small{add $3$ tetrahedra}} \psfrag{7}{\small{add $7$ tetrahedra}}
            \psfrag{v}{\scriptsize{$v$}}
            \psfrag{-}{\scriptsize{$-$}}\psfrag{+}{\scriptsize{$+$}}
            \psfrag{D}{\small{$D_E^-$}}

        \vspace{0 in}
        \begin{center}
\includegraphics[width= 4.5 in]{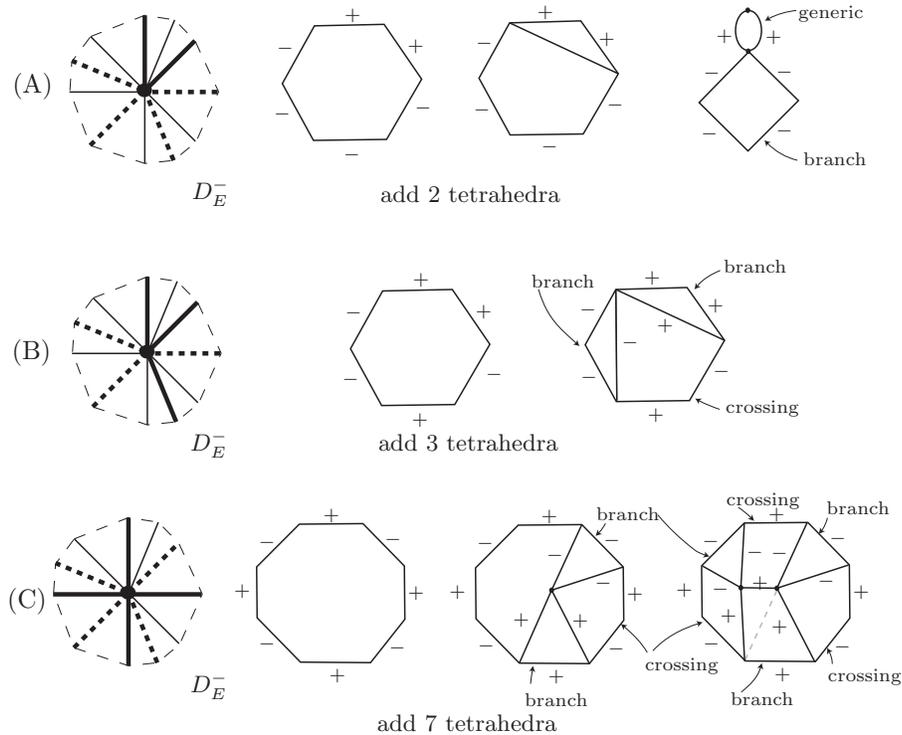} 
\caption{Three examples giving the steps in the algorithm to subdivide a marked planar polygon into marked sub-polygons having generic, crossing, or branch configurations.} 
\label{f-subdivide-poly}
\end{center}
\end{figure}

\begin{proof} The proof is by induction on the number of edges of
the polygon $P$.

By hypothesis there is not just one edge of a given mark; hence, the
induction begins with $P$ a bi-gon and the configuration is generic.

We assume the conclusion is true for any marked polygon satisfying
the hypothesis and having fewer than $p$ edges, $p\ge 3$. Now,
suppose $P$ is a marked $p$-gon.

Consider the various links (maximal collections of adjacent edges
having the same mark) in the boundary of $P$.

If there is only one link, then we have a branch configuration.

If there are no links of length at least $2$, then the markings, as
we proceed about the edges of $P$, alternate. Since $P$ must have an
even number of edges (by alternating markings) and $p\ge 3$, we have
that $P$ must have at least $4$ edges. If $P$ has precisely $4$
edges, then its configuration is a crossing. So, we may assume $P$
has at least $6$ edges. Hence, there are edges of the boundary of
$P$ marked with a ``$+$" and ``$-$" and which are separated by at
least two other edges; we denote these edges $e^+$ and $e^-$,
respectively. Let $v$ be a point in the interior of $P$ and consider
the two triangles $v\ast e^+$ and $v\ast e^-$ formed by the join of
the edges $e^+$ and $e^-$ with $v$. Mark the edges of $v\ast e^+$
with a ``$+$" extending the mark on $e^+$ and mark the edges $v\ast
e^-$ with a ``$-$" extending the mark on $e^-$. This subdivides $P$
into four sub-polygons each with marked edges that satisfy our
hypothesis and each having fewer than $p$ edges. In fact the
triangles $v\ast e^+$ and $v\ast e^-$ are marked with a branch
configuration. The other two sub-polygons have alternating markings.
See Figure \ref{f-subdivide-poly}.

So, we may assume there are links having length at least two. If
there are just two such links, then we can pinch $P$ along a segment
interior to $P$ that separates the two links in its boundary. This
gives two polygons; these polygons have either generic or branch
configurations. Hence, we may assume there are links having length
at least two and together these do not determine the totality of
marks on the boundary of the configuration polygon.

Now, for a link having length at least two, draw a straight line
through the interior of $P$ from one of its end points to the other.
Mark this edge with the same mark as the link that determined it. In
this case we have no ambiguity since we have more than two links in
the boundary. Repeating this for links in the boundary, we subdivide
the polygon $P$ into a number of sub-polygons, all but one having
the configuration of a branch and the exceptional marked sub-polygon
having alternating marks and it  at least $4$ edges. Hence, the
induction step follows.\end{proof}

We exhibit some of the steps in the algorithm from Lemma
\ref{divide-poly} in the following three examples.

In Example (A) there are only two links in the boundary of the
marked polygon. The link of length $2$ in this case leads to a
generic configuration. The link of length $4$ leads to a branch
configuration, requiring the addition of $2$ tetrahedra.

In Example (B), there are two links having length at least two. We
draw line segments in the polygon separating off these links and
mark the new edges. This gives two branches, each of index $3$, and
a new polygon. The new polygon has a crossing configuration. Each
branch of index $3$ requires that we add a tetrahedron and the
crossing configuration requires that we add a tetrahedron; hence, we
add $3$ tetrahedra.

In Example (C), each link in the boundary has length one, the
markings are alternating. Hence, we must add a vertex to the
interior of the polygon and cone on two edges in the boundary that
are length at least two apart and have distinct markings. We mark
the new edges, getting two branch configurations of index $3$ and
two new marked polygons. In this example one of these is a crossing
configuration. The other has alternating markings and so we must
again add a vertex to its interior. The result is two more branch
configurations of index $3$ and a crossing.  We must add $7$
tetrahedra to resolve the configuration of Example (C).

 Face identifications for any of the resulting
decompositions of the configuration polygons are induced by edges
and transverse directions in the boundary of the configuration
polygon. A configuration is either generic, a crossing, branched or
is decomposed into a collection of generic, branch and crossing
configurations. A branch configuration with edges labeled ``$+$" has
its cone vertex identified with the ideal vertex of $E$ at the end
labeled $E^-$ and all of its other vertices at the ideal vertex of
$E$ labeled $E^+$. Each crossing configuration adds a tetrahedron
with two of its vertices identified to  the ideal vertex at one end
of $E$ and the remaining two vertices identified to the ideal vertex
at the other end of $E$. If $c = (c)(02)\ast(c)(13)$ is a crossing,
then our conventions have the vertices $(c)(0)$ and $(c)(2)$
identified to the ideal vertex at the end designated $E^+$; and,
then, $(c)(1)$ and $(c)(3)$ identified to the ideal vertex at the
end designated $E^-$. Every configuration polygon is associated with
a unique edge of $\T^*$.

\begin{figure}[htbp]

            \psfrag{b}{\footnotesize{branch}}
            \psfrag{g}{\footnotesize{generic}}
            \psfrag{c}{\footnotesize{crossing}}
            \psfrag{-}{\footnotesize{$-$}}\psfrag{+}{\footnotesize{$+$}}
            \psfrag{D}{$D_E^+$}\psfrag{d}{$D_E^-$}
             \psfrag{1}{\begin{tabular}{c} index $2$ vertex\\
no tetrahedron added\\
\end{tabular}}
 \psfrag{2}{\begin{tabular}{c} index $3$ branch\\
add one tetrahedron\\
\end{tabular}}
 \psfrag{3}{\begin{tabular}{c} index $4$ branch\\
add two tetrahedra\\
\end{tabular}}
 \vspace{.5  in}
        \begin{center}
\includegraphics[width= 4 in]{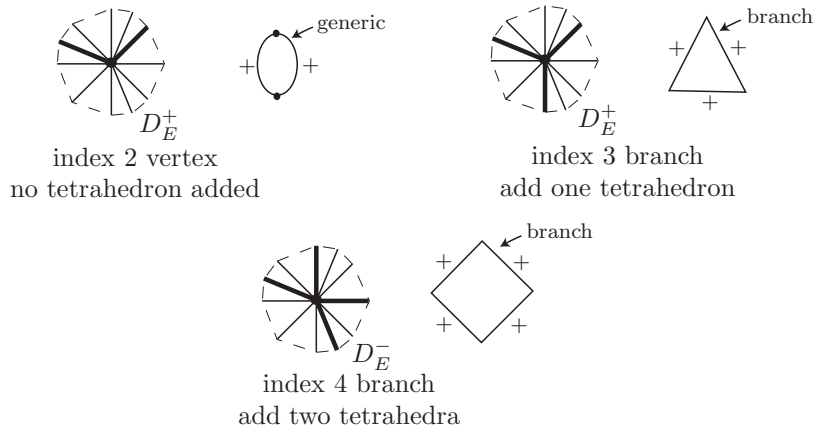} 
\caption{The edge $E$ meets the frame in one point.} 
\label{f-config-poly-1pt}
\end{center}
\end{figure}

\begin{figure}[htbp]

            \psfrag{b}{\footnotesize{branch}}\psfrag{A}{(A)}
\psfrag{B}{(B)}\psfrag{C}{(C)}
            \psfrag{g}{\footnotesize{generic}}
            \psfrag{c}{\footnotesize{crossing}}
            \psfrag{v}{\footnotesize{$v$}}
            \psfrag{-}{\footnotesize{$-$}}\psfrag{+}{\footnotesize{$+$}}
            \psfrag{D}{$D_E^+$}\psfrag{d}{$D_E^-$}
             \psfrag{4}{\begin{tabular}{c} index $2$ over index $2$\\
no crossing\\
no tetrahedron added\\
\end{tabular}}
 \psfrag{5}{\begin{tabular}{c} index $2$ over index $2$\\
crossing\\
\end{tabular}}
 \psfrag{6}{\begin{tabular}{c} index $3$ branch\\
over index $2$\\no crossing\\
\end{tabular}}
\psfrag{7}{\begin{tabular}{c} index $3$ branch\\
over index $2$\\crossing\\
\end{tabular}}
        \begin{center}
\includegraphics[width= 5 in]{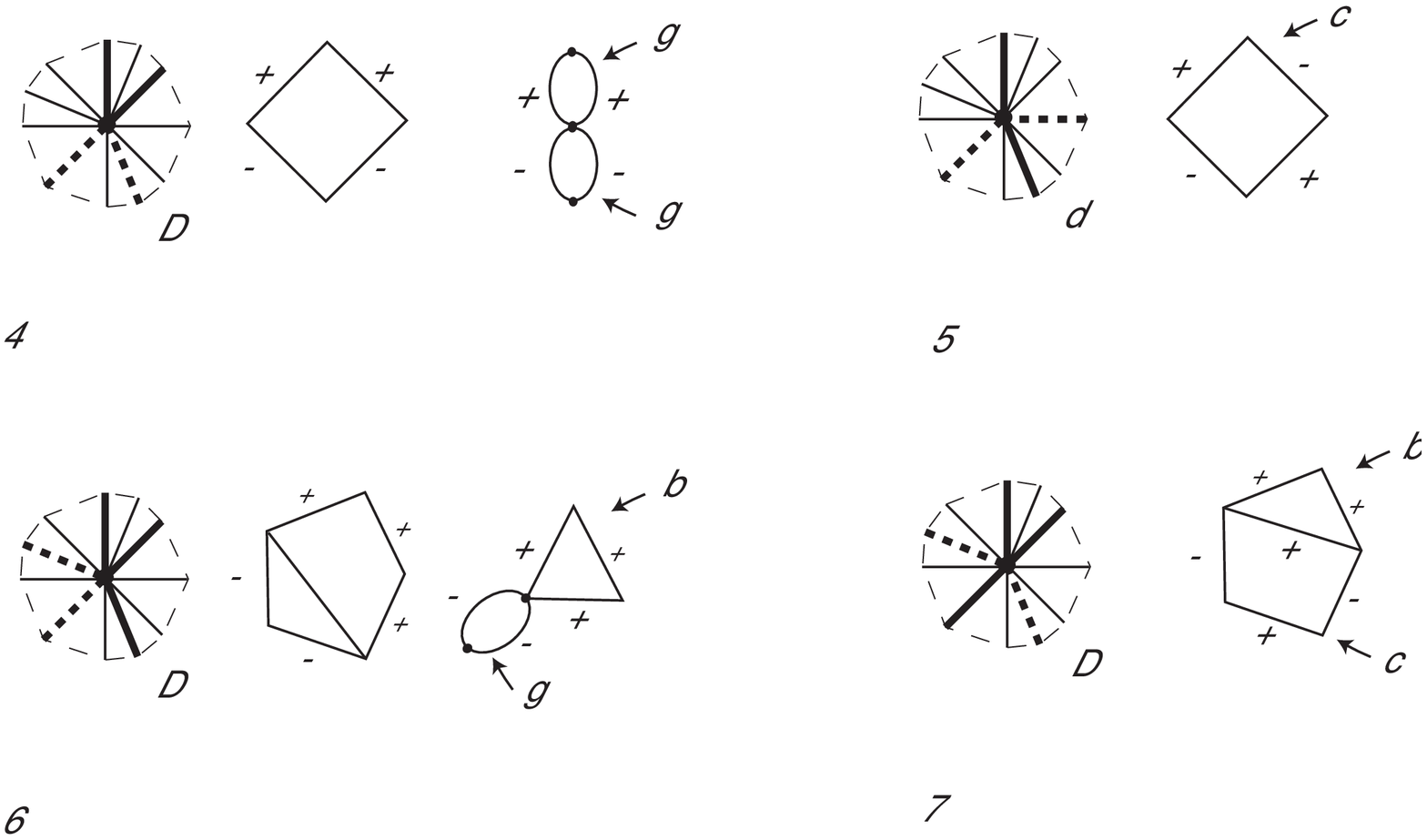} 
\caption{The edge $E$ meets the frame in two points; Part I.}
\label{f-config-poly-2pt-4-5}
\end{center}
\end{figure}

\begin{figure}[htbp]

            \psfrag{b}{\footnotesize{branch}}\psfrag{A}{(A)}
\psfrag{B}{(B)}\psfrag{C}{(C)}
            \psfrag{g}{\footnotesize{generic}}
            \psfrag{c}{\footnotesize{crossing}}
            \psfrag{v}{\footnotesize{$v$}}
            \psfrag{-}{\footnotesize{$-$}}\psfrag{+}{\footnotesize{$+$}}
            \psfrag{D}{$D_E^+$}\psfrag{d}{$D_E^-$}
 \psfrag{1}{\begin{tabular}{c} index $4$ branch\\ over index $2$\\
no crossing\\
\end{tabular}}
 \psfrag{2}{\begin{tabular}{c} index $4$ branch\\
over index $2$\\crossing\\
\end{tabular}}
\psfrag{3}{\begin{tabular}{c} index $4$ branch\\
over index $2$\\crossing\\
\end{tabular}}
\psfrag{4}{\begin{tabular}{c} index $3$ branch\\over index $3$ branch\\
no crossing\\
\end{tabular}}
 \psfrag{5}{\begin{tabular}{c} index $3$ branch\\
over index $3$ branch\\one crossing\\
\end{tabular}}
\psfrag{6}{\begin{tabular}{c} index $3$ branch\\
over index $3$ branch\\ two crossings\\
\end{tabular}}
        
        \begin{center}
\includegraphics[width= 5 in]{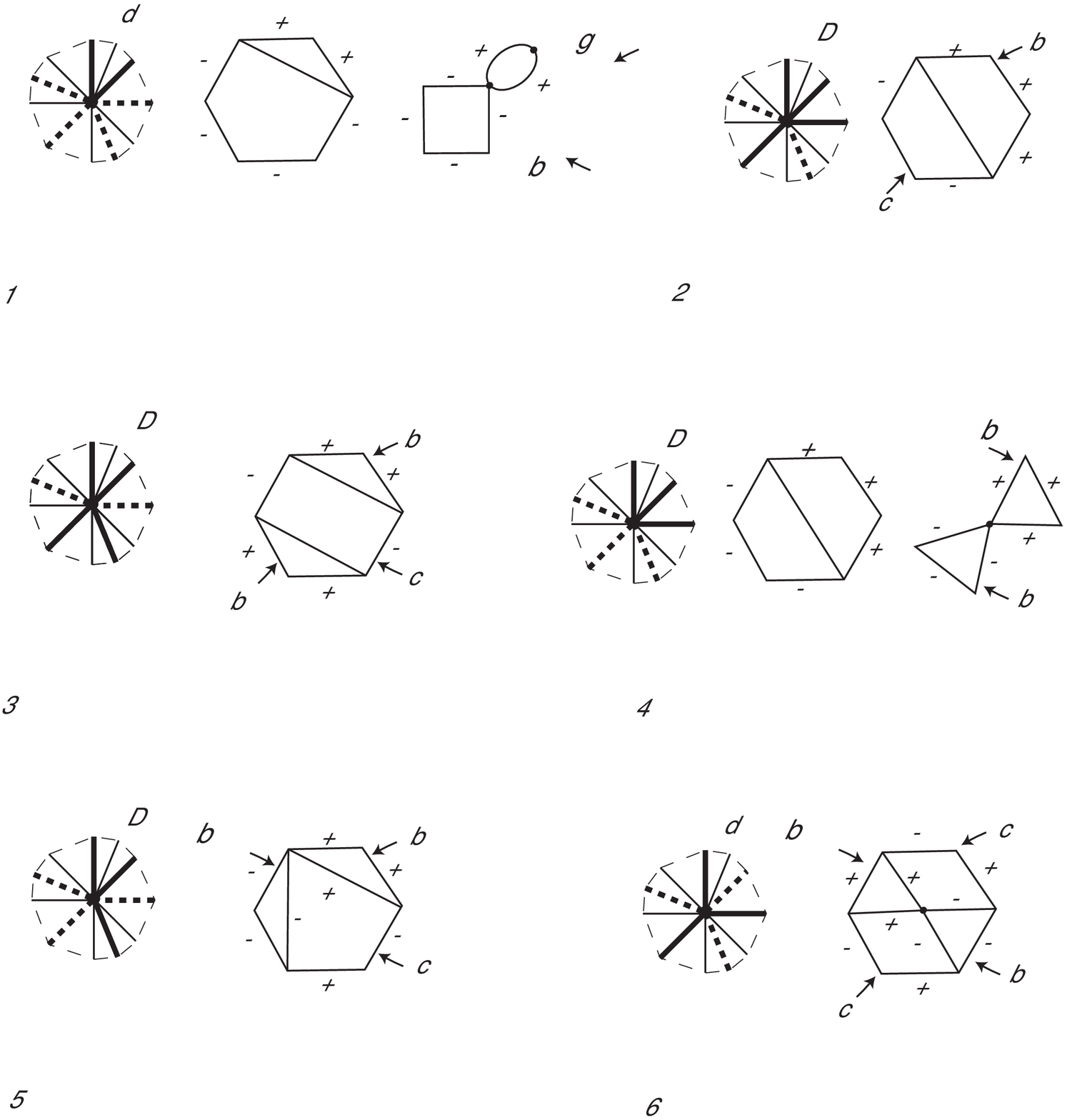} 
\caption{The edge $E$ meets the frame in two points; Part II.}
\label{f-config-poly-2pt-4-6}
\end{center}
\end{figure}

\vspace{.125 in}\noindent{\it Demonstration.} As a demonstration of
what we might have for the configuration polygons, we catalog in
Figures \ref{f-config-poly-1pt}-\ref{f-config-poly-2pt-4-6} all the
possibilities for an ideal triangulation of the interior of a
knot-manifold (one ideal vertex of index one).

\vspace{.25 in}

In the previous construction, we start with an ideal triangulation
$\T^*$ of the interior of a compact $3$--manifold $X$. We select a
collection $\Lambda$ of frames in the vertex-linking surfaces of
$\T^*$. Using these frames as guides, we add tetrahedra to the
collection of tetrahedra in $\T^*$; we discard some of the face
identifications in $\T^*$ and add new face identifications. The new
collection of tetrahedra with some of the face identifications from
$\T^*$ and some new face identifications give us a triangulation
$\T_\Lambda$ of a $3$--complex. Below in Theorem \ref{thm:inflate}
we show that the underlying point set of the 3-complex $\T_\Lambda$
is a compact 3--manifold with boundary, homeomorphic with $X$, and
$\T_\Lambda$ is an inflation of $\T^*$. First, we describe two
important surfaces constructed during the inflation.

\subsection{Boundary and normal boundary in $\T_\Lambda$} For each
ideal vertex $v_j^*$ of $\T^*$ we have an associated vertex-linking
surface $V_j^*$ and a frame $\xi_j$ in the 1--skeleton of the
induced triangulation on $V_j^*$, which is in the collection
$\Lambda$.

\vspace{.125 in}\noindent{\it Boundary surfaces.} The inflation
along $\xi_j$ creates what we are calling free edges and free faces,
as was pointed out in the remarks above following various steps in
the inflation construction. Free edges are introduced during
inflations at the faces of $\T^*$ that contain an edge of $\xi_j$.
Other free edges and free faces are created when subdividing
branching polygons. All the free edges created along a branch of
$\xi_j$ are identified and are identified with boundary edges in the
branching polygons. The identification of free faces and free edges
is equivalent to pairwise identification of the boundary edges of
the branching polygons. For the inflation along the frame $\xi_j$,
these identifications give a connected surface, $B_j$, with a
one-vertex triangulation and $B_j$ is a subcomplex of $\T_\Lambda$.
It follows that $B_j$ is homeomorphic to $V^*_j$.  Since a face in
$B_j$ is contained in only one tetrahedron, it is a boundary face of
$\T_\Lambda$; below we will have that $B_j$ is a component of $\bdy
M$, where $M$ is the underlying point set of $\T_\Lambda$.

\vspace{.125 in}\noindent{\it Boundary-linking normal surfaces.} The
inflation along $\xi_j$ also leads to the construction of a
boundary-linking normal surface in $\T_\Lambda$, which can be viewed
as an inflation of the vertex-linking surface $V_j^*$. While the
significate elements of the inflation of the vertex-linking surface
$V_j^*$ are related to the inflation along the frame $\xi_j$, it is
possible that inflations along other frames in $\Lambda$ also
contribute to the inflation of $V_j^*$.

\begin{figure}[htbp]

            \psfrag{b}{$(b^*)$}\psfrag{P}{$P_b$}
\psfrag{v}{$v^*$}\psfrag{o}{\large {or}}

        \vspace{0  in}
        \begin{center}
\includegraphics[width= 2.5 in]{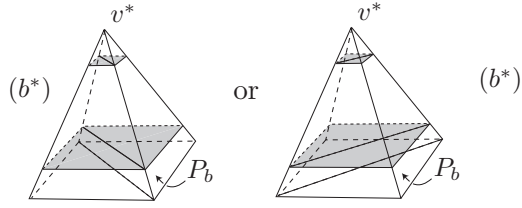}
\caption{Triangles added to the boundary-linking surface at a branch point.}
\label{fr-triangles-in-branch}
\end{center}
\end{figure}

At each edge of $\xi_j$, we have an inflation in a face of $\T^*$
which removes a face identification in the triangulation and adds a
tetrahedron and two new face identifications. The affect on the
vertex-linking surface $V_j^*$ is to remove edge identifications of
$V_j^*$ in the faces containing the edges of $\xi_j$ and for each
tetrahedron added in this fashion, we can add to $V_j^*$ a
quadrilateral around the added free edge. It also is necessary to
add two vertex-linking triangles at each end of the edge that is
opposite the free edge in the added tetrahedron; these
vertex-linking triangles may or may not be added to the
vertex-linking surface $V^*_j$ but do add to a vertex-linking
surface as part of its inflation. These additions are exhibited in
Figures \ref{f-blow-up-face1}, \ref{f-blow-up-face2}, and
\ref{f-blow-up-face3}. At a crossing, we add a tetrahedron and make
identifications between its faces and four unidentified faces of
previously added tetrahedra. As part of the inflation of the
vertex-linking surfaces, we add two quadrilaterals, which might be
added to distinct vertex-linking surfaces; these quadrilaterals are
added about the free edges of the crossing tetrahedron. See Figure
\ref{f-crossing-config}. At a branching, we have a cone over a
planar polygonal $b$--gon, $P_b$, which we denoted above
$(b^*)=0\ast P_b$. The cone faces of $(b^*)$ are then identified
with $b$ unidentified faces of previously added tetrahedra. We
subdivide $P_b$ into $b-2$ triangles and extend this to a
subdivision of $(b^*)$ into $b-2$ tetrahedra by coning each triangle
from the cone point (which as noted above is one of the ideal
vertices of $\T^*$). This gives a triangulation of $(b^*)$ with
$b-2$ free faces corresponding to the faces in the subdivided
$b$--gon, $P_b$. In this case we add two triangles in each of the
$b-2$ tetrahedra, one is a vertex-linking triangle about the cone
point and the other a face-linking triangle about the free face of
the tetrahedron, which is in $P_b$. We exhibit this is Figure
\ref{fr-triangles-in-branch} for a valence 4 branch point.

Hence, starting with a vertex-linking surface $V^*_j$, and by
removing various edge identifications in $V^*_j$ that are dictated
by the inflation construction and adding quadrilaterals and
triangles, we arrive at a closed normal surface in the triangulation
$\T_\Lambda$. We  denote this surface by $V_j$ to indicate its
relationship with $V^*_j$. The normal surface $V_j$ is the frontier
of a small neighborhood of the boundary subcomplex $B_j$ from above.
Further, we observe that the normal surface $V_j$ maps (``crushes")
via a cell-like map to the vertex-linking surface $V^*_j$;  and
therefore, the boundary-linking normal surface $V_j$ is homeomorphic
to the vertex-linking surface $V^*_j$.

\begin{figure}[htbp]

            \psfrag{v}{\small $v^*$}
            \psfrag{w}{\small $w^*$}
            \psfrag{E}{\small $E$}\psfrag{+}{\footnotesize $E^+$}\psfrag{-}{\footnotesize $E^-$}\psfrag{1}{\footnotesize $\td{1}$}\psfrag{0}{\footnotesize $1$}
\psfrag{2}{\footnotesize $2$}
 \psfrag{3}{\footnotesize $\td{2}$}
 \psfrag{4}{\footnotesize $4$}
 \psfrag{5}{\footnotesize $\td{4}$}
\psfrag{6}{\footnotesize $6$} \psfrag{7}{\footnotesize $\td{6}$}
\psfrag{8}{\footnotesize $8$} \psfrag{9}{\footnotesize $\td{8}$}

\psfrag{G}{generic point}\psfrag{C}{crossing point}\psfrag{B}{branch
point}

        \begin{center}
\includegraphics[width= 4.25 in]{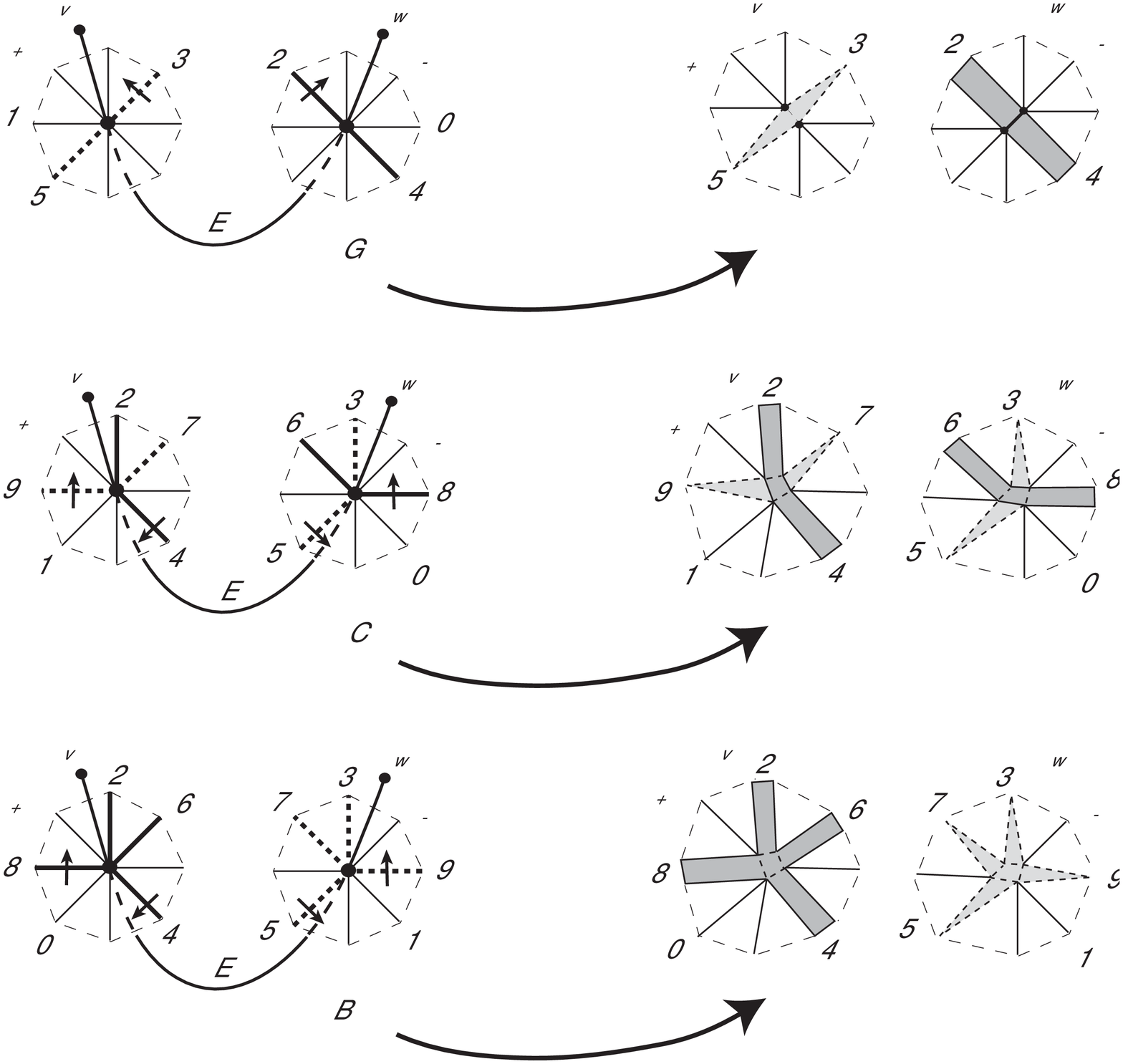}
\caption{The local view of the inflation induced on a vertex-linking surface at a generic, crossing, and branch point of the frame.} 
\label{f-local-inflate-v-surface}
\end{center}
\end{figure}

In Figure \ref{f-local-inflate-v-surface}, we exhibit local results
of the inflation of the vertex-linking surfaces coming from
inflations at faces of $\T^*$ along with a contiguous inflation
along an edge of $\T^*$. We show the latter at a generic, a
crossing, and a branch point of the frame. An inflation in a face
creates an inflation in the vertex-linking surfaces at each vertex
of the face and an inflation along an edge $E$ of $\T^*$ creates an
inflation of the vertex-linking surfaces at both ends of the edge
$E$.

\vspace{.125 in}

\begin{thm}\label{thm:inflate} Suppose $X$ is a compact 3--manifold
with boundary, no component of which is a 2--sphere, and $\T^*$ is
an ideal triangulation of the interior of $X$. Then for any
collection of frames $\Lambda$, one frame in each of the
vertex-linking surfaces of $\T^*$, the underlying point set of
$\T_\Lambda$ is a compact $3$--manifold $M$ homeomorphic to $X$ and
the triangulation $\T_\Lambda$ of $M$ is an inflation of the ideal
triangulation $\T^*$.\end{thm}

\begin{proof} First we shall show that the construction gives a
minimal vertex triangulation of a compact $3$-manifold with
boundary.

Let $\T_\Lambda$ denote the triangulation constructed and let $M$ be
the underlying point set. $M $ is a $3$-manifold at the image of
each interior point of a tetrahedron and at each interior point of a
face. Having chosen a fixed direction on each edge of $\T^*$ and
respecting this orientation with face identifications in both the
inflations at the faces and at the edges, the construction gives a
$3$-manifold at the image of each edge having a complete cycle of
face identifications about it. The edges that do not have a complete
cycle of face identifications about them correspond to the ``free
edges" added when we inflate at a face of $\T^*$ or the ``free
edges" added in the subdivision of the $b$-gons, $b\ge 3$, for
branches in the configuration polygons. Starting at any one of the
free edges and following face identifications in both directions,
there are chains of face identifications that end only when we come
to a face of a $b$-gon in one of the configuration polygons. Hence,
an interior point $p$ of any of these edges has a $3$--cell
neighborhood with the point $p$ in the boundary of the $3$--cell.

It follows that $M $ is a $3$--manifold at each point, except
possibly at the vertices of $\T_\Lambda$. We shall show that $M $ is
a $3$--manifold by showing that each vertex-linking surface in
$\T_\Lambda$ is a disk.

For $v_j$ a vertex of $\T_\Lambda$, we have $v_j$ in $B_j$ and have
denoted the boundary-linking normal surface along $B_j$ by $V_j$. We
shall show that we can modify the boundary-linking normal surface
$V_j$ to get the vertex-linking surface about $v_j$.

The normal surface $V_j$ is the frontier of a small regular
neighborhood of $B_j$. Hence, we can catalog how $V_j$ meets the
tetrahedra of  $\T_\Lambda$.

For $\Delta$ a tetrahedron of $\T_\Lambda$, we have $V_j$ meets
$\Delta$ in only one of the following:
\begin{itemize}\item[-] a subset of vertex-linking triangles and $\Delta$ meets $B_j$ only in
those vertices where there is a vertex-linking triangle in
$V_j$,\item[-] an edge-linking quadrilateral about an edge $e$ of
$\Delta$ and a subset of vertex-linking triangles about the vertices
of the edge of $\Delta$ opposite $e$ and  $\Delta$ only meets $B_j$
in the edge $e$ and those vertices where $V_j$ has a vertex-linking
triangle, \item[-] two parallel edge-linking quadrilaterals and
$\Delta$ only meets $B_j$ in those edges about which $V_j$ has the
edge-linking quadrilaterals, or\item[-]  a face-linking triangle and
possibly a vertex-linking triangle at the opposite vertex and
$\Delta$ only meets $B_j$ in the face at which we have the
face-linking triangle and  at the opposite vertex and then only if
$V_j$ has a vertex-linking triangle at that vertex.\end{itemize}

We now show how to construct the vertex-linking surface at $v_j$
from the normal surface $V_j$. See Figure \ref{fr-replace-link}.
First, remove each quadrilateral, $Q$, in $V_j$ and replace it with
two vertex-linking disks, one at each end of the edge the
quadrilateral $Q$ was linking; see Figure \ref{fr-replace-link}(A).
Then remove each face-linking triangle, $\sigma$,  in $V_j$ and
replace it with three vertex-linking triangles, one at each vertex
of the face in $B_j$ that the triangle $\sigma$ in $V_j$ was
linking; see \ref{fr-replace-link}(B).

\begin{figure}[htbp]
             \psfrag{a}{\footnotesize{$Q$ in $V_j$}}\psfrag{c}{\footnotesize{$\sigma$ in $V_j$}}\psfrag{b}{\footnotesize in $B_j$}\psfrag{A}{\begin{tabular}{c} (A) quad replaced \\
by two triangles\\
\end{tabular}}
 \psfrag{B}{\begin{tabular}{c} (B) triangle replaced\\
  by three triangles\\
\end{tabular}}

 \vspace{0 in}
        \begin{center}
\includegraphics[width= 3.5 in]{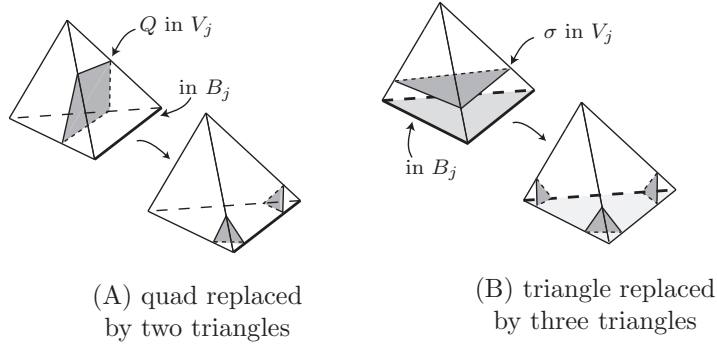} 
\caption{In (A) a quad in $V_j$ is replaced by two triangles. In (B) a boundary parallel triangle in $V_j$ is replaced by three triangles.}
\label{fr-replace-link}
\end{center}
\end{figure}

As observed above, the surfaces $V^*_j, B_j$, and $V_j$ are all
homeomorphic. So, if they have genus $g$, then each has Euler
characteristic, $\chi = 2-2g$. Starting with $V_j$, removing an
edge-linking quadrilateral and adding two vertex-linking disks does
not change the Euler characteristic, if we do not count the new
vertices added to the edge in $B_j$.  Removing a face-linking
triangle and replacing it with three vertex-linking triangles
changes the Euler characteristic by $-\frac{5}{2}$ if we again do
not count the new vertices at the edges in $B_j$. Note in this case,
we have included in the count all edges, including those in the face
in $B_j$.

Thus the Euler characteristic of $V_j$ is modified by $-\frac{5}{2}$
for each face of $B_j$ and by $+2$ for each edge in $B_j$. Let $D_j$
denote the surface we get from these modifications of $V_j$. There
are $4g-2$ faces in $B_j$ and $6g-3$ edges. It follows that
$$\chi(D_j) = 2-2g -5(2g-1)+2(6g-3) = 1$$
and $D_j$ is a disk.

It follows that $M$, the underlying point set of $\T_\Lambda$ is a
compact $3$--manifold with boundary, $\T_\Lambda$ has normal
boundary, and is a minimal vertex-triangulation.

If we crush the triangulation $\T_\Lambda$ along the boundary of
$M$, we get the ideal triangulation $\T^*$ and so the manifold $M$
is homemorphic to $X$. Thus we have that $\T_\Lambda$ is an
inflation of the ideal triangulation $\T^*$.\end{proof}

We have the following summary to the construction of an inflation.

\begin{enumerate}\item Any ideal triangulation $\T^*$ of the interior of a
compact 3--manifold $M$ with boundary, no component of which is a
2--sphere, admits an inflation to a minimal-vertex triangulation
$\T$ of $M$ and $\T$ has normal boundary. An inflation depends on
the choice of frames in the vertex-linking surfaces of $\T^*$;
however, even for fixed frames, there are choices and a finite
number of distinct inflations.

\item The ideal vertices ``inflate" to one-vertex-triangulations of the
boundary components of $M$, which are the triangulations induced on
$\bdy M$ by $\T$.
\item The vertex-linking surfaces in $\T^*$ inflate to
boundary-linking normal surfaces in $\T$.
\item The tetrahedra of $\T^*$ become tetrahedra of $\T$ and some of
the face identifications of $\T^*$ are taken as face identifications
of $\T$.
\item The triangulation $\T$ combinatorially crushes to $\T^*$ along
the boundary of $M$. A combinatorial crushing is
unique.\end{enumerate}

\subsection{Inflations with complexity} We have established that
any ideal triangulation of the interior of a compact $3$--manifold
$X$ has an inflation to a minimal vertex triangulation of $X$. We
now show that we can precisely determine the number of tetrahedra in
the inflation triangulation as a function of the number of
tetrahedra in the given ideal triangulation and a complexity
(number) assigned to the frame of the inflation.

Suppose $\T^*$ is an ideal triangulation of the interior of the
compact $3$--manifold $X$ and $\Lambda = \{\xi_1\ldots,\xi_V\}$ is a
collection of frames in the vertex-linking surfaces of $\T^*$.
Above, in the inflation of $\T^*$ for the collection of frames
$\Lambda$, we did not take particular notice of the number of
tetrahedra added; however, without purposely adding extra tetrahedra
at the various steps of the inflation, the only choices that effect
the number of tetrahedra added are during an inflation at a face of
$\T^*$, and then only if there is more than one edge of the frames
in the face. At this step and in the case of two edges or three
edges in the face, we can arbitrarily choose the order in which we
add the tetrahedra associated to these edges. (Recall if we have two
or three edges in a face, then the inflation in that face is
equivalent to adding a pyramid or a prism, respectively, and our
choices then are in choosing diagonals in the quadrilateral faces of
these pyramids or prisms when we subdivide them into tetrahedra.)
These choices are reflected in the configuration polygons and
necessitate the addition of crossings in the configuration polygons.

\begin{figure}[htbp]
\vspace{.5in}
 \psfrag{e}{\footnotesize $E^+$}\psfrag{f}{\footnotesize
$E^-$}\psfrag{g}{\footnotesize $F^+$}\psfrag{h}{\footnotesize
$F^-$}\psfrag{1}{\footnotesize $1$}\psfrag{2}{\footnotesize
$2$}\psfrag{3}{\footnotesize $3$}\psfrag{4}{\footnotesize
$4$}\psfrag{5}{\footnotesize $5$}\psfrag{6}{\footnotesize
$6$}\psfrag{7}{\footnotesize $7$}\psfrag{8}{\footnotesize
$8$}\psfrag{9}{\footnotesize $9$}\psfrag{a}{\footnotesize
$10$}\psfrag{b}{\footnotesize $11$}\psfrag{c}{\footnotesize
$12$}\psfrag{x}{\footnotesize $\td{3}$}\psfrag{y}{\footnotesize
$\td{6}$}\psfrag{z}{\footnotesize $\td{9}$}\psfrag{p}{\footnotesize
$\td{12}$}\psfrag{w}{\footnotesize $\td{11}$}\psfrag{D}{$D_E^+$}
\psfrag{E}{$E$} \psfrag{F}{$F$}\psfrag{d}{$D_F^+$}
\psfrag{B}{\begin{tabular}{l}
$e(\xi_1) = 5; \hspace{.25 in}\times(\xi_1) = 1;$\\
\\
$ b(\xi_1)=2; \hspace{.25 in}v(\xi_1) = 1$\\
\\
$\mathcal{C}(\xi_1) = 8$\\
\end{tabular}}
\psfrag{+}{$+$}\psfrag{-}{$-$}\psfrag{J}{\begin{tabular}{l}
$e(\xi_2) = 5; \hspace{.25 in}\times(\xi_2) = 2;$\\
\\
$b(\xi_2)=2; \hspace{.25 in}v(\xi_2) =1$\\
\\
$\mathcal{C}(\xi_2) = 9$
\end{tabular}}
\psfrag{M}{\Large \bf Example A.\hspace{.25 in}$\xi_1 = \langle
1\rangle\cup \langle
9,3,\overline{6},\overline{4}\rangle$}\psfrag{N}{\Large \bf Example
B. \hspace{.25 in}$\xi_2 = \langle 1\rangle\cup \langle
4,6,12,11\rangle$}\psfrag{v}{\footnotesize$v^*$}
        \vspace{0 in}
        \begin{center}
\includegraphics[width= 4.5 in]{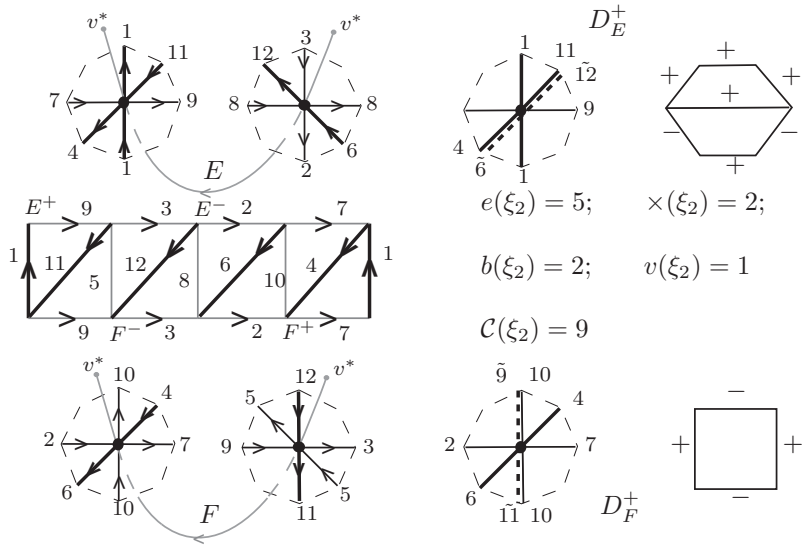}
\caption{Examples of two frames with different complexities in the vertex-linking torus of the two-tetrahedra ideal triangulation of the {\it figure-eight} knot complement in $S^3$.}
\label{f-frames-in-8}
\end{center}
\end{figure}

So, suppose for $\Lambda$ we consider all possible orders for
inflations at the faces of $\T^*$ and then count the number of
crossing configurations in the configuration polygons. The
\underline{minimal} such number will be called the {\it crossing
number for $\Lambda$} and is denoted $\times(\Lambda)$. We let
$b(\Lambda)$ denote the {\it number of branches of $\Lambda$} (the
branches are the components left after removing all branch points of
the various frames in $\Lambda$); and we let $v_b(\Lambda)$ denote
the {\it number of branching points of $\Lambda$}. Finally, we let
$e(\Lambda)$ denote the {\it number of edges in the frames in
$\Lambda$}.

In the above inflation construction, the frames in $\Lambda$
determine various branch configurations in the configuration
polygons. If a branch configuration has branching index $b$, then we
add $b-2$ tetrahedra for this branch configuration. It follows that
if $\xi_j$ is a frame in $\Lambda$,
 $b(\xi_j)$ is the number of branches of $\xi_j$, and
$v_b(\xi_j)$ is the number of branching points, then the number of
tetrahedra needed for the branching configurations coming from
$\xi_j$ is $2[b(\xi_j)-v_b(\xi_j)]$. These are additive functions
over the various frames in $\Lambda$; hence, the number of
tetrahedra added in the inflation over all branching configurations
for the frames in $\Lambda$ is $2[b(\Lambda)-v_b(\Lambda)]$.

We define the {\it complexity of $\Lambda$} as
$$\mathcal{C}(\Lambda) = e(\Lambda) + \times(\Lambda) +
2[b(\Lambda)-v_b(\Lambda)].$$ In particular, the complexity for a
frame $\xi$ in the vertex-linking torus of a cusped 3--manifold with
one cusp is $\mathcal{C}(\xi) = e(\xi)+\times(\xi)+2$.

\begin{thm}\label{blow-up-complexity} Suppose $X$ is a compact
3--manifold with boundary, no component of which is a 2--sphere, and
$\T^*$ is an ideal triangulation of the interior of $X$. If
$\Lambda$ is a collection of frames in the vertex-linking surfaces
of $\T^*$, one from each vertex-linking surface of $\T^*$. An
inflation triangulation $\T_\Lambda$ of $\T^*$ has $card(\T^*) +
\mathcal{C}(\Lambda)$ tetrahedra, where $card(\T^*)$ is the number of
tetrahedra in $\T^*$ and $\mathcal{C}(\Lambda)$ is the complexity of
$\Lambda$.\end{thm}

We give some examples of frames and compute their complexities.

\begin{figure}[htbp]
\vspace{.5 in}
 \psfrag{e}{\footnotesize
$E^+$}\psfrag{f}{\footnotesize $E^-$}\psfrag{g}{\footnotesize
$F^+$}\psfrag{h}{\footnotesize $F^-$}\psfrag{m}{\footnotesize
$G^+$}\psfrag{n}{\footnotesize $G^-$}\psfrag{1}{\footnotesize
$1$}\psfrag{2}{\footnotesize $2$}\psfrag{3}{\footnotesize
$3$}\psfrag{4}{\footnotesize $4$}\psfrag{5}{\footnotesize
$5$}\psfrag{6}{\footnotesize $6$}\psfrag{7}{\footnotesize
$7$}\psfrag{8}{\footnotesize $8$}\psfrag{9}{\footnotesize
$9$}\psfrag{a}{\footnotesize $10$}\psfrag{b}{\footnotesize
$11$}\psfrag{c}{\footnotesize $12$}\psfrag{j}{\footnotesize
$15$}\psfrag{x}{\footnotesize $\td{1}$}\psfrag{y}{\footnotesize
$\td{7}$}\psfrag{z}{\footnotesize $\td{6}$}\psfrag{t}{\footnotesize
$\td{12}$}\psfrag{p}{\footnotesize $\td{3}$}\psfrag{D}{$D_E^+$}
\psfrag{E}{$E$} \psfrag{F}{$F$}\psfrag{d}{$D_F^+$}
\psfrag{G}{$G$}\psfrag{i}{$D_G^+$}\psfrag{B}{\begin{tabular}{l}
$e(\Lambda) = 8; \hspace{.25 in}\times(\Lambda) = 2;$\\
\\
$ b(\Lambda)=4; \hspace{.25 in}v(\Lambda) = 2$\\
\\
$\mathcal{C}(\Lambda) = 14$\\
\end{tabular}}
\psfrag{A}{\begin{tabular}{l}
\\
\\
{\Large\bf Example C.}\\
\hspace{.5 in}\Large{$\xi_{V^*} = \langle 2\rangle\cup\langle 6,12\rangle$}\\
\\
\hspace{.5 in}\Large{$\xi_{W^*} = \langle 15\rangle\cup\langle 7,9,\overline{3},\overline{1}\rangle$}\\
\end{tabular}}
\psfrag{+}{\footnotesize$+$}\psfrag{-}{\footnotesize$-$} \psfrag{v}{\footnotesize$v^*$}
\psfrag{w}{\footnotesize$w^*$}\psfrag{r}{$\xi_{V^*}$}\psfrag{s}{$\xi_{W^*}$}\psfrag{W}{\Large{$W^*$}}\psfrag{V}{\Large{$V^*$}}
        \vspace{0 in}
        \begin{center}
\includegraphics[width= 5 in]{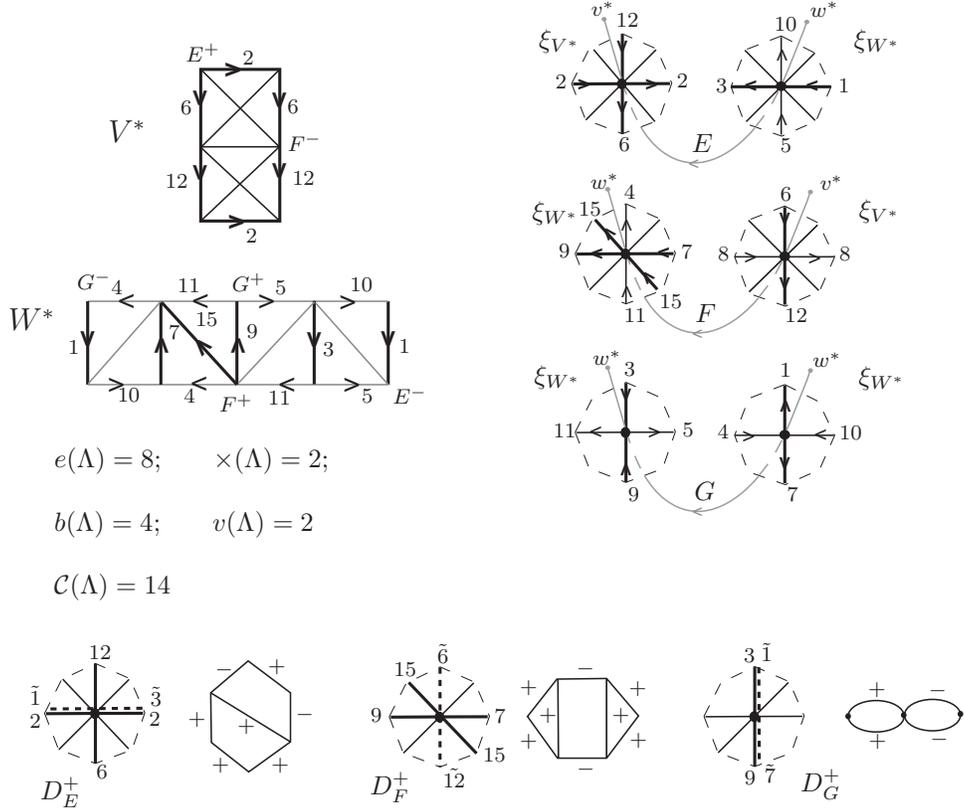}
\caption{Frames in the two vertex-linking tori for a four tetrahedra ideal triangulation of the complement  of the Whitehead Link in $S^3$. The complexity of $\Lambda=\{\xi_{V^*},\xi_{W^*}\}$ is greater than the sum of the complexities of $\xi_{V^*}$ and $\xi_{W^*}$. }
\label{f-frames-in-whitehead}
\end{center}
\end{figure}

In the first two examples (Figure \ref{f-frames-in-8}), we have the
two-tetrahedron ideal triangulation of the {\it figure-eight} knot
complement in the $3$--sphere. This ideal triangulation has one
ideal vertex and the vertex-linking surface is a torus. There are
two edges. In Example (A) the frame is given as $\xi_1 = \langle
1\rangle\cup\langle 9,3,\overline{6},\overline{4}\rangle$, where
$\overline{e}$ indicates the reverse direction to the arrow on the
edge $e$, and has two branches; one is the standard meridian,
labeled $\langle 1\rangle$, and the other is the standard
homological longitude, labeled $\langle
9,3,\overline{6},\overline{4}\rangle$. This frame $\xi_1$ has
complexity $8$. Hence, the inflation of this two-tetrahedron ideal
triangulation of the {\it figure-eight} knot complement, using the
frame $\xi_1$, gives a one-vertex triangulation of the {\it
figure-eight} knot exterior having $10$ tetrahedra. In this case, we
conjecture the inflation triangulation is also a minimal
triangulation. In the next section, we construct this inflation. In
Example (B) the frame is given as $\xi_2 = \langle 1\rangle\cup
\langle 4,6,12,11\rangle$ and again has two branches; one branch is
the meridian but the other, while a longitude, is not the
homological longitude. The complexity in this example is $9$. Hence,
we have two frames in the same vertex-linking surface having
different complexities.

In the next example (Figure \ref{f-ideal-whitehead}) in the
Appendix, we have a four-tetrahedron ideal triangulation of the
complement of the Whitehead link in the $3$--sphere. This ideal
triangulation has two ideal vertices and the vertex-linking surface
at each is a torus. There are three edges, $E, F$ and $G$. We denote
the ideal vertices $v^*$ and $w^*$ and their vertex-linking tori $V^*$
and $W^*$, respectively. In $V^*$ the frame is $\xi_{V^*} = \langle
2\rangle\cup \langle 6,12\rangle$ and in $W^*$ the frame is $\xi_{W^*} =
\langle 15\rangle\cup \langle 7,9,\overline{3},\overline{1}\rangle$.
We set $\Lambda = \{\xi_{V^*},\xi_{W^*}\}$. The collection $\Lambda$ has
four branches and two vertices; its complexity is $14$. Note that
the edges $E$ and $F$ meet different frames at their ends.
Individually, the frame $\xi_{V^*}$ has complexity $5$ and the frame
$\xi_{W^*}$ has complexity $7$; however, we can not add these
complexities to get the complexity of $\Lambda$. This example
displays that we must consider the frames at each end of an edge in
an inflation. If we should first inflate along the frame $\xi_{V^*}$,
then we change the induced triangulation on the vertex-linking torus
at $w^*$ and change the frame and its complexity there.


\section{Examples of the Inflation Construction.}\label{blow-up-examples} In this section we
give two examples of the inflation construction. The first is an
inflation of the two-tetrahedra ideal triangulation of the
complement of the {\it figure-eight} knot in $S^3$ and the second is
an inflation of the one-tetrahedron triangulation of the Gieseking
manifold. As we have remarked elsewhere, we considered using other
examples as these manifolds have been extensively studied and are
repeatedly used as examples; however, we finally decided that
familiarity with these examples may be useful in introducing the
ideas of the inflation construction.

\vspace{.25 in}\noindent {\bf Example. Inflation of figure-eight
knot complement.}

\vspace{.15 in}\noindent {\bf Step 1.} Given an ideal triangulation.

\begin{figure}[htbp]
\vspace{.25 in}
            \psfrag{p}{$(p)$}\psfrag{q}{$(p')$}

            \psfrag{1}{\scriptsize{$1$}}\psfrag{2}{\scriptsize{$2$}}
            \psfrag{3}{\scriptsize{$3$}}\psfrag{4}{\scriptsize{$4$}}
            \psfrag{5}{\scriptsize{$5$}}\psfrag{6}{\scriptsize{$6$}}
            \psfrag{7}{\scriptsize{$7$}}\psfrag{8}{\scriptsize{$8$}}
            \psfrag{9}{\scriptsize{$9$}}\psfrag{a}{\scriptsize{$10$}}
            \psfrag{b}{\scriptsize{$11$}}\psfrag{c}{\scriptsize{$12$}}
\psfrag{w}{$0$}\psfrag{x}{$1$}\psfrag{y}{$2$}\psfrag{z}{$3$}
            \psfrag{F}{\footnotesize $F$}\psfrag{E}{\footnotesize $E$}\psfrag{e}{\small $2$}
            \psfrag{f}{\small $3$}\psfrag{g}{\small $4$}\psfrag{h}{\small $6$}
            \psfrag{j}{\small $7$}\psfrag{k}{\small $9$}\psfrag{l}{\small $1$}
            \psfrag{i}{\begin{tabular}{c}
$(p)(012)\leftrightarrow (p')(012);$\\
\\
$(p)(013)\leftrightarrow (p')(312);$\\
\\
$(p)(023)\leftrightarrow (p')(310);$\\
\\
$(p)(123)\leftrightarrow (p')(320)$\\
\end{tabular}}
        \vspace{0 in}
        \begin{center}
\includegraphics[width= 5 in]{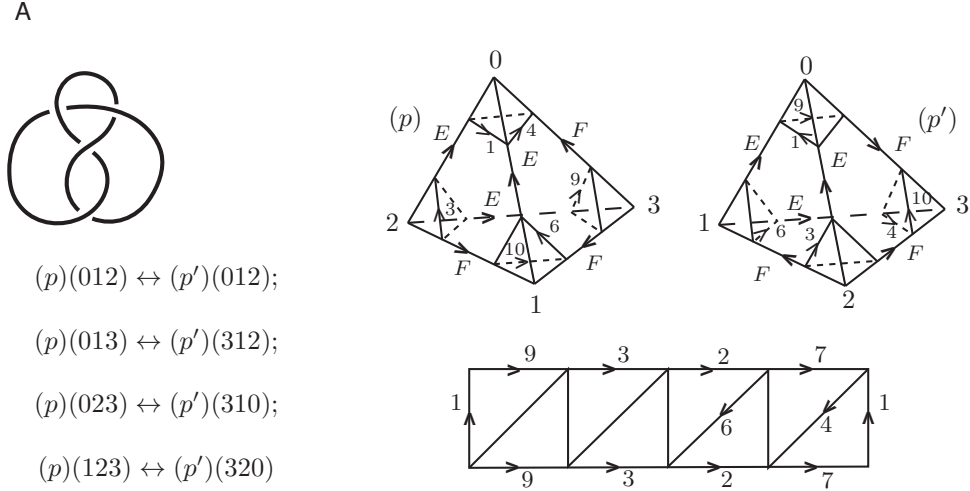} 
\caption{The two-tetrahedra ideal triangulation of the {\it figure-eight} knot complement in $S^3$, along with the induced triangulation on its vertex-linking torus.} 
\label{f-ideal-8}
\end{center}
\end{figure}

For this example, the given ideal triangulation $\T^*$ is the
two-tetrahedra ideal triangulation of the {\it figure-eight} knot
complement in $S^3$ given in  Figure \ref{f-ideal-8}.

We have two tetrahedra, $(p)$ and $(p')$ with face identifications:

\vspace{.15 in}
\begin{tabbing}
\hspace{1 in}\=$(p)(012)\leftrightarrow (p')(012)$\hspace{.5 in}\=$(p)(013)\leftrightarrow (p')(312)$\\

\>$(p)(023)\leftrightarrow (p')(310)$\>$(p)(123)\leftrightarrow (p')(320)$\\
\end{tabbing}

\vspace{.25 in}\noindent {\bf Step 2.} Construct the vertex-linking
surface and choose a frame.

The vertex-linking surface also is shown in Figure \ref{f-ideal-8}.
We shall use as a frame $\xi
=<1>\cup<9,3,\overline{6},\overline{4}>$, which is given as an
example in Figure \ref{f-example-frames}. It is the standard
meridian/longitude pair for the {\it figure-eight} knot in $S^3$.

\vspace{.25 in}\noindent {\bf Step 3.} Direct each branch,
successively label edges in the branches, and determine the
transverse direction for each branch.

\begin{figure}[htbp]

            \psfrag{F}{\footnotesize $F$}\psfrag{E}{\footnotesize $E$}
            \psfrag{e}{\small $2$}
            \psfrag{f}{\small $3$}\psfrag{g}{\small $4$}\psfrag{h}{\small $6$}
            \psfrag{j}{\small $7$}\psfrag{k}{\small $9$}\psfrag{l}{\small $1$}

        \vspace{0 in}
        \begin{center}
\includegraphics[width= 3 in]{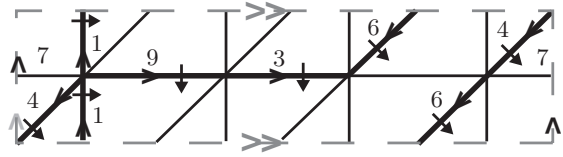} 
\caption{A one vertex, index $4$ frame with two branches: $\langle 1\rangle$ and $\langle 9,3,\overline{6},\overline{4}\rangle$, along with transverse directions. This frame has complexity $8$.}
\label{f-transverse-frame-8}
\end{center}
\end{figure}

There are two branches for this frame. One branch is the single edge
$1$, the meridian. In this example, we chose directions on the edges
in the induced triangulation of the vertex-linking torus to aid the
reader in the face identifications; we shall utilize these labels
and directions. So, the branch $\langle 1\rangle$ is given the same
direction as the edge $1$ and we have indicated a transverse
direction for this branch in Figure \ref{f-transverse-frame-8}.

The second branch and direction is $\langle
9,3,\overline{6},\overline{4}\rangle$, where $\overline{e}$ means
the edge $e$ taken in the opposite direction to that used in the
face identifications and given in Figure \ref{f-ideal-8}. This
branch corresponds to the homological longitude slope in $S^3$. The
transverse direction is given in Figure \ref{f-transverse-frame-8}.

\noindent {\bf Step 4.} Determine the configuration polygons; using
the transverse directions, determine the directions on the boundary
edges of the configuration polygons.

\begin{figure}[htbp]

            \psfrag{F}{$D_F^+$}\psfrag{D}{$D_E^+$}\psfrag{a}{\small $1^0$}
            \psfrag{b}{\small $1^1$}\psfrag{c}{\small $9^0$}
            \psfrag{d}{\small $\td{9}^1$}\psfrag{e}{\small $\td{3}^0$}
            \psfrag{f}{\small $\td{3}^1$}\psfrag{h}{\small $\overline{4}^0$}
            \psfrag{k}{\small $\overline{6}^1$}\psfrag{j}{\small $\td{\overline{6}}^0$}
\psfrag{i}{\small $\overline{4}^1$}\psfrag{u}{\scriptsize
$\overline{6}$}\psfrag{w}{\scriptsize
$\td{9}$}\psfrag{x}{\scriptsize $\td{3}$}\psfrag{y}{\scriptsize
$\td{\overline{6}}$}\psfrag{z}{\scriptsize
$\overline{4}$}\psfrag{m}{\scriptsize $7$}\psfrag{n}{\scriptsize
$11$}\psfrag{2}{\scriptsize $2$}\psfrag{7}{\scriptsize
$7$}\psfrag{g}{\scriptsize $10$}\psfrag{1}{\scriptsize $1$}

        \vspace{0 in}
        \begin{center}
\includegraphics[width= 5 in]{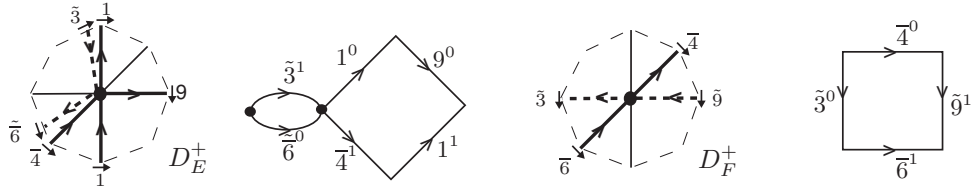}
\caption{Labeled configuration polygons for the inflation of the {\it figure eight} knot complement.}
\label{f-config-polygon-labels-8}
\end{center}
\end{figure}

The configuration polygons are given in Example A in Figure
\ref{f-frames-in-8}; we give them here in Figure
\ref{f-config-polygon-labels-8} with labels and transverse
directions. Note we use $\overline{6}$ and $\overline{4}$ and give
them the direction induced by that of the branch $\langle
9,3,\overline{6},\overline{4}\rangle$.

\noindent {\bf Step 5.} Add a tetrahedron for each edge in the
frame.

\vspace{.125 in}For this example we have 5 edges, giving 5
tetrahedra: $(1)(0123), (9)(0123),$\\$(3)(0123),(\overline{6})(0123)$
and $(\overline{4})(0123)$.

\vspace{.25 in}\noindent {\bf Step 6.} Inflation at the faces of
$\T$.

\begin{figure}[htbp]

            \psfrag{p}{$(p)$}
            \psfrag{q}{$(p')$}\psfrag{1}{\small $1$}\psfrag{0}{\small $0$}
\psfrag{2}{\small $2$}\psfrag{3}{\small $3$}
\psfrag{x}{\scriptsize{$3$}}\psfrag{y}{\scriptsize{$1$}}\psfrag{X}{$(3)$}
\psfrag{Y}{$(1)$} \psfrag{i}{\begin{tabular}{c}
           $(p)(012)\leftrightarrow(3)(320)$;\\
           \\
$(3)(321)\leftrightarrow(1)(032)$;\\
\\
            $(1)(132)\leftrightarrow(p')(012)$\\
            \end{tabular}}

        \begin{center}
\includegraphics[width= 3.75 in]{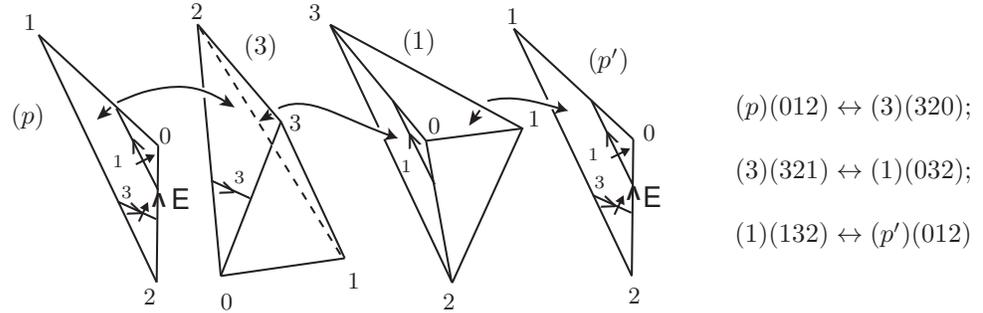}
\caption{Inflation at the face $(p)(012)\leftrightarrow (p')(012)$, which meets the frame in two edges $1$ and $3$.}
\label{f-blow-up-face-8-1}
\end{center}
\end{figure}

\vspace{.15 in}\noindent {\bf -} Inflation at the face
$(p)(012)\leftrightarrow (p')(012)$; this face contains the two
edges  $1$ and $3$. See Figure \ref{f-blow-up-face-8-1}. Notice that
we have made a choice of order in adding the tetrahedra $(1)$ and
$(3)$ (same as selecting a diagonal if we think of the inflation at
the face as adding a prism). Our choice was made by looking at the
configuration polygon in Figure \ref{f-config-polygon-labels-8} and
seeing that we avoid a crossing if we make the selection by first
attaching the tetrahedron $(p)$ to the tetrahedron $(3)$; then
attaching $(3)$ to $(1)$; and, finally, ending by attaching $(1)$ to
$(p')$. This is discussed above where we discussed inflations with
complexity.

\begin{figure}[htbp]

            \psfrag{p}{$(p)$}\psfrag{E}{$E$}
            \psfrag{q}{$(p')$}\psfrag{1}{\small $1$}\psfrag{0}{\small $0$}
\psfrag{2}{\small $2$}\psfrag{3}{\small $3$}
\psfrag{x}{\scriptsize{$\overline{6}$}}\psfrag{y}{\scriptsize{$\overline{4}$}}
\psfrag{X}{$(\overline{6})$} \psfrag{Y}{$(\overline{4})$}
\psfrag{i}{\begin{tabular}{c}
           $(p)(013)\leftrightarrow(\overline{4})(132)$;\\
           \\
$(\overline{4})(032)\leftrightarrow(\overline{6})(213)$;\\
\\
            $(\overline{6})(203)\leftrightarrow(p')(312)$\\
            \end{tabular}}

        \begin{center}
\includegraphics[width= 3.5 in]{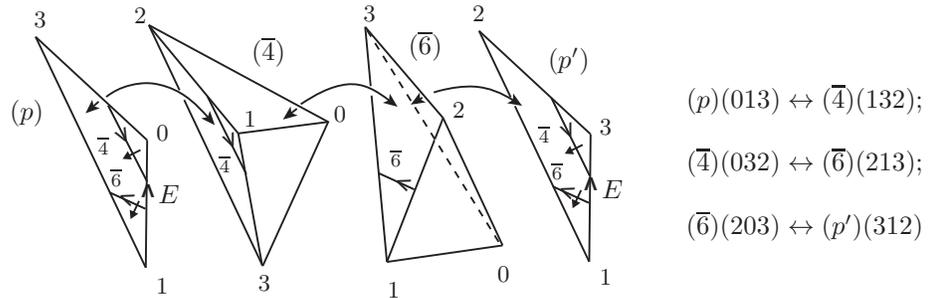}
\caption{Inflation at the face $(p)(013)\leftrightarrow (p')(312)$, which meets the frame in two edges $\overline{6}$ and $\overline{4}$.} 
\label{f-blow-up-face-8-2}
\end{center}
\end{figure}

\vspace{.15 in}\noindent {\bf -} Inflation at the face
$(p)(013)\leftrightarrow (p')(312)$; this face contains the two
edges $\overline{6}$ and $\overline{4}$. See Figure
\ref{f-blow-up-face-8-2}. Again we have used the configuration
polygon to make a selection for the order we add the tetrahedra
$(\overline{6})$ and $(\overline{4})$ to avoid adding additional
crossings. Here we avoid a crossing if we make the selection by
first attaching the tetrahedron $(p)$ to the tetrahedron
$(\overline{4})$; then attaching $(\overline{4})$ to
$(\overline{6})$; and, finally, ending by attaching $(\overline{6})$
to $(p')$.

\begin{figure}[htbp]

\psfrag{p}{$(p)$}
            \psfrag{q}{$(p')$}\psfrag{1}{\small $1$}\psfrag{0}{\small $0$}
\psfrag{2}{\small $2$}\psfrag{3}{\small $3$}
\psfrag{x}{\scriptsize{$9$}} \psfrag{9}{$(9)$}
 \psfrag{i}{\begin{tabular}{c}
           $(p)(023)\leftrightarrow(9)(320)$;\\
           \\
$(9)(321)\leftrightarrow(p')(310)$\\
            \end{tabular}}
        \begin{center}
\includegraphics[width= 3 in]{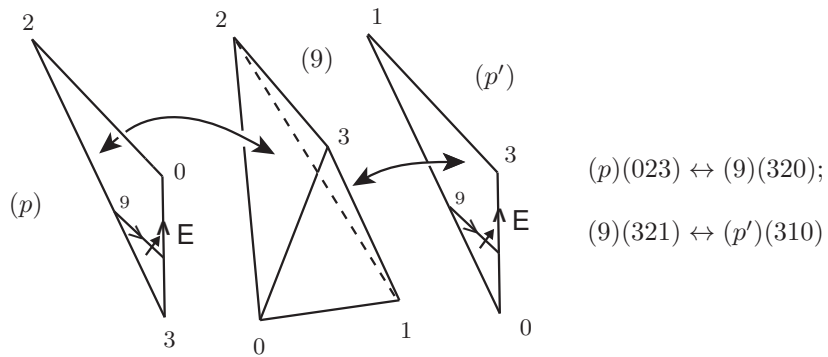} 
\caption{The face $(p)(023)\leftrightarrow(p')(310)$ meets the frame in one edge $9$.} 
\label{f-blow-up-face-8-3}
\end{center}
\end{figure}

\vspace{.15 in}\noindent {\bf -} Inflation at the face
$(p)(023)\leftrightarrow (p')(310)$; this face contains the single
edge $9$. See Figure \ref{f-blow-up-face-8-3}.

\vspace{.25 in}\noindent {\bf Step 7.} Inflation at the edges of
$\T$. There are two edges, $E$ and $F$. The inflation at an edge is
determined by the configuration polygon at that edge. In Step 4,
Figure \ref{f-config-polygon-labels-8}, we give the configuration
polygons for the edges $E$ and $F$.

\begin{figure}[htbp]
\vspace{.25 in}
           \psfrag{e}{\small $1^1$}
            \psfrag{f}{\small $\td{3}^1$}\psfrag{h}{\small $9^0$}
            \psfrag{j}{\small $\td{\overline{6}}^0$}
\psfrag{i}{\small $\overline{4}^1$}\psfrag{g}{\small $1^0$}

            \psfrag{b}{$b$}\psfrag{c}{$(3)$}\psfrag{f}{$(b_1^*)$}\psfrag{d}{$(b_2^*)$}
            \psfrag{9}{$(9)$}\psfrag{6}{$(\overline{6})$}
            \psfrag{a}{$(1)$}\psfrag{4}{($\overline{4})$}
            \psfrag{B}{\large $(b^*)$}\psfrag{0}{\footnotesize $0$}
\psfrag{1}{\footnotesize $1$}\psfrag{2}{\footnotesize
$2$}\psfrag{3}{\footnotesize $3$}\psfrag{m}{\begin{tabular}{c}
           $(\overline{6})(012)\leftrightarrow(3)(013)$\\
            \end{tabular}}\psfrag{n}{\begin{tabular}{c}
           $(1)(012)\leftrightarrow(b^*)(120)$;\hspace{.25in}$(9)(012)\leftrightarrow(b^*)(230)$\\
           \\
$(1)(013)\leftrightarrow(b^*)(b30)$;\hspace{.25 in}$(\overline{4})(013)\leftrightarrow(b^*)(1b0)$\\
            \end{tabular}}

        \vspace{0 in}
        \begin{center}
\includegraphics[width= 5 in]{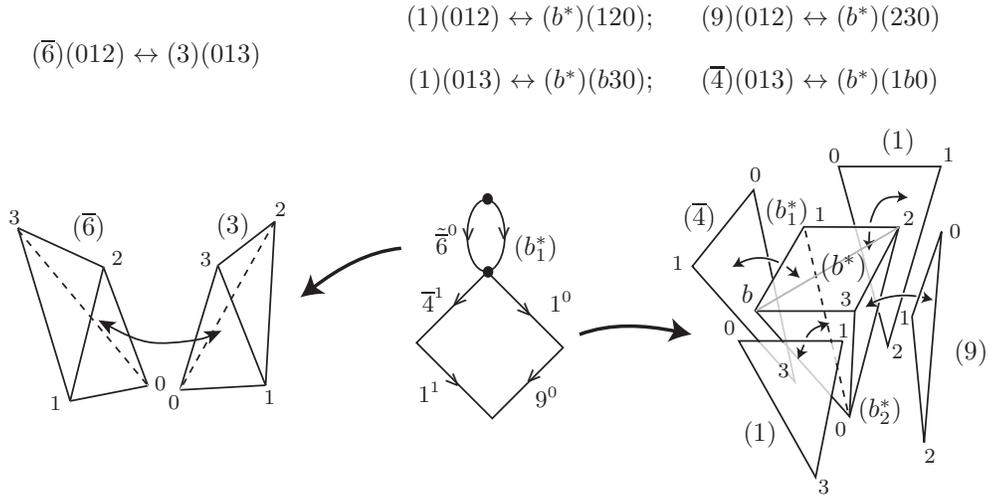}
\caption{Inflation at the edge $E$ where the configuration polygon splits into a generic bi-gon and a branch pyramid, which is subdivided into $(b_1^*),b=3$ and $(b_2^*),b=1$.}
\label{f-blow-up-edge-8-E}
\end{center}
\end{figure}

\vspace{.15 in}\noindent {\bf -} Inflation at the edge $E$. In this
case the configuration polygon splits into two independent polygons;
one is a generic polygon and the other is a branch polygon for a
branch of index $4$. See Figure \ref{f-blow-up-edge-8-E}.

\begin{figure}[htbp]
\vspace{.25 in}
            \psfrag{e}{\small $\td{3}^0$}
            \psfrag{d}{\small $\td{9}^1$}
\psfrag{h}{\small $\overline{4}^0$}\psfrag{k}{\small
$\overline{6}^1$}

            \psfrag{c}{$(3)$}\psfrag{C}{$(c)$}
            \psfrag{9}{$(9)$}\psfrag{6}{$(\overline{6})$}
            \psfrag{4}{($\overline{4})$}
            \psfrag{0}{\footnotesize $0$}
\psfrag{1}{\footnotesize $1$}\psfrag{2}{\footnotesize
$2$}\psfrag{3}{\footnotesize $3$}\psfrag{m}{\begin{tabular}{c}
           $(\overline{6})(013)\leftrightarrow(c)(023)$\\
           $(\overline{4})(012)\leftrightarrow(c)(021)$\\
           $(3)(012)\leftrightarrow(c)(130)$\\
          $(9)(013)\leftrightarrow(c)(132)$\\
\end{tabular}}

        \vspace{0 in}
        \begin{center}
\includegraphics[width=4in]{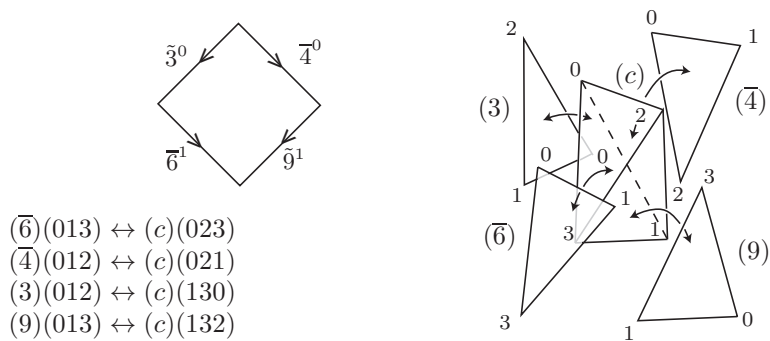}
\caption{Inflation at the edge $F$ where the configuration polygon is a crossing.} 
\label{f-blow-up-edge-8-F}
\end{center}
\end{figure}

\vspace{.15 in}\noindent {\bf -} Inflation at the edge $F$. In this
case the configuration polygon is a crossing. See Figure
\ref{f-blow-up-edge-8-F}.

\vspace{.25 in}\noindent {\bf Step 8.} Finish.

We have face identifications from $\T^*$ which may remain. In this
example, the face identification $(p)(123)\leftrightarrow (p')(320)$
is retained.

We have cones over branch configurations, which may need to be
subdivided into tetrahedra. This is done by first triangulating the
branch polygon, without adding vertices, and then using the cone
structure to extend the triangulation of the polygon to a
triangulation of the cone. For $b$-gons we have $(b-2)$ tetrahedra.
While the number of tetrahedra for a $b$--gon will remain fixed at
$(b-2)$, there are numerous options for subdividing the polygons
(the number of options is the $(b-2)$ Catalan number). In this
example, the branch polygon is a quadrilateral and we have two
choices, depending on the choice of diagonal, for triangulating the
polygon and thus two choices for triangulating the pyramid $(b^*)$.
For example, if we choose the diagonal from the vertex $b$ to the
vertex $2$; then we have a triangulation of the pyramid with two
tetrahedra, $(b_1^*)(012b)$ and $(b_2^*)(023b)$, and a single face
identification $(b_1^*)(02b)\leftrightarrow(b_2^*)(02b)$.

We collect the tetrahedra and face identifications into an array
following  the notational conventions of {\footnotesize\textsc
REGINA} \cite{burton-regina}. A tetrahedron (using the abbreviation
``tet") from our collection is given in the first column, its faces
are given in the top row, and the face identifications are given at
the intersection of a row with a column. For example, to know the
face identification of the face $(013)$ of the tetrahedron $(9)$, go
down the first column to $(9)$ and then cross to the column under
$(013)$ where we fined $(c)(132)$; hence, the face identification is
$(9)(013)\leftrightarrow(c)(132)$. To include the subdivision of the
pyramid $(b^*)$, we arbitrarily choose to subdivide $(b^*)$ into two
tetrahedra $(b_1^*)(012b)$ and $(b_2^*)(0b23)$; hence, in the array,
we indicate the vertex $b$ by $3$ in $(b_1^*)$ and indicate the
vertex $b$ by $1$ in $(b_2^*)$.

Notice that the two faces $(b_1^*)(123)$ and $(b_2^*)(123)$ are not
identified and become the boundary of the compact manifold, which is
the exterior of the {\it figure-eight} knot in $S^3$. If one follows
the identifications of the edges in the boundary of the
quadrilateral, which is the base of the pyramid $(b^*)$, using our
notational conventions,  we have
$$(b_1^*)(12)\rightarrow (1)(01)\rightarrow (b_2^*)(13)$$ and
$$(b_2^*)(23)\rightarrow (9)(01)\rightarrow (c)(13)\rightarrow
(3)(01)\rightarrow (\overline{6})(01)\rightarrow$$ $$\rightarrow(c)(02) \rightarrow
(\overline{4})(01)\rightarrow (b_1^*)(13).$$

The boundary is the one-vertex, two triangle triangulation of the
torus. We see, clearly, here that the free edges added when
inflating at a face along a single branch are all identified and
form an edge in the boundary. This guided our choice of notation for
this edge being $(01)$ along all branches. We also see the crossing
where we have free edges, $c(13)$ and $c(02)$, as opposite edges of
the tetrahedron $(c)$; in this example, they are on the same branch.

\begin{tabbing}\hspace{.5 in}\=\bf {tet}\hspace{.55 in}\=$(012)$\hspace{.55 in}\=$(013)$
\hspace{.55 in}\=$(023)$\hspace{.55 in}\=$(123)$\\
\>$(p)$\>$(3)(320)$\>$(\overline{4})(132)$\>$(9)(320)$
\>$(p')(320)$\\
\>$(p')$\>$(1)(132)$\>$(9)(123)$\>$(p)(321)$
\>$(\overline{6})(032)$\\
\>$(1)$\>$(b_1^*)(120)$\>$(b_2^*)(130)$\>$(3)(312)$\>$(p')(021)$\\
\>$(3)$\>$(c)(130)$\>$(\overline{6})(012)$\>$(p)(210)$
\>$(1)(230)$\\
\>$(\overline{4})$\>$(c)(021)$\>$(b_1^*)(130)$\>$(\overline{6})(231)$
\>$(p)(031)$\\
\>$(\overline{6})$\>$(3)(013)$\>$(c)(023)$\>$(p')(132)$
\>$(\overline{4})(302)$\\
\>$(9)$\>$(b_2^*)(230)$\>$(c)(132)$\>$(p)(320)$
\>$(p')(013)$\\
\>$(c)$\>$(\overline{4})(021)$\>$(3)(201)$\>$(\overline{6})(013)$
\>$(9)(031)$\\
\>$(b_1^*)$\>$(1)(201)$\>$(\overline{4})(301)$\>$(b_2^*)(021)$
\>{\bf bddry}\\
\>$(b_2^*)$\>$(b_1^*)(032)$\>$(1)(301)$\>$(9)(201)$
\>{\bf bddry}\\
\end{tabbing}

We conclude this example by remarking that we suspect this is a
minimal triangulation of the {\it figure-eight} knot exterior in
$S^3$. It has $10$ tetrahedra. There are many other non-isomorphic
$10$-tetrahedra triangulations of the {\it figure-eight} knot
exterior; some of the examples we know are not inflations of the
two-tetrahedron ideal triangulation of the {\it figure-eight} knot
complement but of a three-tetrahedra ideal triangulation of the {\it
figure-eight} knot complement, which is formed by a
$2\leftrightarrow 3$ Pachner move on the two-tetrahedra ideal
triangulation.

\vspace{.25 in}\noindent {\bf Example. Inflation of the Gieseking
manifold.}

\vspace{.15 in}\noindent {\bf Step 1.} Given an ideal triangulation.

\begin{figure}[htbp]
\vspace{.25 in}
            \psfrag{p}{$(p)$}

            \psfrag{1}{\scriptsize{$1$}}
            \psfrag{3}{\scriptsize{$3$}}\psfrag{4}{\scriptsize{$4$}}\psfrag{2}{\scriptsize{$2$}}
            \psfrag{5}{\scriptsize{$5$}}
            \psfrag{a}{$0$}
            \psfrag{b}{$1$}\psfrag{c}{$2$}\psfrag{d}{$3$}
\psfrag{x}{\small $5$}\psfrag{y}{\small $2$}\psfrag{z}{\small
$3$}\psfrag{u}{\small $4$}
            \psfrag{E}{\footnotesize $E$}
            \psfrag{e}{\begin{tabular}{c}
$(p)(012)\leftrightarrow (p)(302)$\\
\\
$(p)(013)\leftrightarrow (p)(123)$\\
\end{tabular}}\psfrag{K}{Klein bottle}
        \vspace{0 in}
        \begin{center}
\includegraphics[width=3.5in]{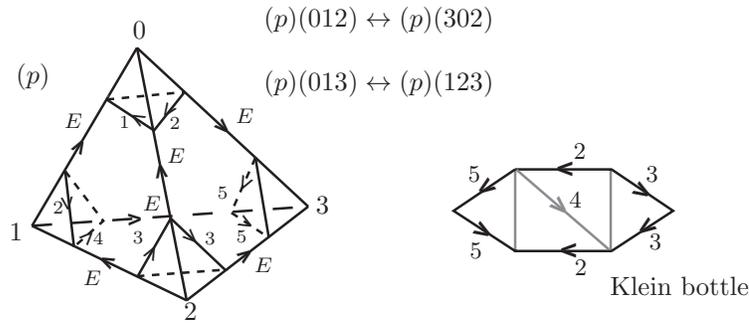} 
\caption{The one-tetrahedron ideal triangulation of the Gieseking manifold (non-orientable), along with the induced triangulation on its vertex-linking Klein bottle.} 
\label{f-geisking}
\end{center}
\end{figure}

For this example, the given ideal triangulation $\T^*$ is the
one-tetrahedra ideal triangulation of the Gieseking manifold given
in  Figure \ref{f-geisking}.

We have one tetrahedra, $(p)$, with face identifications:

\begin{tabbing}
\hspace{1 in}\=$(p)(012)\leftrightarrow (p)(302)$\hspace{.5 in}\=$(p)(013)\leftrightarrow (p)(123)$\\
\end{tabbing}

\noindent {\bf Step 2.} Construct the vertex-linking surface and
choose a frame.

The vertex-linking surface also is shown in Figure \ref{f-geisking};
notice that it is a Klein bottle and the Gieseking manifold is
non-orientable. We shall use as a frame $\xi=<5>\cup<2>\cup<3>$,
which was given in the example in Figure \ref{f-example-frames}.

\noindent {\bf Step 3.} Direct each branch, successively label edges
in the branches, and determine the transverse direction for each
branch.

\begin{figure}[htbp]

            \psfrag{2}{\small $2$}\psfrag{1}{\small $1$}\psfrag{6}{\small $6$}
            \psfrag{3}{\small $3$}\psfrag{4}{\small $4$}\psfrag{5}{\small $5$}

        \vspace{0 in}
        \begin{center}
\includegraphics[width=2.5 in]{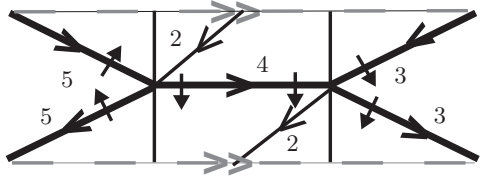}
\caption{A double index three frame with three branches $\langle 5\rangle$, $\langle 2\rangle$, and $\langle 3\rangle$, along with their transverse directions. This frame has complexity $6$.}
\label{f-transverse-frame-geisking}
\end{center}
\end{figure}

There are three branches for this frame. Each branch has just one
edge. As above, we have directions on the edges in the induced
triangulation of the vertex-linking Klein bottle, which were given
to aid the reader in the face identifications; we utilize these
labels and directions. All branches are given the direction of their
edges. However, here when finding the transverse directions, we see
that the orientation reversing edges $3$ and $5$ change the
transverse directions from what we would have in the orientable
case. The transverse directions are given in Figure
\ref{f-transverse-frame-geisking}.

\noindent {\bf Step 4.} Determine the configuration polygons; using
the transverse directions, determine the directions on the boundary
edges of the configuration polygons.

\begin{figure}[htbp]
           \psfrag{D}{$D_E^-$}\psfrag{a}{\small $\td{5}^0$}
            \psfrag{b}{\small $\td{5}^1$}\psfrag{d}{\small $\td{2}^1$}
            \psfrag{c}{\small $2^0$}
            \psfrag{e}{\small $3^0$}\psfrag{f}{\small $3^1$}
            \psfrag{x}{\small $3$}\psfrag{y}{\small $2$}\psfrag{v}{\small $\td{2}$}
\psfrag{4}{\small $4$}\psfrag{5}{\small $\td{5}$}\psfrag{6}{\small
$6$}\psfrag{-}{$-$}\psfrag{+}{$+$}

        \vspace{0 in}
        \begin{center}
\includegraphics[width=4 in]{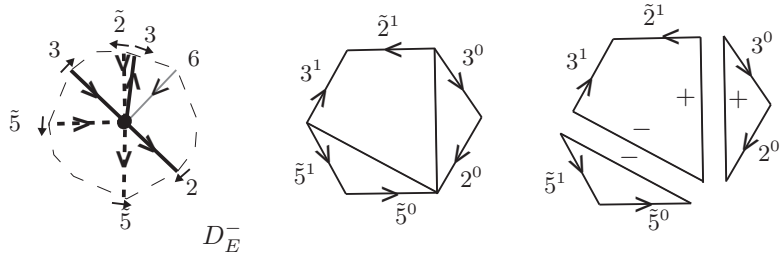}
\caption{Labeled configuration polygons for the inflation of the Gieseking manifold.} 
\label{f-config-polygon-labels-geisking}
\end{center}
\end{figure}

There is only one edge and therefore only one configuration polygon;
it is given  in Figure \ref{f-config-polygon-labels-geisking} with
labels and transverse directions.

\noindent  {\bf Step 5.} Add a tetrahedron for each edge in the
frame.

For this example we have three edges and therefore add three
tetrahedra: $(3)(0123),\\ (2)(0123)$ and $(5)(0123)$.

\noindent {\bf Step 6.} Inflation at the faces of $\T$.

\begin{figure}[htbp]

\psfrag{p}{$(p)$}\psfrag{E}{\small $E$}
            \psfrag{1}{\small $1$}\psfrag{0}{\small $0$}
\psfrag{2}{\small $2$}\psfrag{3}{\small $3$}
\psfrag{x}{\scriptsize{$5$}} \psfrag{X}{$(5)$}
 \psfrag{i}{\begin{tabular}{c}
           $(p)(013)\leftrightarrow(5)(231)$;\\
           \\
$(5)(230)\leftrightarrow(p)(123)$\\
            \end{tabular}}
        \begin{center}
\includegraphics[width=2.75in]{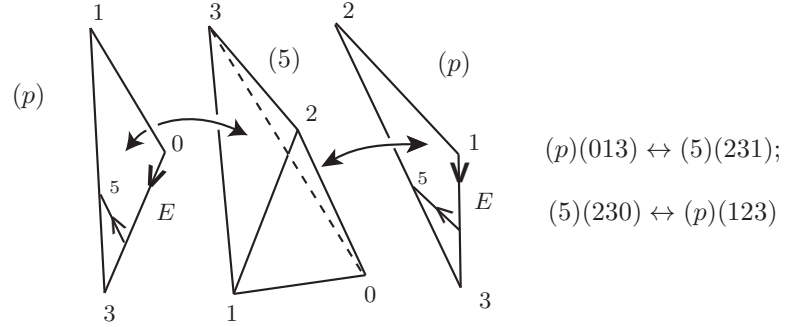}
\caption{The face $(p)(013)\leftrightarrow(p)(123)$ meets the frame in one edge $5$.} 
\label{f-blow-up-face-geisking-1}
\end{center}
\end{figure}

\vspace{.15 in}\noindent {\bf -} Inflation at the face
$(p)(013)\leftrightarrow (p)(123)$; this face contains the single
edge $5$. See Figure \ref{f-blow-up-face-geisking-1}.

\begin{figure}[htbp]

            \psfrag{p}{$(p)$}
            \psfrag{1}{\small $1$}\psfrag{0}{\small $0$}
\psfrag{2}{\small $2$}\psfrag{3}{\small $3$}
\psfrag{x}{\scriptsize{$3$}}\psfrag{y}{\scriptsize{$2$}}\psfrag{X}{$(3)$}
\psfrag{Y}{$(2)$} \psfrag{i}{\begin{tabular}{c}
           $(p)(012)\leftrightarrow(3)(321)$;\\
           \\
$(3)(320)\leftrightarrow(2)(213)$;\\
\\
            $(2)(203)\leftrightarrow(p)(302)$\\
            \end{tabular}}

        \begin{center}
\includegraphics[width=4 in]{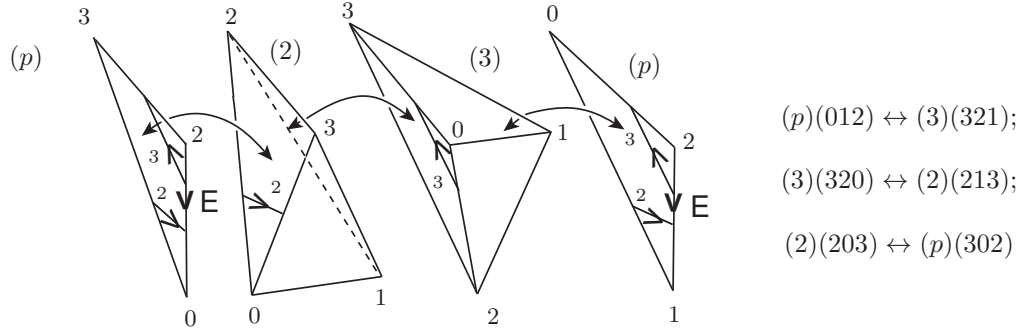}
\caption{Inflation at the face $(p)(012)\leftrightarrow (p)(302)$, which meets the frame in two edges $2$ and $3$.}
\label{f-blow-up-face-geisking-2}
\end{center}
\end{figure}

\vspace{.15 in}\noindent {\bf -} Inflation at the face
$(p)(012)\leftrightarrow (p)(302)$; this face contains the two edges
$2$ and $3$ of the frame $\xi$. See Figure
\ref{f-blow-up-face-geisking-2}. In this example when we remove the
identification $(p)(012)\leftrightarrow(p)302)$, we have a choice of
order in adding the tetrahedra $(2)$ and $(3)$; this is the same as
adding a diagonal in the base of the pyramid we add to give a
triangulation of the pyramid. We make the selection of first adding
the tetrahedron $(3)$ along the face $(p)(012)$ and then adding the
tetrahedron $(2)$.

\noindent {\bf Step 7.} Inflation at the edges of $\T$. There is
only one edge, $E$. The inflation at an edge is determined by the
configuration polygon at that edge. In Step 4, Figure
\ref{f-config-polygon-labels-geisking}, we give the configuration
polygon and the transverse directions.

\begin{figure}[htbp]
\psfrag{C}{$(c)$}
           \psfrag{f}{\small $3^1$}\psfrag{a}{\small $\td{5}^0$}\psfrag{b}{\small $\td{5}^1$}
            \psfrag{e}{\small $3^0$}\psfrag{h}{\small $\td{3}^0$}
            \psfrag{c}{\small $2^0$}
            \psfrag{d}{\small $\td{2}^1$}
            \psfrag{B}{$(2)$}
            \psfrag{A}{$(3)$}
            \psfrag{E}{$(5)$}
            \psfrag{Y}{\large $(b_2^*)$} \psfrag{X}{\large $(b_1^*)$}
            \psfrag{0}{\footnotesize $0$}
\psfrag{1}{\footnotesize $1$}\psfrag{2}{\footnotesize
$2$}\psfrag{3}{\footnotesize $3$} \psfrag{i}{\begin{tabular}{l}
$(3)(012)\leftrightarrow(b_1^*)(210)$;\hspace{.25 in}$(2)(012)\leftrightarrow(b_1^*)(310)$\\
\\
$(5)(012)\leftrightarrow(b_2^*)(210)$;\hspace{.25
in}$(5)(013)\leftrightarrow(b_2^*)(320)$\\
\\
$(c)(012)\leftrightarrow(3)(031)$;\hspace{.275 in}$(c)(013)\leftrightarrow(b_2^*)(031)$\\
\\
$(c)(023)\leftrightarrow(b_1^*)(320)$;\hspace{.2 in}
$(c)(123)\leftrightarrow(2)(031)$\\ \end{tabular}}

        \vspace{0 in}
        \begin{center}
\includegraphics[width=6 in]{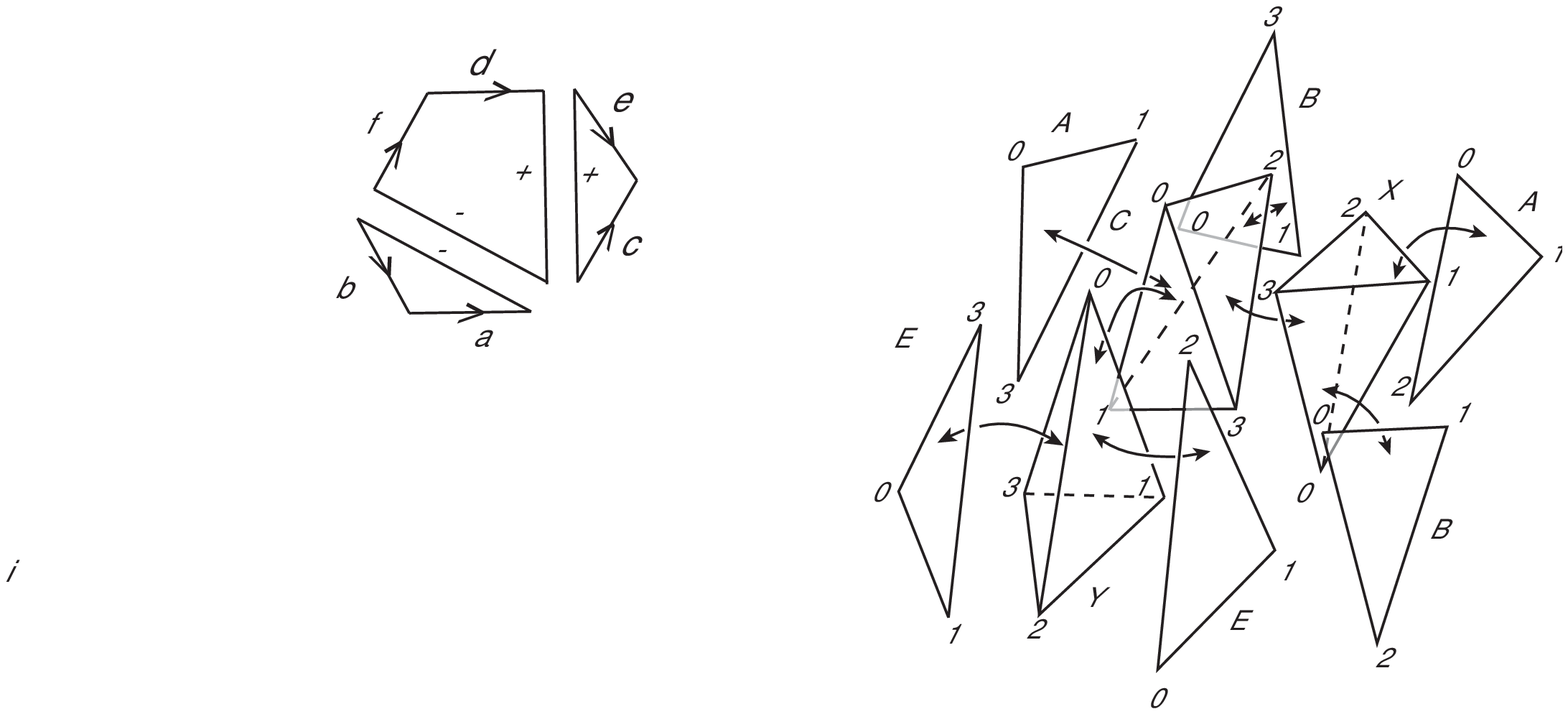}
\caption{Inflation at the edge $E$ where the configuration polygon decomposes into two branch triangles and a crossing.}
\label{f-blow-up-edge-geisking}
\end{center}
\end{figure}

\vspace{.15 in}\noindent {\bf -} Inflation at the edge $E$. The
configuration polygon decomposes into two branch configurations and
a crossing. See Figure \ref{f-blow-up-edge-geisking}.

\vspace{.25 in}\noindent {\bf Step 8.} Finish.

In this example, there are no face identification from $\T^*$
retained and there are no cones that need subdivided.

We collect the tetrahedra and face identifications into an array
(again, following  the notational conventions of
{\footnotesize\textsc REGINA} \cite{burton-regina}).

Here the two faces $(b_1^*)(123)$ and $(b_2^*)(123)$ are not
identified and become the boundary of the compact manifold given by
the inflation triangulation. The identifications of the edges of the
triangles  $(b_1^*)(123)$ and $(b_2^*)(123)$ are determined by the
face identifications and are:
$$(b_1^*)(12)\rightarrow (3)(10)\rightarrow (c)(20)\rightarrow(b_1^*)(23);$$
$$(b_1^*)(13)\rightarrow (2)(10)\rightarrow (c)(31)\rightarrow
 (b_2^*)(13);$$ and
$$(b_2^*)(12)\rightarrow (5)(10)\rightarrow (b_2^*)(23).$$
The identifications give a one-vertex (two-triangle) Klein bottle as
the boundary of the manifold underlying the inflation triangulation.

We collect the face identifications for the inflation triangulation
in the following array.

\vspace{.25 in}
\begin{tabbing}\hspace{.5 in}\=\bf {tet}\hspace{.55 in}\=$(012)$\hspace{.55 in}\=$(013)$
\hspace{.55 in}\=$(023)$\hspace{.55 in}\=$(123)$\\
\>$(p)$\>$(3)(321)$\>$(5)(231)$\>$(2)(032)$
\>$(5)(230)$\\
\>$(2)$\>$(b_1^*)(310)$\>$(c)(132)$\>$(p)(032)$\>$(3)(230)$\\
\>$(3)$\>$(b_1^*)(210)$\>$(c)(021)$\>$(2)(312)$
\>$(p)(210)$\\
\>$(5)$\>$(b_2^*)(210)$\>$(b_2^*)(320)$\>$(p)(312)$
\>$(p)(301)$\\
\>$(c)$\>$(3)(031)$\>$(b_2^*)(031)$\>$(b_1^*)(320)$
\>$(2)(031)$\\
\>$(b_1^*)$\>$(3)(210)$\>$(2)(210)$\>$(c)(320)$
\>{\bf bddry}\\
\>$(b_2^*)$\>$(5)(210)$\>$(c)(031)$\>$(5)(310)$
\>${\bf bddry}$\\
\end{tabbing}

The inflation of the one-tetrahedron Gieseking manifold gives a
$7$-tetrahedron triangulation of a compact, non-orientable
$3$--manifold with a normal boundary and interior homeomorphic to
the Gieseking manifold. The boundary is a Klein bottle.

\vspace{.75 in}


\centerline{\bf APPENDIX}

\begin{figure}[htbp]
\vspace{.25 in}
            \psfrag{p}{$(p)$}

            \psfrag{1}{\scriptsize{$1$}}\psfrag{2}{\scriptsize{$2$}}
            \psfrag{3}{\scriptsize{$3$}}\psfrag{4}{\scriptsize{$4$}}
            \psfrag{5}{\scriptsize{$5$}}\psfrag{6}{\scriptsize{$6$}}
            \psfrag{7}{\scriptsize{$7$}}\psfrag{8}{\scriptsize{$8$}}
            \psfrag{9}{\scriptsize{$9$}}\psfrag{x}{\scriptsize{$10$}}
            \psfrag{y}{\scriptsize{$11$}}\psfrag{z}{\scriptsize{$12$}}
            \psfrag{a}{\footnotesize $E$}
            \psfrag{b}{\footnotesize $F$}\psfrag{d}{\footnotesize $G$}

            \psfrag{A}{ $v^*$}\psfrag{B}{ $w^*$}
\psfrag{X}{\large\bf (A)}\psfrag{Y}{\large\bf
(B)}\psfrag{Z}{\large\bf (C)}
            \psfrag{D}{\begin{tabular}{c}
Diagonal\\Choice\\
\end{tabular}}\psfrag{V}{\begin{tabular}{c}
 vertex-linking\\torus at $v^*, V^*$\\
\end{tabular}}\psfrag{W}{\begin{tabular}{c}
vertex-linking\\torus at $w^*, W^*$\\
\end{tabular}}\psfrag{m}{\footnotesize $\mu_W$}\psfrag{l}{\footnotesize $\lambda_W = 7,9,\overline{3},\overline{1}$}
\psfrag{n}{\footnotesize $\mu_V$}\psfrag{k}{\footnotesize $\lambda_V
= 6,8,12,2 $}

\psfrag{C}{\Large \bf Example A3. Whitehead Link complement.}

        \vspace{0 in}
        \begin{center}
\includegraphics[width= 4.75  in]{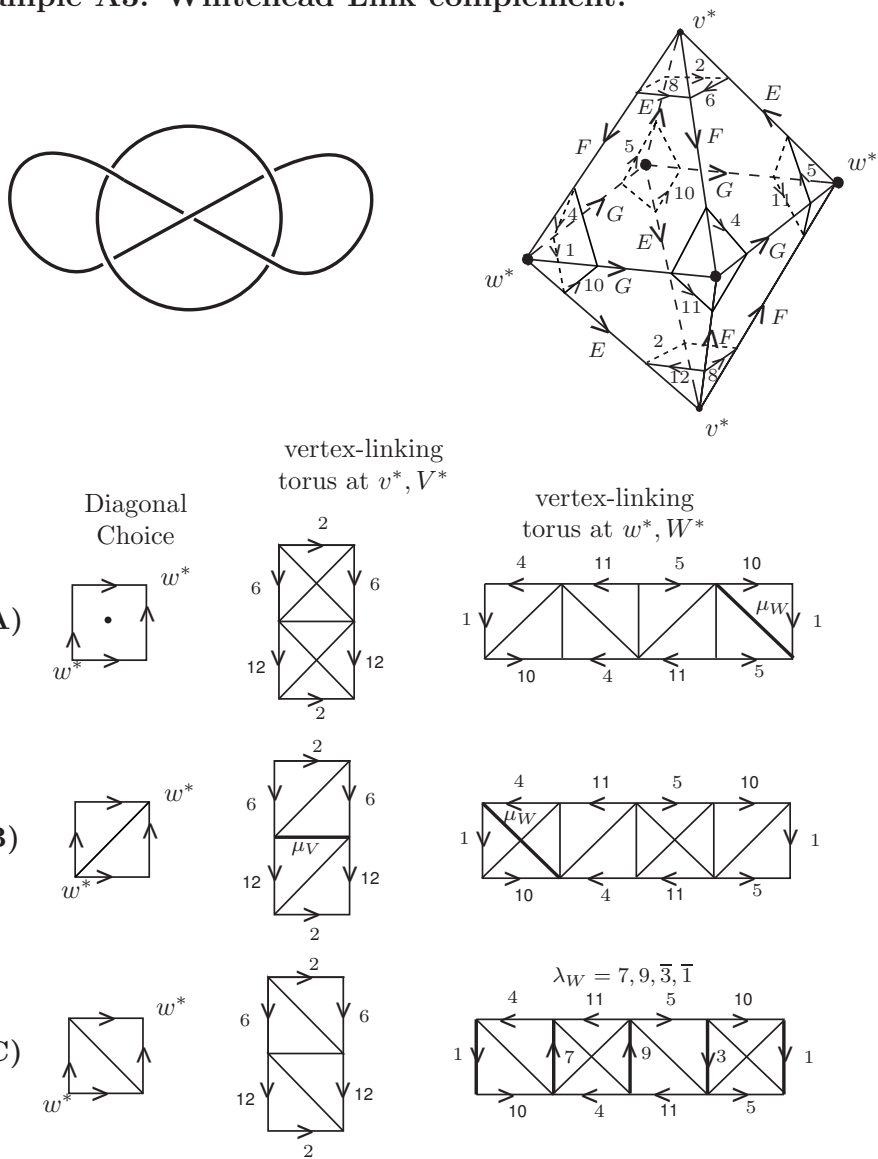} 
\caption{An ideal octagonal decomposition of the complement of the Whitehead Link in $S^3$.  Shown are the vertex-linking tori at the vertices $A$ and $B$, depending on the choice of diagonal in the octagon which subdivide it into a four-tetrahedra ideal triangulation of the Whitehead Link complement.} 
\label{f-ideal-whitehead}
\end{center}
\end{figure}

In Figure \ref{f-ideal-whitehead}, we give the complement of the
Whitehead Link in $S^3$ as an ideal cell-decomposition with just one
$3$-cell, an octahedron. There are two ideal vertices, one labeled
$v^*$ and the other $w^*$; hence, two vertex-linking surfaces, each
a torus and labeled $V$ and $W$, respectively. The vertex-linking
tori have induced cell-decompositions consisting of quadrilaterals.
As is well known, an octahedron can be decomposed into a
triangulation having four tetrahedra by choosing one of the three
possible diagonals. We consider the triangulations from each of
these choices. In $(A)$ the diagonal is between the vertices labeled
$v^*$ in the figure; in $(B)$ the diagonal is between the vertices
labeled $w^*$ in the figure and in $(C)$ the diagonal is between the
two unlabeled vertices in the figure, which are also identified with
 $w^*$.

In each case, we give the induced triangulation on the
vertex-linking torus. The meridian slope on $V$ is designated
$\mu_V$ and on $W$ it is designated $\mu_W$. In all subdivisions,
the meridian slope $\mu_V = \langle 2\rangle$ has length one;
however, the meridian slope $\mu_W$ has length $1$ in $(A)$, length
$2$ in $(B)$ and we can choose the meridian slope, $\mu_W$, to be
either length $1$ or $2$ in $(C)$. The longitudinal slopes,
$\lambda_V = \langle 6,2,8,12\rangle$ and $\lambda_W  = \langle
7,9,\overline{3},\overline{1}\rangle$, (considered as the longitude,
independently, in each component of the link) are circuits in all
the induced triangulations of the vertex-linking tori and in all
cases each has length $4$. The pair $\mu_W, \lambda_W$ forms a
frame; however, the pair $\mu_V, \lambda_V$ does not form a frame.
We can choose as a frame in the vertex-linking torus $V$ the pair
$\lambda'= \langle 6,12\rangle$ and $\mu_V = \langle 2\rangle$ (see
Figure \ref{f-frames-in-whitehead}).

\vspace{.5 in}
\address{Department of Mathematics, Oklahoma State University,
Stillwater, OK 74078}
\vspace{.10 in}

\email{jaco@math.okstate.edu}

\vspace{.25in}
\address{Department of Mathematics and Statistics,
University of Melbourne, Parkville, VIC 3052, Australia}

\vspace{.1 in}\email{rubin@maths.unimelb.edu.au}


\end{document}